\documentclass{amsart}
\usepackage{amssymb}
\usepackage{epsfig}
\usepackage{pstcol}             

\title{Permutohedra, associahedra, and beyond}

\author{Alexander Postnikov}

\address{Department of Mathematics, M.I.T., Cambridge, MA 02139}
\email{apost (at) math (dot) mit (dot) edu}
\def\mytilde{\kern-.015in\hbox{\lower.03in\hbox{\~{}}}\kern-.01in}
\urladdr{http://www.math-mit.edu/\mytilde apost/}


\thanks{The author was supported in part by National Science Foundation 
grant DMS-0201494 and by Alfred P.\ Sloan Foundation research fellowship}


\date{July 7, 2005}

\numberwithin{equation}{section}

\theoremstyle{plain}
\newtheorem{theorem}{Theorem}[section]
\newtheorem{proposition}[theorem]{Proposition}
\newtheorem{lemma}[theorem]{Lemma}
\newtheorem{corollary}[theorem]{Corollary}

\theoremstyle{definition}
\newtheorem{definition}[theorem]{Definition}

\newtheorem{example}[theorem]{Example}

\theoremstyle{remark}
\newtheorem{remark}[theorem]{Remark}

\def\R{\mathbb{R}}
\def\Z{\mathbb{Z}}
\def\C{\mathbb{C}}
\def\CC{\mathcal{C}}
\def\Q{\mathbb{Q}}

\def\ds{\displaystyle}
\def\Vol{\mathrm{Vol}\,}

\def\N{\mathcal{N}}
\def\D{\mathcal{D}}
\def\Todd{\mathrm{Todd}}
\def\J{\mathcal{J}}
\def\I{\mathcal{I}}
\def\diag{\mathrm{diag}}
\def\Ass{\mathrm{Ass}}
\def\Newton{\mathrm{Newton}}
\def\Waff{W_{\mathrm{aff}}}
\def\desc{\mathrm{desc}}
\def\flag{\mathit{flag}}

\begin{document}

\begin{abstract}
The volume and the number of lattice points of the permutohedron $P_n$ are
given by certain multivariate polynomials that have remarkable combinatorial
properties.  We give several different formulas for these polynomials.  We
also study a more general class of polytopes that includes the permutohedron,
the associahedron, the cyclohedron, the Pitman-Stanley polytope, and various
generalized associahedra related to wonderful compactifications of
De~Concini-Procesi.  These polytopes are constructed as Minkowski sums of
simplices.  We calculate their volumes and describe their combinatorial
structure.  The coefficients of monomials in $\Vol P_n$ are certain positive
integer numbers, which we call the mixed Eulerian numbers.  These numbers are
equal to the mixed volumes of hypersimplices.  Various specializations of
these numbers give the usual Eulerian numbers, the Catalan numbers, the
numbers $(n+1)^{n-1}$ of trees, the binomial coefficients, etc.  We calculate
the mixed Eulerian numbers using certain binary trees.  Many results are
extended to an arbitrary Weyl group.  \end{abstract}

\maketitle

\section{Introduction}

The {\it permutohedron\/} $P_n(x_1,\dots,x_{n})$ is the convex hull of the $n!$
points obtained from $(x_1,\dots,x_{n})$ by permutations of the coordinates.
Permutohedra appear in representation theory as {\it weight polytopes} of
irreducible representations of $GL_n$ and in geometry as {\it moment
polytopes}. 

In this paper we calculate volumes of permutohedra and numbers of their integer
lattice points.  Let us give a couple of examples.
It was known before that the volume of the
regular permutohedron $P_n(n,n-1,\dots,1)$ equals the number $n^{n-2}$ of
{\it trees\/} on $n$ labeled vertices and the number of lattice points of this
polytope equals
the number of {\it forests\/} on $n$ labeled vertices.  Another
example is the {\it hypersimplex\/} $\Delta_{k,n} = 
P_n(1,\dots,1,0,\dots,0,)$
(with $k$ ones). 
It is well-know that the volume of $\Delta_{k,n}$ is 
the Eulerian number, that is the number of permutations
of size $n-1$ with $k-1$ descents, divided by $(n-1)!$.
This calculation  dates back to Laplace~\cite{Lap}.  
These examples are just a tip of an iceberg.  They indicate 
at a rich combinatorial structure.  Both the volume and the number of
lattice points of the permutohedron $P_n(x_1,\dots,x_n)$ are given 
by multivariate polynomials in $x_1,\dots,x_n$ that have remarkable
properties.

We present three different combinatorial interpretations of these polynomials
using three different approaches.  Our first approach is based of Brion's
formula that expresses the sum of exponents over lattice points of a
polytope as a rational function.  From this we deduce a formula for the volume
of the permutohedron as a sum on $n!$ polynomials.  Then we deduce a
combinatorial formula for the coefficients in terms permutations
with given descent sets.  We extend the formula for the volume to 
weight polytopes for any Lie type.
There are some similarities between this formula and the Weyl's
character formula.

Our second approach is based on a way to represent permutohedra as a weighted
Minkowski sum $\sum y_I\Delta_I$ of the coordinate simplices.  We extend
our results to a larger class of polytopes that we call {\it generalized
permutohedra}.  These polytopes are obtained from usual permutohedra by
parallel translations of their faces.  

We discuss combinatorial structure of generalized permutohedra.
This class includes many interesting polytopes: associahedra,
cyclohedra, various generalized associahedra related to 
De~Concini-Procesi's wonderful compactifications,
graph associahedra, Pitman-Stanley polytopes, graphical zonotopes, etc.
We describe the combinatorial structure for a class
of generalized permutohedra in terms of {\it nested families.}
This description leads to a generalization of the 
Catalan numbers.

We calculate volumes of generalized
permutohedra by first calculating {\it mixed volumes\/} of various
coordinate simplices using {\it Bernstein's theorem\/} on systems of
algebraic equations.
More generally, we calculate the {\it Ehrhart polynomial\/} of generalized permutohedra,
i.e., the polynomial that expresses their number of lattice points.
Interestingly, the formula for the number of lattice points is obtained
from the formula for the volume by replacing usual powers in monomials
with raising powers.  We also found an interesting new {\it duality\/} 
for generalized permutohedra that preserves the number of lattice points.

We introduce and study {\it root polytopes\/} and their
triangulations.  These are convex hulls of the origin and end-points of 
several positive roots for a type $A$ root system. 
In particular, this class of polytopes includes direct products of
two simplices.  We apply the {\it Cayley trick\/} to show that the volume 
of a root polytope is related to the number of lattice points in a 
certain associated generalized permutohedron.  Each triangulation 
of a root polytope leads to a bijection between lattice points of the 
associated generalized permutohedron and its dual generalized permutohedron.

As an application of these techniques we solve a problem about
combinatorial description of diagonal vectors of shifted
Young tableaux of the triangular shape.

Our third approach is based on a way to represent permutohedra as a Minkowski
sum of the hypersimplices $\sum u_k\Delta_{k,n}$.  We express volumes of
permutohedra in terms of mixed volumes of the hypersimplices.
We call these mixed volume the {\it mixed Eulerian numbers.}
Various specializations of these numbers lead to the usual Eulerian
numbers, the Catalan numbers, the binomial coefficients,
the factorials,  the number $(n+1)^{n-1}$ of trees, and many other
combinatorial sequences.  We prove several identities for the
mixed Eulerian number and give their combinatorial interpretation
in terms of weighted binary trees.
We also extend this approach and generalize mixed Eulerian numbers
to an arbitrary root system. 

A brief overview of the paper follows.
In Section~\ref{sec:perm_and_zono}, we define permutohedra,
give their several known properties, and discuss their
relationship with zonotopes.
In Section~\ref{sec:descent_div_sym}, we give a formula
for volumes of permutohedra (Theorem~\ref{th:f1})
based on Brion's formula
and derive another formula for volumes 
(Theorem~\ref{th:vol=descent_number})
that involves numbers of permutations with given descents sets.
In Section~\ref{sec:weight-polytopes}, we give a formula
for volumes and lattice points enumerators of weight polytopes
for any Lie type (Theorems~\ref{th:f1-W}
and~\ref{th:sum-e-W}). 
In Section~\ref{sec:dragon_marriage_condition}, we give
a formula for volume of permutohedra 
(Theorem~\ref{th:second-formula})
based on our second approach.
In Section~\ref{sec:generalized_permutohedra}, we discuss
generalized permutohedra and several ways to parametrize
this class of polytopes.
In Section~\ref{sec:nested}, we discuss combinatorial
structure for a class of generalized permutohedra in terms
of nested families (Theorem~\ref{th:nested_complex}).
In Section~\ref{sec:examples_of_gen_perm}, we apply this
description to several special cases of generalized permutohedra.
In Section~\ref{sec:vol_via_Bernstein}, 
we extend Theorem~\ref{th:second-formula} to generalized
permutohedra and calculate their
volumes (Theorem~\ref{th:second-formula-generalized})
using Bernstein's theorem.
In Section~\ref{sec:vol_via_Brion}, we give alternative
formulas for volumes
(Theorems~\ref{th:gen-perm-sum-w} and~\ref{th:gen_descents_sum_w})
based on our first approach.
In Section~\ref{sec:generalized_Ehrhart}, we state
a formula for the Ehrhart polynomial of generalized permutohedra 
(Theorem~\ref{th:gen_ehrhrart})
and derive the duality theorem
(Corollary~\ref{cor:duality-lattice-points}).
In Section~\ref{sec:root_polytope}, we discuss root
polytopes and their triangulations for bipartite graphs.  
In Section~\ref{sec:root_non_bipartite}, we treat the case of 
non-bipartite graphs.
In Section~\ref{sec:subdivision}, we show how triangulations of roots 
polytopes are related to lattice points of generalized permutohedra.  
We also finish the proof of Theorem~\ref{th:gen_ehrhrart}.
In Section~\ref{sec:shifted_tableaux}, we describe diagonals
of shifted Young tableaux.
In Section~\ref{sec:Mixed_Eulerian}, we discuss our third
approach based on the mixed Eulerian numbers.  We prove
several properties of these numbers 
(Theorems~\ref{th:mixed_eul_properties} 
and~\ref{th:equivalence-classes}).
In Section~\ref{sec:weighted_binary_trees}, we give 
the third combinatorial formula for  volumes
of permutohedra (Theorem~\ref{th:vol_perm_binary_trees})
and give
a combinatorial interpretation for the mixed Eulerian numbers
(Theorem~\ref{th:mixed_eulerian_binary_trees}).
Finally, in Section~\ref{sec:vol_weight_Phi} we extend
our third approach to weight polytopes for an arbitrary
root system  (Theorems~\ref{th:Phi_volume_trees}
and~\ref{th:Phi_mixed_eulerian_binary_trees}).
In Appendix~\ref{sec:appendix-lattice-points}, we review
and give short proofs of needed general results on enumeration of 
lattice points in polytopes.

Let us give a notational remark about our use of various coordinate
systems.  We use the $x$-coordinates to parametrize permutohedra
expressed in the standard form as convex hulls of $S_n$-orbits
of $(x_1,\dots,x_n)$.  We use the $z$-coordinates to parametrize
(generalized) permutohedra expressed by linear inequalities as
$\{t \mid f_i(t)\geq z_i\}$, i.e., the $z$-coordinates correspond 
to the facets of these polytopes.
We use the $y$-coordinates to parametrize (generalized) permutohedra 
written as weighted Minkowski sums $\sum y_I \Delta_{I}$
of the coordinate simplices.
Finally, we use the $u$-coordinates to parametrize permutohedra written
as weighted Minkowski sums $\sum u_k\,\Delta_{n,k}$
of the hypersimplices. 
For all other purposes we use the $t$-coordinates.  
Throughout the paper, we use the notation 
$[n]:=\{1,2\dots,n\}$ and $[m,n]:=\{m,m+1,\dots,n\}$.

\medskip
\noindent
{\sc Acknowledgments:}
I thank Richard Stanley and Andrei Zelevinsky for helpful discussions.

\section{Permutohedra and zonotopes}
\label{sec:perm_and_zono}

\begin{definition}
For $x_1,\dots,x_{n}\in \R$,
the {\it permutohedron} $P_n(x_1,\dots,x_{n})$ 
is the convex polytope in $\R^{n}$ defined as the convex hull
of all vectors obtained from $(x_1,\dots,x_{n})$
by permutations of the coordinates:
$$
P_n(x_1,\dots,x_{n}):=\mathrm{ConvexHull}((x_{w(1)},\dots,x_{w(n)})
\mid w\in S_{n}),
$$
where $S_{n}$ is the symmetric group.
This polytope lies in the hyperplane 
$H_c=\{(t_1,\dots,t_{n})\mid t_1+\cdots + t_{n} = c\}
\subset \R^{n}$, where $c=x_1+\cdots+x_{n}$.
Thus $P_n(x_1,\dots,x_{n})$ has the dimension at most $n-1$.
\end{definition}

For example, for $n=3$ and distinct $x_1,x_2,x_3$, the permutohedron 
$P_3(x_1,x_2,x_3)$
is the hexagon shown below.
If some of the numbers $x_1,x_2, x_3$ are equal to each other
then the permutohedron degenerates into a triangle, or even a single
point.

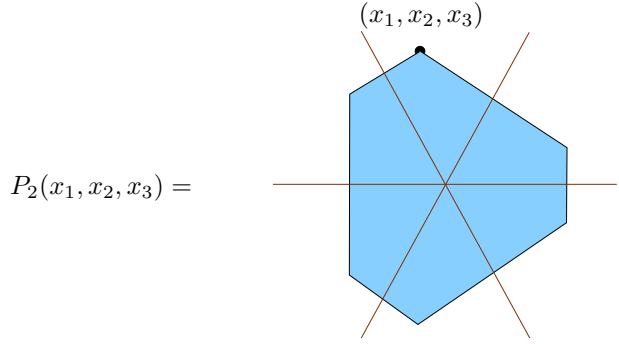
\begin{center}
\input{fig1-2.pstex_t}
\end{center}

For a polytope $P\in H_c$, define its volume $\Vol P$ as
the usual $(n-1)$-dimensional volume of the polytope $p(P)\in\R^{n-1}$, 
where $p$
is the projection $p:(t_1,\dots,t_{n})\mapsto(t_1,\dots,t_{n-1})$.
If $c\in\Z$, then the volume of any parallelepiped formed
by generators of the integer lattice $\mathbb{Z}^{n}\cap H_c$ is 1.

In this paper, we calculate the volume 
$$
V_n(x_1,\dots,x_{n}):=\Vol P_n(x_1,\dots,x_{n})
$$
of the permutohedron. Also, for integer $x_1,\dots,x_{n}$,
its number of lattice points 
$$
N_n(x_1,\dots,x_{n}):=P_n(x_1,\dots,x_{n})\cap \mathbb{Z}^{n}.
$$

We will see that both $V_n(x_1,\dots,x_{n})$ and $N_n(x_1,\dots,x_{n})$ 
are polynomials of degree $n-1$ in the variables $x_1,\dots,x_{n}$. 
The polynomial $V_n$ is the top homogeneous part of $N_n$.
The {\it Ehrhart polynomial} of the permutohedron 
is $E_{P_n}(t)=N_n(tx_1,\dots,tx_n)$.
We will give 3 totally different formulas for these polynomials.

\medskip

The special permutohedron for $(x_1,\dots,x_{n}) = (n-1,n-2,\dots,0)$,
$$
P_n(n-1,\dots,0)=\mathrm{ConvexHull}((w(1)-1,...,w(n)-1)\mid w\in 
S_{n})
$$
is the most symmetric permutohedron.  
It is invariant under the action of the symmetric group $S_{n}$.
For example, for $n=3$, it is the regular hexagon:
\begin{center}
\input{fig2-1.pstex_t}
\end{center}

We will call this special permutohedron $P_n(n-1,\dots,0)$
the {\it regular permutohedron}.
The volume of the regular permutohedron and its Ehrhart polynomial 
can be easily calculated using the general result on graphical zonotopes 
given below.

Recall that the {\it Minkowski sum} of several subsets $A$, \dots, $B$ 
in a linear space is the locus of sums of vectors that belong
to these subsets $A+\cdots + B :=\{a+\cdots + b\mid
a\in A, \dots, b\in B\}$.  If $A,\dots,B$ are convex polytopes then
so is their Minkowski sum.  The Newton polytope $\Newton(f)$
for a polynomial $f=\sum_{a\in\Z^n} \beta_{a}\,t_1^{a_1}\cdots t_n^{a_n}$
is the convex hull of integer points $a\in\Z^n$ such that $\beta_a\ne 0$.
Then $\Newton(f\cdot g)$ is the Minkowski sum $\Newton(f)+\Newton(g)$.
A {\it zonotope} is a Minkowski sum of several line intervals.

\begin{definition}
\label{def:graph_zonotope}
For a graph $\Gamma$ on the vertex set $[n]:=\{1,\dots,n\}$,
the {\it graphical zonotope\/} $Z_\Gamma$ is defined as the Minkowski
sum of the line intervals:
$$
Z_\Gamma :=\sum_{(i,j)\in \Gamma} [e_i,e_j] = \Newton\left(\prod_{(i,j)\in \Gamma}
(t_i-t_j)\right),
$$
where the Minkowski sum and the product are over edges $(i,j)$, $i<j$, 
of the graph $\Gamma$, and
$e_1,\dots,e_{n}$ are the coordinate vectors in $\R^{n}$.  The zonotope
$Z_\Gamma$ lies in the hyperplane $H_{c}$, where $c$ is the number of edges of $\Gamma$.
The polytope $Z_\Gamma$ was first introduced by Zaslavsky (unpublished).
\end{definition}

The following two claims express well-know properties of graphical
zonotopes and permutohedra.

\begin{proposition}
The regular permutohedron $P_n(n-1,\dots,0)$ is the 
graphical zonotope $Z_{K_{n}}$ for the complete graph $K_{n}$.
\end{proposition}

\begin{proof}  
The permutohedron $P_n(n-1,\dots,0)$ is the Newton polytope 
of the Vandermonde determinant 
$\det( t_i^{j-1})_{1\leq i,j\leq n}$.
On the other hand, the Vandermonde determinant
is equal to the product $\prod_{1\leq i<j\leq n} (t_j-t_i)$, 
whose Newton polytope is the zonotope $Z_{K_{n}}$.
\end{proof}

The following claim is given in Stanley~\cite[Exer.~4.32]{EC1}.

\begin{proposition}  
\label{prop:vol-Z-G}
For a connected graph $\Gamma$ on $n$ vertices, the volume 
$\Vol\, Z_\Gamma$ of the graphical zonotope $Z_\Gamma$
equals the number of spanning trees of the graph $\Gamma$.
The number of lattice points of $Z_\Gamma$ equals to the number of
forests in the graph $\Gamma$.

In particular, the volume of the regular permutohedron
is $\Vol P_n(n-1,\dots,0) = n^{n-2}$ and its number of lattice points 
equals the number of forests on $n$ labeled vertices.
\end{proposition}

The zonotope $Z_\Gamma$ can be subdivided into unit parallelepipeds associated
with spanning trees of $\Gamma$, which implies the first claim.

In general, for arbitrary $x_1,\dots,x_{n}$, the permutohedron
$P_n(x_1,\dots,x_{n})$ is not a zonotope.  We cannot easily calculate its
volume by subdividing it into parallelepipeds.

One can alternatively describe the permutohedron 
$P_n(x_1,\dots,x_{n})$ in terms of linear inequalities.

\begin{proposition} 
\label{prop:Rado}
{\rm Rado~\cite{Rad}}
Let us assume that $x_1\geq \cdots \geq x_{n}$.
Then a point $(t_1,\dots,t_{n})\in \R^n$ belongs to the permutohedron 
$P_n(x_1,\dots,x_{n})$ if and only if
$$
t_1+\dots+t_{n} = x_1 + \dots + x_{n}
$$
and, for any nonempty subset $\{i_1,\dots,i_k\}\subset\{1,\dots,n\}$,
we have
$$
t_{i_1}+\cdots +t_{i_k} \leq x_1+\cdots + x_k.
$$
\end{proposition}

The combinatorial structure of the permutohedron $P_n(x_1,\dots,x_{n})$
does not depend on $x_1,\dots,x_{n}$ as long as all these numbers
are distinct.  More precisely, we have the following well-know statement.

\begin{proposition}
\label{prop:comb_structute_perm}
Let us assume that $x_1 > \cdots > x_{n}$.
The $d$-dimensional faces of $P_n(x_1,\dots,x_{n})$ are
in one-to-one correspondence with disjoint subdivisions of the 
set $\{1,\dots,n\}$ into nonempty ordered blocks
$B_1 \cup \dots \cup B_{n-d} = \{1,\dots,n\}$.
The face corresponding to the subdivision into blocks $B_1,\dots,B_{n-d}$
is given by the $n-d$ linear equations
$$
\sum_{i\in B_1\cup \cdots \cup B_k} t_i = x_1 + \dots + 
x_{|B_1\cup\cdots \cup B_k|},\quad\textrm{for } k =1,\dots,n-d.
$$
In particular, two vertices $(x_{u(1)},\dots,x_{u(n)})$
and $(x_{w(1)},\dots,x_{w(n)})$, $u,w\in S_{n+1}$, are connected
by an edge if and only if $w=u\,s_i$, for some adjacent transposition
$s_i = (i,i+1)$.
\end{proposition}

\section{Descents and divided symmetrization}
\label{sec:descent_div_sym}

\begin{theorem} 
\label{th:f1}
Let us fix distinct numbers $\lambda_1,\dots,\lambda_{n}\in \R$.
The volume of the 
permutohedron $P_n = P_n(x_1,\dots,x_{n})$ is equal to
$$
\Vol P_n  = \frac{1}{(n-1)!} \sum_{w\in S_{n}}
\frac{(\lambda_{w(1)} x_1 + \cdots + \lambda_{w(n)} x_{n})^{n-1}}
{(\lambda_{w(1)}-\lambda_{w(2)})
(\lambda_{w(2)}-\lambda_{w(3)})\cdots
(\lambda_{w(n-1)}-\lambda_{w(n)}). }
$$
\end{theorem}

Notice that all $\lambda_i$'s in the right-hand side cancel each other after 
the symmetrization.  Theorem~\ref{th:f1-W} below gives a similar formula 
for any Weyl group.  Its proof is based on the 
Brion's formula~\cite{Bri};
see Appendix~\ref{sec:appendix-lattice-points}.

Theorem~\ref{th:f1} gives an efficient way to calculate the polynomials 
$V_n = \Vol P_n$.
However this theorem does not explain the combinatorial significance of the
coefficients in these polynomials.  The next theorem gives a combinatorial
interpretation for the coefficients.

Given a sequence of nonnegative integers $(c_1,\dots,c_{n})$
such that $c_1+\cdots+c_{n} = n-1$, let us construct the sequence
$(\epsilon_1,\dots,\epsilon_{2n-2})\in \{1,-1\}^{2n-2}$ by replacing
each entry `$c_i$' with `$1,\dots,1,-1$' ($c_i$ `1's followed by one `$-1$'), 
for $i=1,\dots,n$, and then removing the last `$-1$'.
For example, the sequence $(2,0,1,1,0,1)$ gives 
$(1,1,-1,-1,1,-1,1,-1,-1,1)$. This map is actually a bijection
between the sets
$\{(c_1,\dots,c_{n})\in\Z_{\geq 0}^{n}\mid c_1+\cdots +c_{n} =n-1\}$
and $\{(\epsilon_1,\dots,\epsilon_{2n-2})\in\{1,-1\}^{2n-2}\mid 
\epsilon_1+\cdots + \epsilon_{2n-2} =0\}$.
Let us define the set $I_{c_1,\dots,c_{n}}$ by
$$
I_{c_1,\dots,c_{n}}:=\{i\in \{1,\dots,n-1\}\mid 
\epsilon_1 + \cdots + \epsilon_{2i-1}<0\}.
$$

The {\it descent set} of a permutation $w\in S_{n}$ 
is $I(w)=\{i\in\{1,\dots,n-1\} \mid w_i> w_{i+1}\}$.
Let $D_{n}(I)$ be the number of permutations in $S_{n}$
with the descent set $I(w)=I$.

\begin{theorem} 
\label{th:vol=descent_number}
The volume of the permutohedron $P_n = P_n(x_1,\dots,x_{n})$
is equal to
$$
\Vol P_n = 
\sum 
(-1)^{|I_{c_1,\dots,c_{n}}|}
\, D_{n}(I_{c_1,\dots,c_{n}})\,
\frac{x_1^{c_1}}{c_1!}\cdots \frac{x_{n}^{c_{n}}}{c_{n}!}\,,
$$
where the sum is over sequences of nonnegative integers $c_1,\dots,c_{n}$
such that $c_1+\cdots +c_{n} =n-1$.
\end{theorem}
  
We can graphically describe the set $I_{c_1,\dots,c_{n}}$, as follows.
Let us construct the lattice path $P$ on $\Z^2$
from $(0,0)$ to $(n-1,n-1)$ with steps of the two types $(0,1)$ ``up'' and 
$(1,0)$ ``right'' such that $P$ has exactly $c_i$ up steps 
in the $(i-1)$-st column, for $i=1,\dots,n$.
Notice that the $(2i-1)$-th and $2i$-th steps in the path $P$ are either
both above the $x=y$ axis or both below it.  The set $I_{c_1,\dots,c_{n}}$ 
is the set of indices $i$ such that the $(2i-1)$-th and $2i$-th steps in $P$ are
below the $x=y$ axis.
\medskip

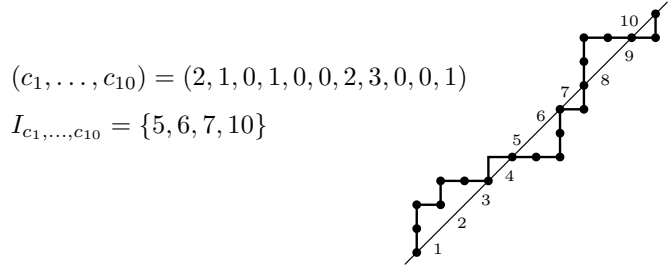
\begin{center}
\input{fig9-1.pstex_t}
\end{center}

\begin{example}
We have
$V_2= x_1 - x_2$ and
$V_3=\frac{x_1^2}{2} + x_1 x_2 - 2 x_1 x_3 - 2 \frac {x_2^2}{2}
+x_2 x_3 + \frac{x_3^2}{2}$.
The following figure shows the paths corresponding to all terms
in $V_2$ and $V_3$.
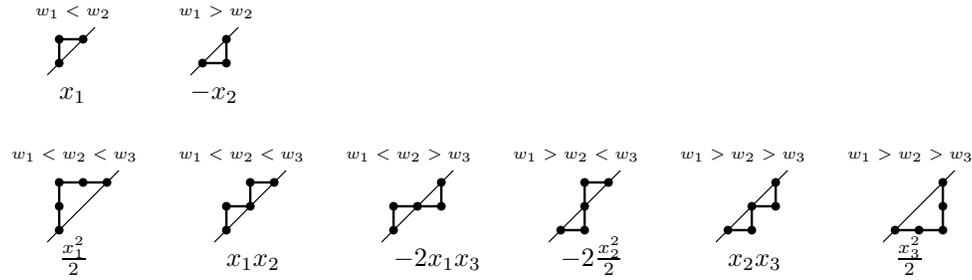
\begin{center}
\input{fig10-1.pstex_t}
\end{center}
For example, $I_{1,0,1} = \{2\}$ and there are 2 permutations $132, 231\in S_3$
with the descent set $\{2\}$.  Thus the coefficient of $x_1 x_3$ in $V_3$
is $-2$.  
\end{example}
\medskip

For a polynomial $f(\lambda_1,\dots,\lambda_{n})$, define its
{\it divided symmetrization} by
$$
\left<f\right>:=
\sum_{w\in S_{n}} w\left(\frac{f(\lambda_1,\dots,\lambda_{n})}
{(\lambda_1-\lambda_2)(\lambda_2-\lambda_3)\cdots (\lambda_{n-1}-\lambda_{n})}
\right),
$$
where the symmetric group $S_{n}$ acts by permuting the variables $\lambda_i$.

\begin{proposition}
\label{prop:<f>=const}
Let $f$ be a polynomial of degree $n-1$ in the variables 
$\lambda_1,\dots,\lambda_{n}$.  Then its divided symmetrization
$\left<f\right>$ is a constant.
If $\deg f<n-1$, then $\left<f\right>=0$.
\end{proposition}

\begin{proof}
We can write $\left<f\right> = g/\Delta$, where 
$\Delta=\prod_{i<j}(\lambda_i-\lambda_j)$ is the common denominator
of all terms in $\left<f\right>$ and $g$ is a certain polynomial of degree 
$\deg \Delta = \binom{n}{2}$.
Since $\left<f\right>$ is a symmetric rational function, $g$ should
be an anti-symmetric polynomial and thus it is divisible by $\Delta$.
Since $g$ and $\Delta$ have the same degree, their quotient is a constant.
If $\deg f<n-1$, then $\deg g< \deg \Delta$ and, thus, $g=0$.
\end{proof}

\begin{proposition}  
\label{prop:<lambda>}
We have
$\left<\lambda_1^{c_1}\cdots \lambda_{n}^{c_{n}}\right>=
(-1)^{|I|} D_{n}(I)$,
where $c_1,\dots,c_{n}$ are nonnegative integers 
with $c_1+\cdots+c_{n} = n-1$ and $I= I_{c_1,\dots,c_{n}}$.
\end{proposition}

\begin{proof}
We can expand the expression $\frac{1}{\lambda_i - \lambda_j}$, $i<j$
as the Laurent series that converges in the region 
$\lambda_1>\cdots > \lambda_{n} > 0$: 
$$
\frac{1}{\lambda_i - \lambda_j} = \lambda_i^{-1} \frac {1}{1-\lambda_j/\lambda_i}
= \sum_{k\geq 0} \lambda_i^{-k-1}\, \lambda_j^k.
$$
Let us use this formula to expand each term 
$w\left(\frac{\lambda_1^{c_1}\cdots \lambda_{n}^{c_{n}}}{(\lambda_1-\lambda_2)
\cdots (\lambda_{n-1}-\lambda_{n})}\right)$
as a Laurent series $f_w$ that converges in this region. 
Let $CT_w$ be the constant term of the series $f_w$.
Then, according to Proposition~\ref{prop:<f>=const}, we have
$\left<\lambda_1^{c_1}\cdots \lambda_{n}^{c_{n}}\right> = 
\sum_{w\in S_{n}} CT_w$.
Equivalently, the number $CT_w$ is the constant term in
the series $w^{-1} (f_w)$, i.e., 
the Laurent series obtained by the expansion of each term 
$\frac{1}{\lambda_i - \lambda_{i+1}}$ in
$\frac{\lambda_1^{c_1}\cdots \lambda_{n}^{c_{n}}}{(\lambda_1-\lambda_2)
\cdots (\lambda_{n-1}-\lambda_{n})}$ as
$$
\frac{1}{\lambda_i - \lambda_{i+1}} = \left\{
\begin{array}{cl}
\ds
\sum_{k\geq 0} \lambda_i^{-k-1}\, \lambda_{i+1}^k, & \textrm{for }w(i)<w(i+1), \\[.2in]
\ds
-\sum_{k\geq 0} \lambda_i^{k}\, \lambda_{i+1}^{-k-1}, & \textrm{for }w(i)>w(i+1).
\end{array}
\right.
$$
Let $I=I(w)$ be the descent set of the permutation $w$.
Then $CT_w$ equals $(-1)^{|I|}$ times the number of 
nonnegative integer sequences $(k_1,\dots,k_{n-1})$ 
such that we have $(c_1,\dots,c_{n}) = v_1+\dots + v_{n-1}$,
where
$$
v_i = \left\{ \begin{array}{cl} 
(k_i + 1)\,e_i - k_i \, e_{i+1}, & \textrm{for }i \not\in I  
\ \ (w_i<w_{i+1}), \\[.1in] 
-k_i\,e_i + (k_i +1) \, e_{i+1}, & \textrm{for }i \in I 
\ \ (w_i>w_{i+1}),
\end{array}
\right.
$$
and the $e_i$ are the coordinate vectors. 
Notice that, for a fixed permutation $w$, there is at most 1 sequence 
$(k_1,\dots,k_{n-1})$
that produces $(c_1,\dots,c_{n})$, as above.  Thus $CT_w\in\{1,-1,0\}$.

Let $P$ be the lattice path from $(1,1)$ to $(n,n)$ constructed from 
the sequence $(c_1,\dots,c_{n})$ as shown after
Theorem~\ref{th:vol=descent_number}.
In other words,  $P$ is the continuous piecewise-linear path obtained by joining
the points 
$$
(0,0) - (0,c_1) - (1,c_1) - (1, c_1 + c_2) - (2, c_1 + c_2) - 
(2,c_1+c_2+ c_3) - \cdots - (n-1,n-1)
$$
by the straight lines.

Let $r$ be the maximal index such that $w(1)< w(2)<\cdots < w(r)$.
Then we have $c_1 = k_1+1$, $c_2 = k_2 + 1 - k_1$, \dots, $c_{r-1} = k_{r-1} + 1 
- k_{r-2}$, $c_r = - k_r - k_{r-1}$.
Thus $k_i = c_1 +\cdots + c_i -i\geq 0 $, for $i = 1,\dots,r-2$,
$k_{r-1} = c_1+ \cdots + c_{r-1} - (r-1) = 0$ and $k_r = c_r = 0$.
This means that the path $P$ stays weakly above 
the $x=y$ axis as it goes from the point $(0,0)$ to the point $(r-1,r-1)$,
then it passes through the point $(r-1,r-1)$, and goes
strictly below the $x=y$ axis (if $r< n+1$).
For $i=1,\dots, r-1$, the number $k_i$ is exactly the distance between
the lowest point of the path $P$ on the line $x=i$ and the point $(i,i)$.

Let $r'$ be the maximal index such $w(r)>w(r+1)>\cdots > w(r')$.
Then we have  $c_r = - k_r=0$, $c_{r+1} = k_r + 1 - k_{r+1}$, \dots,
$c_{r'-1} = k_{r'-2} + 1 - k_{r'-1}$, and 
$c_{r'} = (k_{r'-1} + 1) + (k_{r} + 1)$.   Thus
$k_i = i-r - c_r - \cdots - c_i = i-1 - c_1-\cdots - c_i\geq 0$, 
for  $i=r,\dots,r'-1$,
and $k_{r'} = c_r + \cdots + c_{r'} -r \geq 0$.
This means that the path $P$ stays weakly below the $x=y$ axis as it goes
from the point $(r-1,r-1)$ to the point $(r'-1,r'-1)$,
then it passes through the point $(r'-1,r'-1)$ and goes strictly above the 
$x=y$ axis (if $r'<n+1$).
For $i=r,\dots, r'-1$, the number $k_i$ is the distance between 
the highest point of the path $P$ on the line $x=i-1$ and the point $(i-1,i-1)$.

We can continue working with maximal monotone intervals in the permutation
$w$ in this fashion.
Let $r''$ be the maximal index such that $w(r')< \cdots < w(r'')$.
Similarly to the above argument, we obtain that that path $P'$ stays 
weakly above the $x=y$ axis until it crosses it at the point $(r''-1,r''-1)$,
etc.  

We deduce that the indices $r,r',r'',\dots$ characterizing 
the descent set of $w$ correspond to the points where the path $P$ crosses
the $x=y$ axis.  Thus the descent set of $w$ is uniquely reconstructed
from the sequence $(c_1,\dots,c_n)$ as $I= I_{c_1,\dots,c_{n}}$.
Moreover, for any permutation $w$ with such descent set, 
the nonnegative integer sequence $(k_1,\dots,k_{n-1})$ 
is uniquely reconstructed
from the sequence $(c_1,\dots,c_n)$ as 
$$
k_i = \left\{
\begin{array}{ll}
\min \{y-i\mid  (i,y)\in P\} & \textrm{if } i\not\in I,\\
\min \{i-1-y\mid  (i-1,y)\in P\} & \textrm{if } i\in I,
\end{array}
\right.
$$
and, thus, $CT_w = (-1)^{|I|}$.
This shows that only permutations with the descent set $I=I_{c_1,\dots,c_{n}}$ 
make a contribution
to $\left<\lambda_1^{c_1}\cdots \lambda_{n}^{c_{n}}\right>$, and the
contribution of any such permutation is $(-1)^{|I|}$.
This finishes the proof.
\end{proof}

\begin{proof}[Proof of Theorem~\ref{th:vol=descent_number}]
According to Theorem~\ref{th:f1}, the volume of the permutohedron 
can be written as the divided symmetrization of the power of 
a linear form:
$$
V_n = \frac {1}{(n-1)!} 
\left<(x_1\lambda_1 + \cdots + x_{n}\lambda_{n})^{n-1}\right>
=\sum_{c_1+\cdots + c_{n} = n-1}  
\left<\lambda_1^{c_1}\cdots \lambda_{n}^{c_n}\right>\,
\frac{x_1^{c_1}}{c_1!}\cdots \frac{x_{n}^{c_{n}}}{c_{n}!}\,.
$$
Now apply Proposition~\ref{prop:<lambda>}.
\end{proof}

\section{Weight polytopes}
\label{sec:weight-polytopes}

Theorem~\ref{th:f1} can be extended to any Weyl group, as follows.
Let $\Phi$ be a root system of rank $r$.  Let $\Lambda$ be the 
associated integer {\it weight lattice} and $\Lambda_\R = \Lambda\otimes\R$ 
be the weight space.  The roots in $\Phi$ span the root lattice
$L\subseteq \Lambda$.  
The associated {\it Weyl group} $W$ 
acts on the weight space $\Lambda_\R$.
Let $(x,y)$ be a nondegenerate $W$-invariant inner product on $\Lambda_\R$.

\begin{definition}
\label{def:weight_polytope}
For $x\in \Lambda_\R$, we can define the {\it weight polytope} 
$P_W(x)$ as the convex hull of a Weyl group orbit:
$$
P_W(x) :=\mathrm{ConvexHull}(w(x)\mid w\in W)\subset \Lambda_\R.
$$
For the Lie type $A_r$, the weight polytope 
$P_W(x)$ is the permutohedron $P_{r+1}(x)$.
\end{definition}

Let us fix a choice of {\it simple roots} $\alpha_1,\dots,\alpha_r$ in $\Phi$. 
Let $\Vol$ be the volume form on $\Lambda_\R$ normalized so that
the volume of the parallelepiped generated by the simple roots 
$\alpha_i$ is 1.  
Recall that a weight $\lambda\in\Lambda_\R$ is called 
{\it regular} if $(\lambda,\alpha)\ne 0$ for any root $\alpha\in\Phi$.
A weight $\lambda$ is called {\it dominant} if $(\lambda,\alpha_i)\geq 0$, for 
$i=1,\dots,r$.

\begin{theorem}  
\label{th:f1-W}
Let $\lambda\in \Lambda_\R$ be a regular weight.
The volume of the weight polytope is equal to
$$
\Vol P_W(x) = \frac{1}{r!}\sum_{w\in W} \frac{(\lambda,w(x))^r}
{(\lambda,w(\alpha_1))\cdots (\lambda,w(\alpha_r))}.
$$
\end{theorem}

For type $A_r$, $W=S_{r+1}$ and 
Theorem~\ref{th:f1-W} specializes to Theorem~\ref{th:f1}.

Let $G$ be a Lie group with the root system $\Phi$.
For a dominant weight $\lambda$, let $V_\lambda$ be the irreducible 
representation of $G$ with the highest weight $\lambda$.
The character of $V_\lambda$ is a certain nonnegative linear combination
$ch(V_\lambda)$ of the formal exponents $e^\mu$, $\mu\in\Lambda$. 
(These formal exponents are subject to the relation 
$e^\mu\cdot e^\nu=e^{\mu+\nu}$.)
The weights that occur in the representation $V_\lambda$ with nonzero
multiplicities, i.e., the weights $\mu$ such that $e^\mu$ has a nonzero
coefficient in $ch(V_\lambda)$, are exactly the points of the weight polytope
$P_W(\lambda)$ in the lattice $L+\lambda$ (the root lattice shifted by
$\lambda$).
Let 
$$
S(P_W(\lambda)):=\sum_{\mu\in P_W(\lambda)\cap(L+\lambda)}e^\mu
$$ 
be the sum of formal exponents over these lattice points.
In other words, $S(P_W(\lambda))$ is obtained from the character
$ch(V_\lambda)$ by replacing all nonzero coefficients with 1.
For example, in the type $A$, the expression $S(P_n(\lambda))$
is obtained from the Schur polynomial by erasing the coefficients
of all monomials.

We have the following identity in the field of rational expressions
in the formal exponents.

\begin{theorem} 
\label{th:sum-e-W}
For a dominant weight $\lambda$, the sum of exponents over lattice 
points of the weight polytope $P_W(\lambda)$ equals
$$
S(P_W(\lambda)) = \sum_{w\in W} 
\frac{e^{w(\lambda)}}
{(1-e^{-w(\alpha_1)})\cdots (1-e^{-w(\alpha_r)})}.
$$
\end{theorem}

Notice that if we replace the product over simple roots  $\alpha_i$
in the right-hand side of Theorem~\ref{th:sum-e-W} by 
a similar product over {\it all} positive roots, we obtain exactly
Weyl's character formula for $ch(V_\lambda)$.

Theorems~\ref{th:f1}, \ref{th:f1-W}, and \ref{th:sum-e-W} follow from 
Brion's formula~\cite{Bri} on summation over lattice points in a rational 
polytope.
In Appendix~\ref{sec:appendix-lattice-points}, we give a brief overview
of this result and related results of Khovanskii-Pukhlikov~\cite{KP1, KP2}
and Brion-Vergne~\cite{BV1, BV2}.  The following proof assumes reader's
familiarity with the Appendix.

\begin{proof}[Proof of Theorems~\ref{th:f1}, \ref{th:f1-W}, \ref{th:sum-e-W}]
Let us identify the lattice $L+\lambda$ embedded into $\Lambda_\R$
with $\Z^r\subset \R^r$.
Then (for a regular weight $\lambda$) the polytope $P_W(\lambda)$
is a Delzant polytope, i.e., for any vertex of $P_W(\lambda)$,
the cone at this vertex is generated by an integer basis of the lattice $\Z^r$;
see Appendix~\ref{sec:appendix-lattice-points}.
Indeed, the generators of the cone at the vertex $\lambda$
are $-\alpha_1,\dots,-\alpha_r$.  Thus the generators
of the cone at the vertex $w(\lambda)$, for $w\in W$, 
are $g_{i, w(\lambda)} = -w(\alpha_i)$, $i=1,\dots,r$.
Now Theorem~\ref{th:sum-e-W} is obtained from 
Brion's formula given in Theorem~\ref{th:S-any-polytope}(2).
As we mention in the proof of 
Theorem~\ref{th:Todd-Euler-Maclaurin}(1), this claim remains
true for non-regular weights $\lambda$ when some of the vertices 
$w(\lambda)$ may accidentally merge.
Similarly, Theorems~\ref{th:f1} and \ref{th:f1-W},
are obtained from Theorem~\ref{th:S-any-polytope}(4).
\end{proof}

In a sense, Theorems~\ref{th:f1-W} and \ref{th:f1}
are deduced from Theorem~\ref{th:sum-e-W} in the same way
as Weyl's dimension formula is deduced from Weyl's character
formula, cf.~Appendix~\ref{sec:appendix-lattice-points}.

\section{Dragon marriage condition}
\label{sec:dragon_marriage_condition}

In this section we give a different combinatorial formula
for the volume of the permutohedron.

Let us use the coordinates $y_1,\dots,y_{n}$ related to
$x_1,\dots,x_{n}$ by
$$
\left\{
\begin{array}{l}
y_1 = -x_1\\
y_2 = - x_2 + x_1\\
y_3 = - x_3 + 2 x_2 - x_1 \\
\vdots
\\
y_{n} = - \binom{n-1}{0}\, x_n + \binom{n-1}{1}\, x_{n-1} - \cdots \pm 
\binom{n-1}{n-1}\,x_1
\end{array}
\right.
$$
Write $V_n = \Vol P_n(x_1,\dots,x_{n})$ as a polynomial
in the variables $y_1,\dots,y_{n}$.

\begin{theorem}   
\label{th:second-formula}
We have
$$
\Vol P_n = \frac{1}{(n-1)!}\sum_{(J_1,\dots,J_{n-1})}
y_{|J_1|}\cdots y_{|J_{n-1}|},
$$
where the sum is over ordered collections of subsets 
$J_1,\dots,J_{n-1}\subseteq [n]$
such that, for any distinct $i_1,\dots,i_k$, we have
$|J_{i_1}\cup \cdots \cup J_{i_k}| \geq k+1$. 
\end{theorem}

We will extend and prove Theorem~\ref{th:second-formula} for a larger class of
polytopes called generalized permutohedra; see
Theorem~\ref{th:second-formula-generalized}.
Theorem~\ref{th:second-formula} 
implies that $(n-1)!\,V_n$ is a polynomial in $y_2,\dots,y_{n}$ 
with {\it positive} integer coefficients.

\begin{example}
We have
$V_2 = \Vol([(x_1,x_2),(x_2,x_1)]) = x_1 - x_2 = y_2$
and $2V_3 =
x_1^2 + 2x_1 x_2 - 4 x_1 x_3 - 2 x_2^2 + 2x_2 x_3 + x_3^2
= 6\,y_2^2 + 6\, y_2\, y_3 +  y_3^2$.
\end{example}

\begin{remark}  The condition on subsets $J_1,\dots,J_{n-1}$
in Theorem~\ref{th:second-formula} is similar to the condition in 
Hall's marriage theorem~\cite{Hal}.  One just needs to replace the inequality 
$\geq k+1$ with $\geq k$ to obtain Hall's marriage condition.
\end{remark}

Let us give an analogue of the marriage problem and Hall's theorem.  

\medskip
\noindent
{\large$\mathfrak{Dragon\ marriage\ problem.}$}  
{\it
There are $n$ brides, $n-1$ grooms living in a medieval town, 
and 1 dragon who likes to visit the town occasionally.
Suppose we know all possible pairs of brides and grooms
who do no mind to marry each other.
A dragon comes to the village and takes one of the brides.
When will it be possible to match the remaining brides and grooms
no matter what the choice of the dragon was?
}

\begin{proposition}
\label{prop:dragon-marriage}
Let $J_1,\dots,J_{n-1}\subseteq[n]$.
The following three conditions are equivalent:
\begin{enumerate}
\item 
For any distinct $i_1,\dots,i_k$, we have
$|J_{i_1}\cup \cdots \cup J_{i_k}| \geq k+1$. 
\item For any $j\in[n]$, there is a 
system of distinct representatives
in $J_1,\dots,J_{n-1}$ that avoids $j$.
(This is a reformulation of the dragon marriage problem.)
\item  There is a system of $2$-element representatives
$\{a_i,b_i\}\subseteq J_i$, $i=1,\dots,n-1$,
such that $(a_1,b_1),\dots,(a_{n-1},b_{n-1})$ are edges
of a spanning tree in $K_{n}$.
\end{enumerate}
\end{proposition}

\begin{proof}  It is clear that (2) implies (1).
On the other hand, (1) implies (2) according to usual Hall's theorem. 
We leave it as an exercise for the reader to check that either of these
two conditions is equivalent to (3).
\end{proof}

We will refer to the three equivalent conditions
in Proposition~\ref{prop:dragon-marriage}
as the {\it dragon marriage condition.}

\begin{example}
\label{exam:dragon_marriage}
Let $M_n$ be the number of sequences of subsets 
$J_1,\dots,J_{n-1}\subseteq[n]$ satisfying
the dragon marriage condition.  Equivalently, $M_n$ is the
number of bipartite subgraphs $G\subseteq K_{n-1,n}$ such that
for any vertex $j$ in the second part there is a matching in $G$
covering the remaining vertices.
According to Theorem~\ref{th:second-formula} with $y_1=\dots=y_n=1$,
we have
$M_n = (n-1)!\,\Vol P_n(-1,-2,-4,\dots,-2^{n-1})$.
Let us calculate a few numbers $M_n$ using
Theorem~\ref{th:f1}.

\smallskip
\begin{center}
\begin{tabular}{|c||l|l|l|l|l|l|l|}
\hline
$n$   & 2 & 3  & 4    & 5      & 6         &  7             & 8 \\
\hline
$M_n$ & 1 & 13 & 1009 & 354161 & 496376001 &  2632501072321 & 
52080136110870785 \\
\hline
\end{tabular}
\end{center}
\smallskip
\end{example}

\section{Generalized permutohedra}
\label{sec:generalized_permutohedra}

\begin{definition}
Let us define {\it generalized permutohedra} as deformations
of the usual permutohedron, i.e., as polytopes obtained by moving vertices
of the usual permutohedron so that directions of all edges are preserved
(and some of the edges may accidentally degenerate into a singe point);
see Appendix~\ref{sec:appendix-lattice-points}.
In other words, a generalized permutohedron is the convex hull
of $n!$ points $v_w\in\R^n$ labeled by permutations $w\in S_{n}$ 
such that, for any $w\in S_{n}$ and any adjacent transposition $s_i=(i,i+1)$,
we have $v_{w} - v_{w\,s_i} = k_{w,i} (e_{w(i)} - e_{w(i+1)})$,
for some nonnegative number $k_{w,i}\in\R_{\geq 0}$, where $e_1,\dots,e_{n}$
are the coordinate vectors in $\R^{n}$,
cf.\ Proposition~\ref{prop:comb_structute_perm}.
\end{definition}

Each generalized permutohedron is obtained by parallel translation
of the facets of a usual permutohedron.  Recall that these facets
are given by Rado's theorem (Proposition~\ref{prop:Rado}).
Thus generalized permutohedra are parametrized by collections
$\{z_I\}$ of the $2^{n}-1$ coordinates $z_I$, for nonempty subsets 
$I\subseteq[n]$, that belongs to a certain 
deformation cone $\D_n$.
Each generalized permutohedron has the form
$$
P_n^z(\{z_I\}) = \left\{(t_1,\dots,t_{n})\in \R^{n}\mid 
\sum_{i=1}^{n}t_i = z_{[n]},\
\sum_{i\in I} t_i \geq z_I,\textrm{ for subsets } I\right\},
$$
for $\{z_I\}\in\D_n$.  
If $z_I = z_J$ whenever $|I|=|J|$, then $P_n(\{z_I\})$ is
a usual permutohedron.

The following figure shows examples of generalized permutohedra:
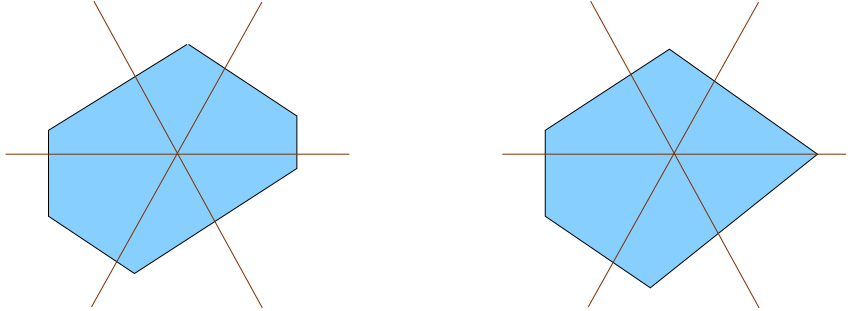
\begin{center}
\end{center}
\begin{center}
\begin{picture}(0,0)%
\includegraphics{fig4-2.pstex}%
\end{picture}%
\setlength{\unitlength}{1579sp}%
\begingroup\makeatletter\ifx\SetFigFont\undefined%
\gdef\SetFigFont#1#2#3#4#5{%
  \reset@font\fontsize{#1}{#2pt}%
  \fontfamily{#3}\fontseries{#4}\fontshape{#5}%
  \selectfont}%
\fi\endgroup%
\begin{picture}(13224,4846)(-11,-4064)
\end{picture}

\end{center}
\medskip

According to Theorem~\ref{th:Todd-Euler-Maclaurin}, we have
the following statement.

\begin{proposition}  The volume of the generalized permutohedron
$P_n(\{z_I\})$ is a polynomial function of the $z_I$'s defined
on the deformation cone $\D_n$.  The number of lattice points
$P_n(\{z_I\})\cap \Z^{n}$ in the generalized permutohedron
is a polynomial function of the $z_I$'s defined
on the lattice points $\D_n\cap\Z^{2^{n}-1}$ of the deformation
cone. 
\end{proposition}

Let us call the multivariate polynomial that expresses the number of 
lattice points in $P_n(\{z_I\})$ the {\it generalized Ehrhart 
polynomial} of the permutohedron.

Let us give a different construction for a class of generalized permutohedra.
Let $\Delta_{[n]} = \textrm{ConvexHull}(e_1,\dots,e_{n})$ be
the standard coordinate simplex in $\R^{n}$.
For a subset $I\subset [n]$, let $\Delta_I = 
\textrm{ConvexHull}(e_i\mid i\in I)$ denote the face
of the coordinate simplex $\Delta_{[n]}$:
$$
\Delta_I = \mathrm{ConvexHull}(e_i\mid i\in I).
$$
Let $\{y_I\}$ be a collection of nonnegative parameters $y_I\geq 0$,
for all nonempty subsets $I\subset [n]$.  
Let us define the polytope
$P_n^y(\{y_I\})$ as the Minkowski sum of the simplices $\Delta_I$
scaled by the factors $y_I$:
$$
P_n^y(\{y_I\}) := \sum_{I\subset[n+1]} y_I\cdot \Delta_I.
$$

\begin{proposition}
\label{prop:Py=Pz}
Let $\{y_I\}$ be a collection of 
nonnegative real numbers for all nonempty subsets $I\subseteq [n]$, and let 
$\{z_I\}$ be the collection of numbers given by 
$$
z_I = \sum_{J\subseteq I} y_J,\quad \textrm{for all nonempty $I\subseteq[n]$}.
$$
Then $P_n^y(\{y_I\} = P_n^z(\{z_I\}$.
\end{proposition}

\begin{proof}  Let us first pick a nonempty subset $I_0\subseteq[n]$
and set $y_I = \delta(I,I_0)$ (Kronecker's delta).  Then
$P_n^y(\{y_I\}) = \Delta_{I_0}$, because the Minkowski contains 
only 1 nonzero term.  In this case, we have $z_I = 1$, if $I\supseteq I_0$,
and $z_I=0$, otherwise.  The inequalities describing the polytope
$P_n^z(\{z_I\})$ give the same coordinate simplex $\Delta_{I_0}$.
The general case follows from the fact that the Minkowski sum
of two generalized permutohedra $P_n^z(\{z_I\})$ and 
$P_n^z(\{z_I'\})$, for $\{z_I\}, \{z_I'\} \in\D_n$,
is exactly the generalized permutohedron $P_n^z(\{z_I+z_I'\})$
parametrized by the coordinatewise sum $\{z_I+z_I'\}\in\D_n$.
This fact is immediate from the definition of $P_n^z(\{z_I\})$.
\end{proof}

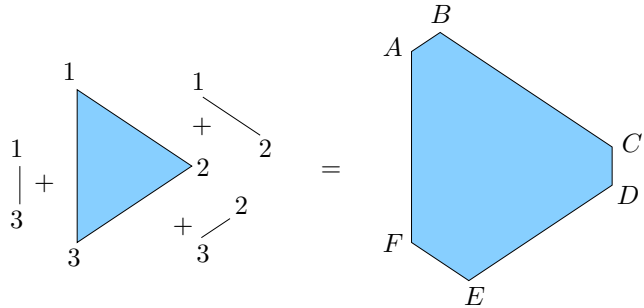
\begin{center}
\input{fig11-1.pstex_t}
\end{center}

\begin{remark}  Not every generalized permutohedron $P_n^z(\{z_I\})$
can be written as a Minkowski sum $P_n^y(\{y_I\})$
of the coordinate simplices. 
For example, for $n=3$, the polytope $P_3^y(\{y_I\})$ (usually a hexagon) 
is the Minkowski sum of the coordinate triangle $\Delta_{[3]}$
and 3 line intervals $\Delta_{\{1,2\}}$, $\Delta_{\{1,3\}}$, $\Delta_{\{2,3\}}$
parallel to its edges (scaled by some factors); see the 
figure above.  For this hexagon we always have $|AB|\leq |DE|$.  
On the other hand, any hexagon with edges parallel to the edges 
of $\Delta_{[3]}$
is a certain generalized permutohedron $P_3^z(\{z_I\})$.

The points $\{z_I\}$ of the deformation cone $\D_n$ that can be
expressed as $z_I=\sum_{J\subseteq I} y_J$ through nonnegative 
parameters $y_I$ form a certain region $\D_n'$ of top dimension in 
the deformation cone $\D_n$.  Since the volume and generalized
Ehrhart polynomial are polynomial functions on $\D_n$, it is
enough to calculate them for the class of polytopes $P_n^y(\{y_I\})$
and then extend from $\D_n'$ to $\D_n$ by the polynomiality.
\end{remark}

In what follows will refer to the polytopes $P_n^y(\{y_I\})$ 
as generalized permutohedra, keeping in mind that they form
a special class of polytopes $P_n^z(\{z_I\})$.

\section{Nested complex}
\label{sec:nested}

The combinatorial structure of the generalized permutohedron
$P_n^y=P_n^y(\{y_I\})$ depends only on the set $B\subset  2^{[n]}$ of nonempty
subsets $I\subseteq[n]$ such that  $y_I>0$.  In this section, we describe the
combinatorial structure of $P_n^y$ when the set $B$ satisfies
some additional conditions.

\begin{definition}
\label{def:building_set}
Let us say that a set $B$ of nonempty subsets in $S$
is a {\it building set\/} on $S$ if it satisfies the 
conditions:
\begin{itemize}
\item[(B1)]
If $I,J\in B$ and $I\cap J\ne\emptyset$,
then $I\cup J\in B$.
\end{itemize}
\begin{itemize}
\item[(B2)]
$B$ contains all singletons $\{i\}$, for $i\in S$.
\end{itemize}
\end{definition}

Condition (B1) is a certain ``connectivity condition'' for building sets.
Note that condition (B2) does not impose any additional restrictions 
on the structure of generalized permutohedra and was added only 
for convenience.  
Indeed, the Minkowski sum of a polytope with $\Delta_{\{i\}}$, 
which is a single point, is just a parallel translation of
the polytope.

Let $B_{\max}\subset B$ be the subset of maximal by inclusion elements in $B$.
Let us say that a building set $B$ is {\it connected\/} if it has a unique
maximal by inclusion element $S$. According to (B1) all elements of $B_{\max}$
are pairwise disjoint.  Thus each building set $B$ is a union
of pairwise disjoint connected building sets, called the {\it connected
components\/} of $B$,
that correspond to elements of $B_{\max}$.

For a subset $C\subset S$, define the {\it induced building set\/} 
as $B|_C=\{I\in B\mid I\subseteq C\}$.

\begin{example}
\label{example:building_set_G}
Let $\Gamma$ be a graph on the set of vertices $S$.
Define the {\it graphical building\/} $B(\Gamma)$ as the set of all nonempty
subsets $C\subseteq S$ of vertices such that the induced graph $\Gamma|_C$
is connected.  Clearly, it satisfies conditions (B1) and (B2).
The building set $B(\Gamma)$ is connected if and only if the graph $\Gamma$
is connected.  The connected components of $B(\Gamma)$ correspond to 
connected components of the graph $\Gamma$.  The induced building set
is the building set for the induced graph: $B(\Gamma)|_C = B(\Gamma|_C)$.
\end{example}

\begin{definition}
\label{def:nested_family}
A subset $N$ in the building set $B$ is called a {\it nested set\/}
if it satisfies the following conditions:
\begin{enumerate}
\item[(N1)]  For any $I,J\in N$, we have either $I\subseteq J$, or
$J\subseteq I$, or $I$ and $J$ are disjoint.
\item[(N2)]  For any collection of $k\geq 2$ disjoint subsets 
$J_1,\dots,J_k\in N$, their union $J_1\cup \cdots \cup J_k$ is not in $B$.
\item[(N3)] $N$ contains all elements of $B_{\max}$.
\end{enumerate}
The {\it nested complex\/} $\N(B)$ is defined as the poset of all nested 
families in $B$ ordered by inclusion.
\end{definition}

Clearly, the collection of all nested sets in $B$ (with elements of $B_{\max}$ removed) 
is a simplicial complex.

\begin{theorem}
\label{th:nested_complex}
Let us assume that the set $B$ associated with a generalized
permutohedron $P_n^y$ is a building set on $[n]$.
Then the poset of faces of $P_n^y$ ordered by reverse inclusion 
is isomorphic to the nested complex $\N(B)$.
\end{theorem}

This claim was independently discovered by E.~M.~Feichtner and
B.~Sturmfels~\cite[Theorem~3.14]{FS}.   They also 
defined objects similar to $B$-forests discussed below; 
see~\cite[Proposition~3.17]{FS}.

\begin{proof}  Each face of an arbitrary polytope can be described
as the set of points of the polytope that minimize a linear function $f$.
Moreover, the face of a Minkowski sum $Q_1+\dots+Q_m$ that minimizes $f$
is exactly the Minkowski sum of the faces of $Q_i$'s that minimize $f$.

Let us pick a linear function $f(t_1,\dots,t_n)=a_1 t_1+ \cdots + a_nt_n$ on
$\R^n$.  It gives an ordered set partition of $[n]$ into a disjoint 
union of nonempty blocks $[n]=A_1\cup \cdots \cup A_s$ such that $a_i=a_j$, 
whenever $i$
and $j$ are in the same block $A_s$, and $a_i<a_j$, whenever $i\in A_s$ and
$j\in A_t$, for $s<t$.
The face of a coordinate simplex $\Delta_I$ that minimizes the linear
function $f$ is the simplex $\Delta_{\widehat{I}}$,
where $\widehat{I}:= I\cap A_{j(I)}$ and $j=j(I)$ is 
the minimal index such that the intersection $I\cap A_j$ is nonempty.
We deduce that the face of $P_n^y$ minimizing $f$ 
is the Minkowski sum $\sum_{I\in B} y_I \,\Delta_{\widehat{I}}$.

We always have $j(I)\geq j(J)$, for $I\subset J$.  
Let $N\subseteq B$ 
be the collection of elements $I\in B$ such that 
$j(I)\gneq  j(J)$, for any $J\supsetneq I$, $J\in B$.
We can also recursively construct the subset $N\subseteq B$, as follows.
First, all maximal by inclusion elements of $B$ should be in $N$.
According to (B1), all other elements of $B$ should belong to one 
of the maximal elements $I_m$.
For each maximal element $I_{m}\in B$, all elements $I\subsetneq I_{m}$
such that $j(I)=j(I_{m})$, i.e., the elements $I$ that have 
a nonempty intersection with $\widehat{I}_m$, do not belong to $N$.
The remaining elements $I\subsetneq I_m$ are exactly the elements 
of the induced building set 
$B|_{I_m\setminus\widehat{I}_m}$. 
Let us repeat the above procedure for each of the induced building sets.  
In other words, find all maximal by inclusion elements $I_{m'}$ in 
$B|_{I_m\setminus\widehat{I}_m}$.  
These maximal elements should be in $N$.  Then, for each maximal element 
$I_{m'}$, construct the induced building set 
$B|_{I_{m'}\setminus \widehat{I}_{m'}}$,
etc.  Let us keep on doing this branching
procedure until we arrive to building sets that consist of singletons,
all of which should be in $N$.

It follows from this branching construction that $N$
is a nested set in $B$.  It is immediate that $N$ satisfies
conditions (N1) and (N3).  If $J_1,\dots,J_k\in N$ are disjoint 
subsets and $J_1\cup \cdots \cup J_k\in B$, $k\geq 2$, 
then we should have included
$J_1\cup \cdots \cup J_k$ in $N$ in the recursive construction,
and then the $J_i$ cannot all belong to $N$.  This implies condition (N2).
It is also clear that, given $N$, 
we can uniquely reconstruct the subset $\widehat{I} \subseteq I$, 
for each $I\in B$.
Indeed, find the minimal by inclusion element $J\in N$ such that 
$J\supseteq I$.  Then $\widehat{J} = J\setminus 
\bigcup_{K\subsetneq J, K\in N}K$ and $\widehat{I}$ is the
intersection of the last set with $I$.
Thus the nested set $N$ uniquely determines the face 
$\sum_{I\in B} y_I\, \Delta_{\widehat{I}}$ of $P_n^y$ that minimizes $f$.

Let us show that, for any nested set $N\in\N(B)$, there exists 
a face of $P_n^y$ associated with $N$.  Indeed, 
let $A_I= I\setminus \bigcup_{J\subsetneq I,J\in N} J$,
for any $I\in N$.  Then $\bigcup_{I\in N} A_I$ is a disjoint decomposition
of $[n]$ into nonempty blocks.  Let us pick any linear order of 
$A_1<\dots<A_s$ of the blocks $A_I$ such that $A_I<A_J$, for $I\subsetneq J$,
and any linear function $f$ on $\R^n$ that gives this set partition,
for example, $f(t_1,\dots,t_n) = \sum_{i,j\in A_i} i\, t_j$.
Then the function $f$ minimizes a certain face $F_N$ of $P_n^y$ and
if we apply the above procedure to $F_N$ we will recover the nested set 
$N$.  We also see from this construction that the face $F_N$ 
contains the face $F_{N'}$ if and only if $N\subseteq N'$.
\end{proof}

We can express the generalized permutohedron 
$P_n^y(\{y_I\})$ as $P_n^z(\{z_I\})$,
where $z_I = \sum_{i: I_i\subseteq I} y_i$; see 
Section~\ref{sec:generalized_permutohedra}.
Let us give an explicit description of its faces.

\begin{proposition}
\label{prop:faces_F_N}
As before, let us assume that $B$ is a building set.
The face $P_N$ of $P_n^y(\{y_I\}) = P_n^z(\{z_I\})$ associated
with a nested set $N\in \N(B)$ is given by
$$
P_N = \{(t_1,\dots,t_n)\in\R^n\mid \sum_{i\in I} t_i = z_I,\textrm{ for }
I\in N; \ \sum_{i\in J} t_i \geq z_J,\textrm{ for }J\in B\}.
$$
The dimension of the face $P_N$ equals $n-|N|$.
In particular, the dimension of $P_n^y(\{y_I\})$ is $n-|B_{\max}|$.
\end{proposition}

\begin{proof}
According to the proof of Theorem~\ref{th:nested_complex},
for a nested set $N\in\N(B)$, we have the disjoint
decomposition $[n]=\bigcup_{I\in N} A_I$ into nonempty blocks, and 
the corresponding face of $P_n^y$ is given by
$$
P_N = \sum_{I\in N,\, J\in B,\, J\cap A_I \ne \emptyset} y_J\,
\Delta_{J\cap  A_I}.
$$
This Minkowski sum involves the terms $\Delta_{A_I}$, among others.
Thus $\dim P_N \geq \dim (\sum_{I\in N} \Delta_{A_I}) = n-|N|$.
It also follows from the construction that $J\cap A_I\ne\emptyset$
implies that $J\subseteq I$.  Thus we have the equality
$\sum_{i\in I} t_i = z_I$, for $I\in N$ and any point
$(t_1,\dots,t_n)\in P_N$.  It follows that the codimension of $P_N$
in $\R^n$ is at least $|N|$.  Together with the inequality for 
the dimension, this implies that $\dim P_N = n-|N|$ and the face 
$P_N$ is described by the above $|N|$ linear equations, as needed.
\end{proof}

Theorem~\ref{th:nested_complex} implies that vertices of $P_n^y$ are
in a bijective correspondence with maximal by inclusion elements 
of the nested complex $\N(B)$.  We will call these elements 
{\it maximal nested families.}  
The following proposition gives their description.

\begin{proposition}  
\label{prop:max_nested_bijection}
A nested set $N\in \N(B)$ is maximal
if and only if, for each $I\in N$, we have
$|A_I| = 1$, where $A_I = I\setminus\bigcup_{J\subsetneq I, J\in N}J$.
For a maximal nested set $N$, the map $I\mapsto i_I$, where $\{i_I\}=A_I$,
is a bijection between $N$ and $[n]$.  
\end{proposition}

\begin{proof}  
According to the proof of Theorem~\ref{th:nested_complex} and
Proposition~\ref{prop:faces_F_N}, a nested set $N\in \N(B)$ is maximal
(and $F_N$ is a point) if and only if $\dim(\sum_{I\in N}\Delta_{A_I}) =
\sum_{I\in N} (|A_I|-1) = 0$, i.e., all $A_I$ should be one elements sets.  The
map $I\mapsto i_I$ is clearly an injection.  On the other hand, for 
any $i\in [n]$ and the minimal by inclusion element $I$ of $N$ that 
contains $i$, we have $I\mapsto i$.
\end{proof}

For a maximal nested set $N\in\N(B)$, let us partially order
the set $[n]$ by $i\geq_N j$ whenever $I\supseteq J$.  
The Hasse diagram of the order ``$\geq_N$'' is a rooted forest,
i.e., a forest with a chosen root in each connected component
and edges directed away from the roots.
The set of such forests can be described, as follows.

For two nodes $i$ and $j$ in a rooted forest, we
say that $i$ is a {\it descendant\/} of $j$
if the node $j$ belongs to the shortest chain connecting $i$ and the root of
its connected component.  In particular, each node is a descendant
of itself.  Let us say that two nodes $i$ and $j$ are 
{\it incomparable\/} if neither $i$ is a descendant of $j$, nor
$j$ is a descendant of $i$.

\begin{definition}
\label{def:B_forests}
For a rooted forest $F$ and a node $i$, let $\desc(i,F)$ be the set of all
descendants of the node $i$ in $F$ (including the node $i$ itself).
Define a {\it $B$-forest\/} as a rooted forest $F$ on
the vertex set $[n]$ such that 
\begin{enumerate}
\item[(F1)] For any $i\in[n]$, we have $\desc(i,F)\in B$.
\item[(F2)] There are no $k\geq 2$ distinct incomparable nodes $i_1,\dots,i_k$ in $F$
such that $\bigcup_{j=1}^k\desc(i_j,F)\in B$.
\item[(F3)] The sets $\desc(i,F)$, for all roots $i$ of $F$, are
exactly the maximal elements of the building set $B$.
\end{enumerate}
Condition (F3) implies that the number of connected components in a $B$-forest equals
the number of connected components of the building set $B$.
We will call such graphs {\it $B$-trees\/} in the case when $B$ is connected.
\end{definition}

\begin{proposition}
\label{prop:max_nested_trees}
The map $N\mapsto F_N$ is a bijection between maximal nested
families $N\in \N(B)$ and $B$-forests.
\end{proposition}

\begin{proof}
The claim is immediate from the above discussion.
Indeed, note that each maximal nested set $N\in\N(B)$ can be reconstructed
from the forest $F=F_N$ as $N=\{\desc(1,F),\dots,\desc(n,F)\}$.
\end{proof}

Let us describe the vertices of the generalized permutohedron
in the coordinates.  

\begin{proposition}
\label{prop:vertices-coord}
The vertex $v_F=(t_1,\dots,t_n)$ of the generalized permutohedron 
$P_n^y$ associated with a $B$-forest $F$
is given by
$t_i = \sum_{J\in B:\,i\in J\subseteq \desc(i,F)} y_{J}$, 
for $i=1,\dots,n$.
\end{proposition}

\begin{proof}
Let $N$ be the maximal nested set associated with the $B$-forest $F$.
By Proposition~\ref{prop:faces_F_N}, the associated vertex 
$v_F= (t_1,\dots,t_n)$ is
given by the $n$ linear equations $\sum_{i\in I} t_i= z_I$,
for each $I\in N$.  
Notice that, for each $J\in B$, there exists a unique $i\in J$
such that $i\in J\subseteq \desc(i,F)$.  Indeed, $\desc(i,F)$
should be the minimal element of $N$ containing $J$.
Thus, for the numbers $t_i$ defined as in 
Proposition~\ref{prop:vertices-coord} and any $I\in N$,  we have 
$$
\sum_{i\in I} t_i = \sum_{i\in I}\  \sum_{J\in B:\,i\in J\subseteq 
\desc(i,F)} 
y_J = \sum_{J\subseteq I} y_J = z_I,
$$
as needed.
\end{proof}

\begin{proposition} 
Let $F$ be a $B$-forest and let $v_F$ be the associated vertex of the generalized 
permutohedron $P_n^y$.
For each nonrooted node $i$ of $F$, define the $n$-vector $g_{i,F} = e_i - e_j$,
where the node $j$ is the parent of the node $i$ in $F$.
(Here $e_1,\dots,e_n$ are the coordinate vectors in $\R^n$.)
Then the integer vectors $g_{i,F}$ generate the local cone of the 
polytope $P_n^y$ at the vertex $v_F$.
In particular, the generalized permutohedron $P_n^y$ is a simple Delzant polytope;
see Appendix~\ref{sec:appendix-lattice-points}.
\end{proposition}

\begin{proof}  Let $N$ be the maximal nested set associated with the forest $F$.
Then each edge of $P_n^y$ incident to $v_F$ correspond to a nested sets 
obtained from $N$ by removing an element $I\in N\setminus B_{\max}$.
There are $n-|B_{\max}|$ such edges and 
Proposition~\ref{prop:faces_F_N} implies that they are generated by the vectors
$g_{i,F}$.
\end{proof}

Let $f_B(q)$ be the {\it $f$-polynomial} of the generalized permutohedron
$P_n^y$.  According to 
Theorem~\ref{th:nested_complex} is is given by
$$
f_B(q) = \sum_{i=0}^{n-1} f_i \,q^i = \sum_{N\in\N(B)} q^{n-|N|},
$$
where $f_i$ is the number of $i$-dimensional faces of $P_n^y$.
The recursive construction of nested families
implies the following recurrence relations fort the $f$-vector.

\begin{theorem}
\label{th:f-recurrence}
The $f$-polynomial $f_B(q)$ is determined by
the following recurrence relations:
\begin{enumerate}
\item If $B$ consists of a single singleton, then $f_B(q)=1$.
\item If $B$ has connected components $B_1,\dots,B_k$, then
$$
f_B(q) = f_{B_1}(q)\cdots f_{B_k}(q).
$$
\item If $B$ is a connected building on $S$, then
$$
f_{B}(q) = \sum_{C\subsetneq S} q^{|S|-|C|-1} f_{B|_C}(q).
$$
\end{enumerate}
\end{theorem}

\begin{definition}
\label{def:gen_Catalan}
We define the {\it generalized Catalan number},
for a building set $B$, as the number $C(B)=f_B(0)$ of vertices 
of the generalized permutohedron $P_n^y$, or, equivalently,
the number of maximal nesting families in $\N(B)$,
or, equivalently, the number of $B$-forests.
\end{definition}

The reason for this name will become apparent from examples in the next
section.  The generalized Catalan numbers $C(B)$ are determined by the
recurrence relations similar to the ones in Theorem~\ref{th:f-recurrence},
where in (3) we sum only over subsets $C\subset S$ of cardinality $|S|-1$.

In the following section we show that the associahedron is a special
case of generalized permutohedra.  Thus we can also call this class
of polytopes {\it generalized associahedra}.  However this name is already 
reserved for a different generalization of the associahedron studied by Chapoton, 
Fomin, and Zelevinsky~\cite{CFZ}.

Even though Chapoton-Fomin-Zelevinsky's generalized associahedra are different from 
our ``generalized associahedra,'' there are some similarities between these two
families of polytopes.  In~\cite{Zel} Zelevinsky gives an alternative
construction for generalized permutohedra associated with building sets which is 
parallel to the construction from~\cite{CFZ}.  He first constructs the dual fan
for the nested complex $\N(B)$ and then shows that it has a polytopal realization.

A natural question to ask is how to find a common generalization of
Chapoton-Fomin-Zelevinsky's generalized associahedra and generalized
permutohedra discussed in this section.

\section{Examples of generalized permutohedra}
\label{sec:examples_of_gen_perm}

\subsection{Permutohedron}

Let us assume that building set $B=B_{all}=2^{[n]}\setminus\{\emptyset\}$ 
is the set of all nonempty
subsets in $[n]$.  Then $P_n^y$ is combinatorially equivalent to
the usual permutohedron, say, $P_n(n,n-1,\dots,1)$.  This
is the generic case of generalized permutohedra.
In this case, nested families are flags of subsets
$J_1\subsetneq J_2\subsetneq \cdots \subsetneq J_s = [n]$.
Indeed, two disjoint subsets $I$ and $J$ cannot belong to a nested
set because their union $I\cup J$ is in $B$.
The maximal nested families are complete flags on $n$ subsets.
Clearly, there are $n!$ such flags, which correspond to the $n!$ vertices
of the permutohedron.  In this case, $B_{all}$-trees are directed chains
of the form $(w_1,w_2),(w_2,w_3),\dots,(w_{n-1},w_n)$, where
$w_1,\dots,w_n$ is a permutation in $S_n$.  The generalized 
Catalan number in this case is $C(B_{all}) = n!$.  

\subsection{Associahedron}
\label{ssec:associahedron}

Assume that the building set  $B=B_{int}=\{[i,j]\mid 1\leq i\leq j\leq n\}$ 
is the set of all continuous intervals in $[n]$.
In this case, the generalized permutohedron is combinatorially equivalent
to the {\it associahedron,} also known as the {\it Stasheff polytope,}
which first appeared in the work of Stasheff~\cite{Sta}.

A nested set $N\subseteq B_{int}$ is a collection of intervals such that,
for any $I, J\in N$, we either have $I\subseteq J$,
$J\subseteq I$, or $I$ and $J$ are disjoint {\it non-adjacent\/} intervals,
i.e., $I\cup J$ is not an interval.
Let us describe $B_{int}$-trees.

Recall that a {\it plane binary tree\/} is an tree such 
that each node has at most 1 left child and at most
one right child.  (If a node has only one child, we specify
if it is the left or the right child.)  It is well
known that there are the Catalan number $C_n = \frac{1}{n+1}
\binom{2n}{n}$ of plane binary trees on $n$ unlabeled nodes.

For a node $v$ in such a tree, let $L_v$ be the left branch and $R_v$ be the
right branch at this node, both of which are smaller plane binary trees.
If $v$ has no left child, then $L_v$ is the empty
graph, and similarly for $R_v$.  For any plane binary tree on $n$ nodes,
there is a unique way to label the nodes by the numbers $1,\dots,n$ so that,
for any node $v$, all labels in $L_v$ are less than the label of $v$ and all
labels in $R_v$ are greater than the label of $v$.  Indeed, label each node
$v$ by the number $|L_v|+1$.  

We can also describe this labeling using the {\it depth-first search}.
This is the walk on the nodes of a tree that starts at the root
and is determined by the rules: (1) if we are at a some node and have never
visited its left child, then go to the left child; (2) otherwise, if we have
never visited its right child, then go to the right child; (3) otherwise, if
the node has the parent, then go to the parent; (4) otherwise stop.  
Let us mark the nodes by the integers $1,\dots,n$ in the order
of their appearance in this walk, as follows.
Each time when we visit an unmarked vertex and {\it do not\/} apply rule (1), 
we mark this node.
The labeling of nodes defined by any of these equivalent ways 
is called the {\it binary search labeling}.
It was described by Knuth in~\cite[6.2.1]{Knu}.
Example~\ref{ex:plane_binary_trees} below shows a plane binary tree 
with the binary search labeling.

\begin{proposition} The $B_{int}$-trees are exactly plane binary 
trees on $n$ nodes with the binary search labeling.
\end{proposition}

\begin{proof}  Let $N$ be a maximal nested set.
Suppose that the maximal element $[n]\in N$ corresponds to
$i=i_{[n]}$ under the bijection in Proposition~\ref{prop:max_nested_bijection}.
Then $N\setminus [n]$ is the union of two maximal nested families
on $[1,i-1]$ and on $[i+1,n]$.  Equivalently, each $B_{int}$-tree
is a rooted tree with root labeled $i$ and two branches which
are $B_{int}$-trees on the vertex sets $[1,i-1]$ and $[i+1,n]$.
This implies the claim.
\end{proof}

Thus in this case the generalized permutohedron has the Catalan number
$C_n$ vertices associated with plane binary trees.  
Proposition~\ref{prop:vertices-coord} implies the following 
description of the vertices of $P_n^y(\{y_{ij}\})$, 
where $y_{ij} = y_{[i,j]}$ for each interval $[i,j]\subseteq [n]$.
For a plane binary trees $T$ with binary search labeling,
let $\desc(k,T)=[l_k,r_k]$, for $k=1,\dots,n$.
Then the left branch of a vertex $k$ is $L_k = [l_k,k-1]$
and the right branch is $R_k = [k+1,r_k]$.

\begin{corollary}
\label{cor:ass_realization}
The vertex $v_T=(t_1,\dots,t_n)$ associated
with a plane binary tree $T$ is given by
$t_k = \sum_{l_k\leq i\leq k\leq j\leq r_k} y_{ij}$.
In particular, in the case when $y_{ij}=1$, for any $1\leq i\leq j\leq n$,
we have
$$
v_T=((|L_1|+1)(|R_1|+1),\cdots, ((|L_n|+1)(|R_n|+1)).
$$
\end{corollary}

The polytope $\Ass_n$ with the $C_n$ vertices given by the second part of
Corollary~\ref{cor:ass_realization} is exactly the realization of the
{\it associahedron\/} described by Loday~\cite{Lod}.  
The will refer to this particular
geometric realization of the associahedron as the {\it Loday realization.}
This polytope can be equivalently defined as the Newton polytope
$\Ass_n := \Newton\left(\prod_{1\leq i\leq j\leq n} 
(t_i+t_{i+1}+\cdots +t_j)\right)$.
We will calculate volumes and numbers of lattice points in 
$\Ass_n$, for $n=1,\dots,8$, in 
Examples~\ref{exam:assoc_volume} and~\ref{exam:lattice_points_ass}.

We can also describe the Loday realization, as follows.
There are $C_n$ subdivisions of the triangular
Young diagram of the shape $(n,n-1,\dots,1)$ into a disjoint
union of $n$ rectangles; see Thomas~\cite[Theorem~1.1]{Tho}
and Stanley's Catalan addendum~\cite[Problem 6.19($u^5$)]{St3}.
These subdivisions are in simple a bijective correspondence with
plane binary trees on $n$ nodes.  
The $i$-th rectangle in such a subdivision is the rectangle that
contains the $i$-th corner of the triangular shape.
Then, for a vertex $v_T=(t_1,\dots,t_n)$
of the associahedron in the Loday realization, the $i$-the coordinate 
$t_i$ equals the number of boxes in the $i$-th rectangle of the
associated subdivision;
see Example~\ref{ex:plane_binary_trees} below.

\begin{example}
\label{ex:plane_binary_trees}
Here is an example of a plane binary tree $T$ with the binary search labeling:
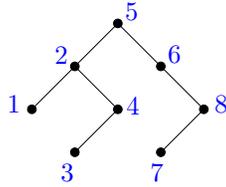
\begin{center}
\input{fig12-1.pstex_t}
\end{center}
This tree is associated with the maximal nested set
$$
N=\{\desc(1,T),\dots,\desc(8,T)\}=
\{[1,1],[1,4],[3,3],[3,4],[1,8],[6,8],[7,7],[7,8]\}.
$$
This tree corresponds to the following subdivision of the triangular
shape into rectangles.  (Here we used shifted Young diagram notation
for a future application; see Section~\ref{sec:shifted_tableaux}.)
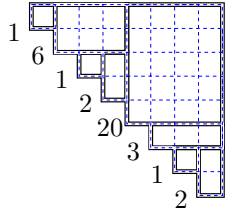
\begin{center}
\input{fig13-1.pstex_t}
\end{center}
The corresponding vertex of the associahedron in the Loday realization
is
$$
(1\cdot 1 \,,\ 2\cdot 3\,,\ 1\cdot 1\,,\ 2\cdot 1\,,\ 5\cdot 4\,,\ 
1\cdot 3\,,\ 1\cdot 1\,,\ 2\cdot 1).
$$
\end{example}

\begin{example}
The next figure shows the Loday realization of the associahedron for $n=3$:
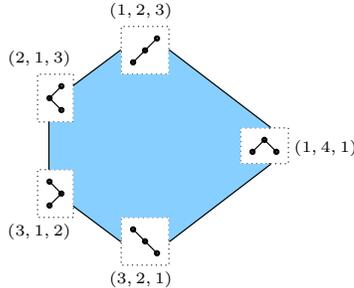
\begin{center}
\input{fig14-1.pstex_t}
\end{center}
\end{example}

\subsection{Cyclohedron}

Let $B=B_{cyc}$ be the set of all {\it cyclic intervals\/} in $[n]$,
i.e., subsets of the form $[i,j]$ and $[1,i]\cup [j,n]$, for 
$1\leq i\leq j \leq n$.
In this case, the generalized permutohedron is the 
{\it cyclohedron\/} that was also introduced by Stasheff~\cite{Sta}.
If we restrict the building set $B_{cyc}$ to $[n]\setminus \{i\}$,
then we obtain the building set isomorphic to the set $B_{int}$
of usual intervals in $[n-1]$.  Thus we obtain the following
description of $B_{cyc}$-trees.

\begin{proposition} The set of $B_{cyc}$-trees is exactly the set
of trees that have a root at some vertex $i$ attached to a plane binary 
tree on $n-1$ nodes with the binary search labeling by integers in 
$[n]\setminus\{i\}$ with respect to the order 
$i+1<i+2<\cdots <n < 1<\cdots<i-1$. 
\end{proposition}

The generalized Catalan number in this case is
$C(B_{cyc}) = n\cdot C_{n-1} = \binom{2n-2}{n-1}$.

\subsection{Graph associahedra}
\label{ssec:graph_ass}

Let $\Gamma$ be a graph on the vertex set $[n]$.  Let us assume that $B=B(\Gamma)$ is 
is the set of subsets $I\subseteq [n]$ such
that the induced graph $\Gamma|_I$ is connected; see
Example~\ref{example:building_set_G}.  In this case, the generalized
permutohedron $P_n^y$ is called the {\it graph associahedron.}  
The above examples are special cases of graph associahedra.
If $\Gamma=A_n$ is the chain with $n$ nodes, i.e., the type $A_n$ Dynkin 
diagram, then we obtain the usual associahedron discussed above.   
In the case of the complete graph $\Gamma=K_n$ we obtain the
usual permutohedron.  If $\Gamma$ is the $n$-cycle, then we obtain 
the cyclohedron.

Various graph associahedra, especially those graph associahedra
that correspond to Dynkin diagrams and extended Dynkin diagrams, 
came up earlier in the work of De~Concini and Procesi~\cite{DP1}
on wonderful models of subspace arrangements and then in the 
work on Davis-Januszkiewitz-Scott~\cite{DJS}.
The class of graph associahedra was independently
discovered by Carr and Devadoss in~\cite{CD}. 
They constructed these polytopes using blow-ups, cf.~\cite{DJS}.
These polytopes also recently appeared in 
the paper by Toledano-Laredo~\cite{Tol}
under the name De~Concini-Procesi associahedra.
We borrowed the term graph associahedra from~\cite{CD}.
Since they are special cases of our generalized permutohedra,
we can also call them {\it graph permutohedra.}

In the case of graph associahedra it is enough to require condition (N2) of 
Definition~\ref{def:nested_family} and condition (F2) of Definition~\ref{def:B_forests}
only for $k=2$.
Indeed, if we have several disjoint subsets
$I_1,\dots,I_k\in B(\Gamma)$ such that $\Gamma|_{I_1\cup \cdots \cup I_k}$ is connected,
then $\Gamma|_{I_i\cup I_j}$ is connected for some pair $i$ and $j$.

\begin{definition}
For a graph $\Gamma$, let us define the {\it $\Gamma$-Catalan number\/} as $C(\Gamma) =
C(B(\Gamma))$, i.e., it the number of vertices of the graph associahedron, or,
equivalently, the number of $B(\Gamma)$-trees; see
Definition~\ref{def:gen_Catalan}.  
\end{definition}

For the $n$-chain $\Gamma=A_n$, i.e., the Dynkin diagram of the type $A_n$, the
$A_n$-Catalan number is the usual Catalan number $C(A_n)=C_n$.  For the
complete graph, we have $C(K_n)=n!$.  Let us calculate several other
$G$-Catalan numbers.

Let $T_{n_1, \dots ,n_r}$ be the star graph that has a central node with $r$ 
attached chains with $n_1,\dots,n_r$ nodes.
For example, $T_{1,1,1}$ 
is the Dynkin diagram of the type $D_4$.  

\begin{proposition}  For a positive integer $r$, the generating
function $\tilde C(x_1,\dots,x_r)$
for the $T_{n_1,\dots,n_r}$-Catalan numbers is given by
$$
\sum_{n_1,\dots,n_r\geq 0} C(T_{n_1,\dots,n_r})\, x_1^{n_1} \dots x_r^{n_r}= 
\frac{C(x_1)\cdots\, C(x_r)}
{1-x_1\,C(x_1)-\cdots - x_r C(x_r)},
$$
where $C(x)=\sum_{n\geq 0} C_n x^n = \frac{1-\sqrt{1-4x}}{2x}$
is the generating function for the usual Catalan numbers.
\end{proposition}

\begin{proof}
According to the recurrence relation in Theorem~\ref{th:f-recurrence},
we have
\begin{equation}
\label{ex:CTnnn}
C(T_{n_1,\dots, n_r}) = C_{n_1}\cdots C_{n_r} +
\sum_{k=1}^r \sum_{i=1}^{n_k} 
C(T_{n_1,\dots ,n_{k-1}, n_k-i ,n_{k+1},\dots ,n_r})
\cdot C_{i-1}.
\end{equation}
Indeed, the first term corresponds to removing the central node
and splitting the graph $T_{n_1,\dots,n_k}$ into $r$ chains.
The remaining terms correspond to removing a node in one of the chains
and splitting the graph into two connected components.
This relation can be written in terms of generating functions as
$$
\tilde C(x_1,\dots,x_r) = C(x_1)\dots C(x_r) + \sum_{k=1}^r x_k\cdot
\tilde C(x_1,\dots,x_r)\cdot C(x_k),
$$
which is equivalent to the claim.
\end{proof}

Let us calculate $\Gamma$-Catalan numbers for
a class of graphs which includes all Dynkin diagrams.
Let $D_n = T_{1,1,n-3}$, $\hat A_n$ be the $(n+1)$-cycle,
$E_n= T_{1,2,n-4}$.

\begin{proposition} The $\Gamma$-Catalan numbers for these graphs are given by
$$
\begin{array}{l}
C(A_n) = C_n = \frac{1}{n+1}\binom{2n}{n},
\textrm{ for } n\geq 1,
\\[.05in]
C(\hat A_n) = (n+1)C_n = \binom{2n}{n},
\textrm{ for } n\geq 3,
\\[.05in]
C(D_n) = 2\,C_n -2\,C_{n-1}-C_{n-2}, \textrm{ for } n\geq 3,
\\[.05in]
C(E_n) = 3\,C_n -4\,C_{n-1}-3\,C_{n-2} - 2\,C_{n-3}, \textrm{ for } n\geq 4.
\end{array}
$$
\end{proposition}

\begin{proof}
We already proved that $C(A_n)=C_n$.
Using Theorem~\ref{th:f-recurrence}, we deduce that
$C(\hat A_n) = (n+1)C(A_n)$.
According Theorem~\ref{th:f-recurrence} or~(\ref{ex:CTnnn}), 
we deduce that the numbers $C(D_n)$ can be calculated
using the recurrence relations
$C(D_n) = C_{n-3} + 2C_{n-1} + \sum_{i=1}^{n-3} C(D_{n-i})\, C_{i-1}$,
for $n\geq 4$, and $C(D_3) = 5$.
In order to prove that $C(D_n) = 2C_n - 2C_{n-1} - C_{n-2}$
it is enough to check that the right-hand side satisfy the recurrence
this relation and that $2C_3 -2 C_2 - C_1 = 5$.
We can easily do this using the recurrence relation for the
Catalan numbers
$C_n = \sum_{i=1}^n C_{n-i} C_{i-1}$, for $n\geq 1$.
Similarly, the numbers $C(E_n)$ are given by the recurrence relation
$C(E_n) = C_{n-1} + C(D_{n-1}) + C_{n-2} + 2\,C_{n-4} +
\sum_{i=1}^{n-4} C(E_{n-i})\,C_{i-1} = 
3\, C_{n-1} - C_{n-2} - C_{n-3} + 2\,C_{n-4} +
\sum_{i=1}^{n-4}  C(E_{n-i})\,C_{i-1}$, for $n\geq 5$,
and $C(E_4)=14$.
Again, we can easily check that the right hand side 
of $C(E_n) = 3\,C_n - 4\,C_{n-1} -3 \,C_{n-2} - 2 \,C_{n-3}$
satisfies this relation, and that $3\,C_4 - 4\,C_3 - 3\,C_2 - 2\,C_1 = 14$.
\end{proof}

Similarly, for any fixed $n_1,\dots,n_{k-1}$, the number
$f_n=C(T_{n_1,n_2,\dots,n_{k-1},n})$ can be expressed as a linear 
combination of several Catalan numbers.

\begin{remark}
One can define the generalized Catalan number for any Lie type.
However this number does not depend on multiplicity of edges in the Dynkin
diagram.  The Catalan number for the Lie types $B_n$ and $C_n$ is the 
usual Catalan number $C_n$.
\end{remark}

\subsection{Pitman-Stanley polytope}
\label{ssect:Pitman-Stanley}

All above examples are special cases of graph associahedra.
Let us consider an example that does not belong to this class.

Let $B=B_{\flag} =\{[1],[2],\dots,[n]\}$
be the complete flag of subsets in $[n]$, and 
let $z_i=\sum_{j=1}^i y_{[j]}$, for $i=1,\dots,n$.
According to Proposition~\ref{prop:Py=Pz},
the generalized permutohedron is this case is the polytope
given by the inequalities:
$$
\{(t_1,\dots,t_n)\mid t_i\geq 0,\ t_1+\cdots + t_i\geq z_i,
\textrm{ for }i=1,\dots,n-1,\ t_1+\dots+t_n = z_n\}
$$
This is exactly the polytope studied by Pitman and Stanley~\cite{PiSt}.
We will call it the {\it Pitman-Stanley polytope.}

Let $B_{\flag}^+ = B_{\flag}\cup \{\{1\},\dots,\{n\}\}$
be the set obtained from $B_{\flag}$ by adding all singletons.
The generalized permutohedron for $B_{\flag}^+$ is just a parallel
translation of the Pitman-Stanley polytope.
The set $B_{\flag}^+$ is a building set.  Nested families 
$N\in \N(B_{\flag}^+)$ are the subsets 
$N\subset B_{\flag}^+$ such that (1) if $[i]\in N$ then $\{i+1\}\not\in N$,
and (2) $[n]\in N$.
Let us encode a nested set $N$ by a word $u_1,\dots,u_{n-1}$
in the alphabet $\{0,1,*\}$ such that, for $i=1,\dots,n-1$, if 
$[i]\in N$ then $u_i=0$, if $\{i+1\}\in N$ then $u_i=1$,
otherwise $u_i=*$.  This gives a bijection between nested families
and $3^{n-1}$ words of length $n-1$ with these $3$ letters.
A nested set $N$ contains a nested set $N'$ 
whenever the word for $N$ is obtained from the word for $N'$
by replacing some $*$'s with $0$'s and/or $1$'s.
In particular, a nested set is maximal if its words
contains only $0$'s and $1$'s.
Thus the nested complex $\N(B_{\flag}^+)$ is isomorphic to the face lattice 
of the $(n-1)$-dimensional hypercube.

\begin{proposition}  The Pitman-Stanley polytope has $2^{n-1}$ 
vertices and it is combinatorially equivalent to the $(n-1)$-dimensional
hypercube.
\end{proposition}

Thus the generalized Catalan number in this case is 
$C(B_{\flag}^+) = 2^{n-1}$.

\begin{example}
The following figure shows the combinatorial structure of the 
Pitman-Stanley polytope for $n=3$ in terms of nested families.
\medskip

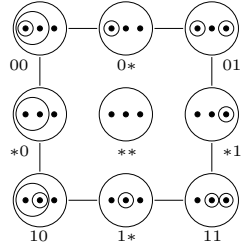
\begin{center}
{\tiny
\input{fig16-1.pstex_t}
}
\end{center}
Note that, as a geometric polytope, the Pitman-Stanley
polytope is a {\it non-regular\/} quadrilateral, as shown on the following figure.
\medskip

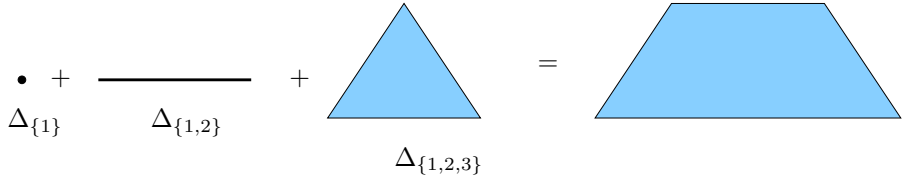
\begin{center}
\input{fig17-1.pstex_t}
\end{center}
\end{example}

\subsection{Graphical zonotope}

Let $\Gamma$ be a graph on the vertex set $[n]$, and 
let $B$ be the set of all pairs $\{i,j\}\subset [n]$ such
that $(i,j)$ is an edge of $\Gamma$.  The set $B$ {\it does not\/} satisfy 
the axioms of a building set; see Definition~\ref{def:building_set}.
The minimal building set that contains $B$ is
the graphical building set $B(\Gamma)$; see Example~\ref{example:building_set_G}.
The generalized permutohedron for the set $B$ is the graphical zonotope
$Z_\Gamma$; see Definition~\ref{def:graph_zonotope}.
In this case, we can not describe combinatorial structure of $Z_\Gamma$ using
nested families.  
However it is well-known that the vertices of $Z_\Gamma$ correspond to 
acyclic orientations of the graph $\Gamma$.  It is not hard to describe the
faces of this polytope as well.  Note that the polytope $Z_\Gamma$ is dual to the 
graphic arrangement for the graph $\Gamma$.

\section{Volume of generalized permutohedra via Bernstein's theorem}
\label{sec:vol_via_Bernstein}

Let $G\subseteq K_{m,n}$ be a bipartite graph with no isolated vertices.
(This graph should not be confused with graphs used in 
Section~\ref{sec:examples_of_gen_perm}.)
We will label the vertices $G$ by $1,\dots,m,\bar 1,\dots,\bar n$ 
and call $1,\dots,m$ the {\it left vertices\/} and  $\bar 1,\dots,\bar n$
the {\it right vertices.}
Let us associate this graph with the collection $\I_G$ of subsets 
$I_1,\dots,I_m\subseteq[n]$ such that
$j\in I_i$ if and only if $(i,\bar j)$ is an edge of $G$.
Let us define the polytope $P_G(y_1,\dots,y_m)$ is the Minkowski sum
$$
P_G(y_1,\dots,y_m) = y_1\Delta_{I_1} + \cdots + y_m \Delta_{I_m}.
$$
The polytope $P_G(y_1,\dots,y_m)$ is exactly the generalized permutohedron
$P_{n}^y(\{y_I\})$, where $y_I = \sum_{i\mid I_i = I} y_i$.

\begin{remark}
\label{rem:P_G(1,1,1)}
The class of polytopes $P_G(1,\dots,1)$
is as general as $P_G(y_1,\dots,y_m)$ for arbitrary nonnegative
integers $y_1,\dots,y_m$.  Indeed, we can always replace 
a term $y_i\Delta_{I_i}$ with $y_i$ terms $\Delta_{I_i}$.
We  use the notation $P_G(y_1,\dots,y_m)$ in order to
emphasize dependence of this class of polytopes on 
the parameters $y_1,\dots,y_m$.
\end{remark}

\begin{definition}
Let us say that a sequence of nonnegative integers $(a_1,\dots,a_m)$ is a 
{\it $G$-draconian sequence\/} if $\sum a_i = n-1$ and, for any subset
$\{i_1<\cdots <i_k\}\subseteq[m]$, we have $|I_{i_1} \cup \cdots \cup I_{i_k}|
\geq a_{i_1}+\cdots + a_{i_k} + 1$.  Equivalently, $(a_1,\dots,a_m)$ is an {\it
$G$-draconian\/} sequence of integers if the sequence of subsets
$I_1^{(a_1)},\dots,I_m^{(a_m)}$, where $I^{(a)}$ means $I$ repeated $a$ times,
satisfies the dragon marriage condition; see
Proposition~\ref{prop:dragon-marriage}.  
\end{definition}

Theorem~\ref{th:second-formula} can be extended to generalized
permutohedra, as follows.

\begin{theorem}  
\label{th:second-formula-generalized}
The volume of the generalized permutohedron
$P_G(y_1,\dots,y_m)$ equals
$$
\Vol P_G(y_1,\dots,y_m) = \sum_{(a_1,\dots,a_m)} 
\frac {y_1^{a_1}}{a_1!} \cdots
\frac {y_m^{a_m}}{a_m!},
$$
where the sum is over all $G$-draconian sequences 
$(a_1,\dots,a_m)$.
\end{theorem}

We can also reformulate Theorem~\ref{th:second-formula-generalized}, 
as follows.

\begin{corollary}
\label{cor:vol_gen_perm}
The volume of the generalized permutohedron 
$P_{n}^y(\{y_I\})$ is given by 
$$
\Vol P_{n}^y(\{y_I\})= \frac{1}{(n-1)!}\sum_{(J_1,\dots,J_{n-1})}
y_{J_1}\cdots y_{J_{n-1}},
$$
where the sum is over ordered collections of nonempty subsets 
$J_1,\dots,J_{n-1}\subset [n]$
such that, for any distinct $i_1,\dots,i_k$, we have
$|J_{i_1}\cup \cdots \cup J_{i_k}| \geq k+1$. 
\end{corollary}

\begin{proof} 
Assume in Theorem~\ref{th:second-formula-generalized}
that $G$ is the bipartite graph associated with the collection 
$I_1,\dots,I_m$, $m=2^n-1$, of all nonempty subsets in $[n]$.
Then replace the summation over $G$-draconian sequences
$(a_1,\dots,a_m)$ by the summation over $\binom {n-1}{a_1,\dots,a_m}$
rearrangements $(J_1,\dots,J_{n-1})$ of the sequence
$(I_1^{(a_1)},\dots,I_m^{(a_m)})$. 
\end{proof}

\begin{example}  
\label{exam:regular_permuitohedron-K_n}
Suppose that $I_1,\dots,I_m$, $m=\binom n2$, is the collection of all 
$2$-element
subsets in $[n]$ and $G\subset K_{m,n}$ is the associated bipartite graph.
Then $P_G(1,\dots,1)$ is the regular permutohedron $P_{n-1}(n-1,n-2,\dots,0)$.
In this case, there are $n^{n-2}$ $G$-draconian sequences $(a_1,\dots,a_m)$, 
which are in a bijective correspondence with trees on $n$ vertices.
For a tree $T\subset K_n$, the $a_i$'s corresponding to the edges of $T$
are equal to $1$ and the remaining $a_i$' are zero, 
cf.\ Proposition~\ref{prop:dragon-marriage}.
Thus we recover the result that $\Vol P_{n}(n-1,n-2,\dots,0)  = n^{n-2}$.
\end{example}

\begin{definition} A sequence of positive integers $(b_1,\dots,b_m)$ 
is called a {\it parking function\/} if its increasing rearrangement
$c_1\leq c_2\leq \cdots \leq c_m$ satisfies $c_i\leq i$, for $i=1,\dots,m$.
\end{definition}

Recall that there are $(m+1)^{m-1}$ parking functions of the length $m$.

\begin{example}
\label{exam:pitman-stanley-volume}
Suppose that $I_i=[n+1-i]$, for $i=1,\dots,m$,
where $m=n-1$.
In this case, the generalized permutohedron $P_{G}(y_1,\dots,y_{m})$
is the Pitman-Stanley polytope; see Subsection~\ref{ssect:Pitman-Stanley}.
A $G$-draconian sequence is a nonnegative integer sequence
$(a_1,\dots,a_{m})$ such that $a_1+\cdots + a_i \geq i$, for $i=1,\dots,m$,
and $a_1+\dots+a_{m}=m$.  There are the Catalan number 
$C_{m} = \frac{1}{m+1}\binom {2m}{m}$ such sequences.
Let us call them {\it Catalan sequences.}
A collection of intervals $I_{b_1},\dots,I_{b_{m}}$ satisfies 
the dragon marriage condition if and only if $(b_1,\dots,b_{m})$ 
is a parking function.
We recover the following two formulas for the volume of the Pitman-Stanley
polytope proved in~\cite{PiSt}:
$$
\Vol P_G(y_1,\dots,y_m)= 
\sum_{(a_1,\dots,a_m)} \frac{y_1^{a_1}}{a_1!}\cdots \frac{y_m^{a_m}}{a_m!}
=\frac{1}{m!} \sum_{(b_1,\dots,b_m)} y_{b_1} \cdots y_{b_m},
$$
where the first sum is over Catalan sequences
$(a_1,\dots,a_m)$ and the second sum is over parking functions
$(b_1,\dots,b_m)$.
In particular, $\Vol P_G(1,\dots,1) = \frac{ (m+1)^{m-1}}{m!} =
\frac{n^{n-2}}{(n-1)!}$.
\end{example}

The proof of Theorem~\ref{th:second-formula-generalized}
relies on Bernstein's theorem
on systems of polynomial equations.  
Let us first recall the definition of the {\it mixed volume}
$\Vol(Q_1,\dots,Q_n)$ of $n$ polytopes $Q_1,\dots,Q_n\subset \R^n$.
It is based on the following proposition.

\begin{proposition}
\label{prop:mixed_volume}
There exists a unique function 
$\Vol(Q_1,\dots,Q_n)$ defined on $n$-tuples of polytopes in $\R^n$
such that, for any collection of $m$ polytopes $R_1,\dots,R_m\subset\R^n$, 
the usual volume of the Minkowski sum 
$y_1 R_1+ \cdots + y_m R_m$, for nonnegative factors $y_i$, 
is the polynomial in $y_1,\dots,y_m$ given by
$$
\Vol(y_1 R_1+ \cdots + y_m R_m) =
\sum_{(i_1,\dots,i_n)} \Vol(R_{i_1},\dots,R_{i_n})\,
y_{i_1}\cdots y_{i_n},
$$
where the sum is over ordered sequences $(i_1,\dots,i_n)\in[m]^n$.
\end{proposition}

For a finite subset $A\subset\Z^n$, let $f_A(t_1,\dots,t_n)
=\sum_{a\in A} \beta_a\, t_1^{a_1}\cdots t_n^{a_n}$
be a Laurent polynomial in $t_1,\dots,t_n$ with 
some complex coefficients $\beta_a$.

\begin{theorem}
\label{th:bernstein}
{\rm Bernstein~\cite{Ber}}  \
Fix $n$ finite subsets $A_1,\dots,A_n\subset\Z^n$.
Let $Q_i$ be the convex hull of $A_i$, for $i=1,\dots,n$.
Then the system 
$$
\left\{
\begin{array}{c}
f_{A_1}(t_1,\dots,t_n) = 0,\\
\vdots \\
f_{A_n}(t_1,\dots,t_n) = 0
\end{array}
\right.
$$
of $n$ polynomial equations in the $n$ variables $t_1,\dots,t_n$
has exactly $n!\,\Vol(Q_1,\dots,Q_n)$ isolated solutions 
in $(\C\setminus\{0\})^n$ whenever the collection
of all coefficients of the polynomials $f_{A_i}$ 
belong to a certain Zariski open set in $\C^{\sum|A_i|}$.
\end{theorem}

Bernstein's theorem is usually used for finding the number 
of solutions of a system of polynomial equations by calculating
the mixed volume.
We will apply Bernstein's theorem in the opposite direction.
Namely, we will calculate the mixed volume by solving a system
of polynomial equations.  Actually, in our case we need
to solve a system of linear equations.

\begin{proof}[Proof of Theorem~\ref{th:second-formula-generalized}]
According to Proposition~\ref{prop:mixed_volume}
and the definition of the polytope $P_G(y_1,\dots,y_m)$ as the 
Minkowski sum of simplices, we have
$$
\Vol P_G(y_1,\dots,y_m) = \sum_{i_1,\dots,i_{n-1}}
\Vol(\Delta_{I_{i_1}},\dots,\Delta_{I_{i_{n-1}}})\,y_{i_1}\cdots y_{i_{n-1}},
$$
where the sum is over all $i_1,\dots,i_{n-1}\in[m]$.
Here we can define $(n-1)$-dimensional (mixed) volumes of polytopes
embedded into $\R^{n}$ as (mixed) volumes of their projections 
into, say, the first $n-1$ coordinates.
It remains to show that the mixed volume of several coordinate simplices 
is equal to
$$
\Vol(\Delta_{J_1},\dots,\Delta_{J_{n-1}}) = 
\left\{
\begin{array}{cl}
\frac{1}{(n-1)!} &\textrm{if $J_1,\dots,J_{n-1}$ satisfy DMC}, \\[.1in]
0 &\textrm{otherwise,}
\end{array}
\right.
$$ 
where ``DMC'' stands for the dragon marriage condition; see
Proposition~\ref{prop:dragon-marriage}.
Consider the following system of $n-1$ linear equations in
the variables $t_1,\dots,t_{n-1}$:
$$
\left\{
\begin{array}{c}
\sum_{j\in J_1} \beta_{1,j}\, t_j = 0,\\
\vdots\\
\sum_{j\in J_{n-1}} \beta_{n-1,j}\, t_j = 0,
\end{array}
\right.
$$
where we assume that $t_{n} = 1$.  According to 
Bernstein's theorem, this system has exactly 
$(n-1)!\, \Vol(\Delta_{J_1},\dots,\Delta_{J_{n-1}})$ isolated solutions 
in $(\C\setminus\{0\})^{n-1}$ for generic coefficients $\beta_{i,j}\in\C$,
for $j\in J_i$.

Of course, we can always solve this linear system using Cramer's rule.
Let $B=(\beta_{ij})$ be the $(n-1)\times n$-matrix formed by
the coefficients of the system, where we assume that $\beta_{i,j}=0$,
for $j\not\in I_i$; and let $|B^{(i)}|$ be the $i$-th maximal minor 
of this matrix.  The system in nondegenerate if and only if 
$|B^{(n)}|\ne 0$.  In this case, its only solution is given by
$t_i = (-1)^i |B^{(i)}|/|B^{(n)}|$, for $i=1,\dots,n-1$.
Thus the system has a single isolated solution in $(\C\setminus\{0\})^{n-1}$ 
if and only if {\it all} $n$ maximal minors of $B$ are nonzero.  
Otherwise, the system has no isolated solutions in $(\C\setminus\{0\})^{n-1}$.

The matrix $B=(\beta_{i,j})$ is 
subject to the only constraint $\beta_{i,j}=0$, for $j\not\in J_i$.
For generic values of $\beta_{i,j}$, the $k$-th maximal minor of 
this matrix is nonzero if and only if there is a system of 
distinct representatives of $J_1,\dots,J_{n-1}$ that avoids $k$.
According to Proposition~\ref{prop:dragon-marriage}, these 
conditions are equivalent to the needed condition. 
This finishes the proof.
\end{proof}

\section{Volumes via Brion's formula}
\label{sec:vol_via_Brion}

Let us give a couple of alternative formulas for volume of generalized
permutohedra that extend results of Section~\ref{sec:descent_div_sym}.
It is more convenient to expresses generalized permutohedra
in the form $P_n^y(\{y_I\})$;  
see Section~\ref{sec:generalized_permutohedra}.

\begin{theorem}
\label{th:gen-perm-sum-w}
For any distinct $\lambda_1,\dots,\lambda_{n}$, we have
$$
\Vol P_n^y(\{y_I\}) = \frac{1}{(n-1)!} \sum_{w\in S_{n}} 
\frac{\left(\sum_{I\subseteq [n]} \lambda_{w(\min(I))} y_{w(I)}\right)^{n-1}}
{(\lambda_{w(1)}-\lambda_{w(2)})\cdots (\lambda_{w(n-1)}-\lambda_{w(n)})}.
$$
\end{theorem}

This theorem is deduced from Brion's formula 
(see Appendix~\ref{sec:appendix-lattice-points})
in exactly the same way as Theorems~\ref{th:f1-W} 
and~\ref{th:f1}. 

For example, we have
$$
\Vol P_2^y(\{y_I\})= 
\frac{\lambda_1 y_{\{1\}} + \lambda_2 y_{\{2\}} + \lambda_1 y_{\{1,2\}}}
{\lambda_1-\lambda_2} 
+
\frac{\lambda_2 y_{\{2\}} + \lambda_1 y_{\{1\}} + \lambda_2 y_{\{2,1\}}}
{\lambda_2-\lambda_1}  = y_{\{1,2\}}
$$
Note the terms  $\lambda_i y_{\{i\}}$ make a zero contribution.
Thus in the summation in Theorem~\ref{th:gen-perm-sum-w} we can skip 
singleton subsets $I$.

For a collection of subsets $J_1,\dots,J_{n-1}\subseteq [n]$, 
construct the integer vector $(a_1,\dots,a_{n}) 
= e_{\min(J_1)} + \cdots + e_{\min(J_{n-1})}$. 
Let $I(J_1,\dots,J_{n-1}) = I_{a_1,\dots,a_{n}}$, defined
as in Section~\ref{sec:descent_div_sym}.
Theorem~\ref{th:vol=descent_number} can be extended as follows.

\begin{theorem}
\label{th:gen_descents_sum_w}
   We have
$$
\Vol P_n^y(\{y_I\}) = 
\sum_{J_1,\dots,J_{n-1} \in[n]} (-1)^{|I(J_1,\dots,J_{n-1})|}
\sum_w y_{w(J_1)} \cdots y_{w(J_{n-1})},
$$
where the second sum is over permutations $w\in S_{n}$ with 
the descent set $I(w)=I(J_1,\dots,J_{n-1})$.
\end{theorem}

This result is deduced from Theorem~\ref{th:gen-perm-sum-w}
using the same argument as in the proof of
Theorem~\ref{th:vol=descent_number}. 

Theorem~\ref{th:gen-perm-sum-w} 
is convenient for explicit calculations
of volumes.  Let us give a couple of examples obtained with some help
of a computer. 

\begin{example}
\label{exam:assoc_volume}
Let $A_n = (n-1)!\, \Vol \Ass_n$, where $\Ass_n$ is the associahedron
in the Loday realization; see Subsection~\ref{ssec:associahedron}.
According to Theorem~\ref{th:gen-perm-sum-w} we have
$$
A_n = \sum_{w\in S_n}
\frac{\left(\sum_{1\leq i\leq j\leq n} \lambda_{m(i,j,w)}\right)^{n-1}}
{(\lambda_{w(1)}-\lambda_{w(2)})\cdots (\lambda_{w(n-1)}-\lambda_{w(n)})},
$$
where $m(i,j,w) = w(\min(w^{-1}([i,j])))= \min\{k\mid w(k)\in[i,j]\}$.
The numbers $A_n$, for $n=1,\dots,8$, are the following:

\smallskip
\begin{center}
\begin{tabular}{|c||l|l|l|l|l|l|l|l|}
\hline
$n$   & 1 & 2 & 3 & 4   & 5    & 6      & 7        & 8          \\
\hline 
$A_n$ & 1 & 1 & 7 & 142 & 5895 & 417201 & 45046558 & 6891812712  \\
\hline
\end{tabular}
\end{center}
\smallskip
\end{example}

\begin{example}
{\rm (cf.~Example~\ref{exam:dragon_marriage})} \
Let us call a subgraph $G\subseteq K_{n,n}$ a {\it Hall graph\/}
if it contains a perfect matching or, equivalently, satisfies
the Hall marriage condition.  Let $H_n$ be the number of Hall
subgraphs in $K_{n,n}$.  According to Corollary~\ref{cor:vol_gen_perm},
$\frac{1}{(n-1)!}\, H_{n-1}$ is the volume of the generalized permutohedron
$P_n^y(\{y_I\})$ with $y_I=1$, for subsets $I\subseteq[n]$ such that
$n\in I$, and $y_I=0$, otherwise.
Using Theorem~\ref{th:gen-perm-sum-w} we can calculate several 
numbers $H_n$.

\smallskip
\begin{center}
\begin{tabular}{|c||l|l|l|l|l|l|l|}
\hline
$n$   &  1 & 2 & 3   & 4     & 5        & 6           &  7 \\
\hline
$H_n$ &  1 & 7 & 247 & 37823 & 23191071 & 54812742655 &  494828369491583 \\
\hline
\end{tabular}
\end{center}
\smallskip
\end{example}

\section{Generalized Ehrhart polynomial}
\label{sec:generalized_Ehrhart}

In this section we give a formula for the number of lattice points
of generalized permutohedra.

Let us define the {\it Minkowski difference\/} of two polytopes
$P,Q\subset \R^n$ as $P-Q = \{x\in\R^n\mid x+ Q \subseteq P\}$.
Its main property is the following.

\begin{lemma}
\label{lem:Minkowski_difference}
For any two polytopes, we have $(P+Q)-Q = P$.
\end{lemma}

\begin{proof}  We need to prove that, for a point $x$, 
we have $x+Q\subseteq P+Q$ 
if and only if $x\in P$.  The ``if'' direction is trivial.
Let us check the ``only if'' direction.  It is enough to assume 
that $x=0$.  We need to show that $Q\subseteq P+Q$ implies that 
$0\in P$.  Suppose that $0\not\in P$. Because of convexity of $P$
we can find a linear form $f$ such that $f(p)>0$, for any point $p\in P$
(and, of course, $f(0)=0$).  Let $q_{\min}\in Q$ be the point of $Q$ with
minimal possible value of $f(q_{\min})$.  Then for any point $p+q\in P+Q$,
where $p\in P$ and $q\in Q$, we have $f(p+q) = f(p)+f(q)> f(q_{\min})$.
Thus $q_{\min}\not\in P+Q$.  Contradiction.
\end{proof}

\begin{definition}
Let us define the {\it trimmed generalized permutohedron}
as the Minkowski difference of $P_G(y_1,\dots,y_m)$ 
and the simplex $\Delta_{[n]}$:
$$
P_G^-(y_1,\dots,y_{m})  = P_G(y_1,\dots,y_m) - 
\Delta_{[n]} = \{x\in \R^n\mid x+\Delta_{[n]} \subseteq P_G\}.
$$
\end{definition}

This is a slightly more general class of polytopes than generalized
permutohedra $P_G$.
Suppose that $I_1=[n]$, i.e., the vertex $1$ in $G$ is connected
with all vertices in the right part.  (If this is not the case,
we can always add such a vertex to $G$.)
According to Lemma~\ref{lem:Minkowski_difference}, we have
$$
P_G(y_1,\dots,y_m) = P_G^{-}(y_1+1,y_2,\dots,y_m).
$$
In other words, if one of the summands in the Minkowski sum for $P_G$ is
$\Delta_{[n]}$ then the trimmed generalized
permutohedron $P_G^-$ equals the (untrimmed) generalized permutohedron
given by a similar Minkowski sum without this summand.
Also notice that the class of polytopes $P_G^{-}(1,\dots,1)$ is as general as
$P_G^{-}(y_1,\dots,y_m)$ for arbitrary nonnegative integer $y_1,\dots,y_m$,
cf.~Remark~\ref{rem:P_G(1,1,1)}.

Let us give a formula for the generalized Ehrhart polynomial 
of (trimmed) generalized permutohedra.
Define raising powers as
$(y)_a := y(y+1)\cdots (y+a-1)$, for $a\geq 1$, and $(y)_0 := 1$.
Equivalently, $\frac{(y)_a}{a!} := \binom{y+a-1}{a}$.

\begin{theorem} 
\label{th:gen_ehrhrart}
For nonnegative integers $y_1,\dots,y_m$,
the number of lattice points  in the trimmed generalized 
permutohedron $P_G^-(y_1,\dots,y_m)$ equals
$$
P_G^-(y_1,\dots,y_m)\cap \Z^{n}=
\sum_{(a_1,\dots,a_m)} 
\frac {(y_1)_{a_1}}{a_1!} \cdots
\frac {(y_m)_{a_m}}{a_m!},
$$
where the sum is over all $G$-draconian sequences 
$(a_1,\dots,a_m)$.
In particular, the number of lattice points in 
$P_G(y_1,\dots,y_m)$ equals the above expression with 
$y_1$ replaced by $y_1+1$, assuming that $I_1=[n]$.

This also implies that the number of lattice points in $P_G^{-}(1,\dots,1)$
equals the number of $G$-draconian sequences.
\end{theorem}

In other words, the formula for the number of lattice points
in $P_G^{-}$ is obtained from the formula for the volume 
of $P_G$ by replacing usual powers in all terms by raising powers.
We will prove this theorem in Section~\ref{sec:subdivision}.

\begin{example}  
\label{exam:lattice_points_reg_permut}
Let $I_1=[n]$ and $I_2,\dots,I_m$, $m=\binom{n}{2}+1$, 
be all $2$-element subsets in $[n]$, cf.\ 
Example~\ref{exam:regular_permuitohedron-K_n}. 
Then the polytope $P_G^-(1,\dots,1)$ is the regular 
permutohedron $P_n(n-1,\dots,0)$ and 
$$
P_G^-(0,1,\dots,1)= P_n(n-1,\dots,0)-\Delta_{[n]} = 
P_n(n-2,n-2,n-3,\dots,0).
$$
In this case, $G$-draconian sequences are in a bijection with
forests $F\subset K_n$.  The $G$-draconian sequence $(a_1,\dots,a_m)$ 
associated with a forest $F$ with $c$ connected components 
is given by $a_1 = c-1$, $a_i = 1$ if $I_i$ is an edge of $F$, 
and $a_i=0$ otherwise, for $i=2,\dots,m$. 
Theorem~\ref{th:gen_ehrhrart} implies that 
the number of lattice points of lattice points in the regular
permutohedron equals the number of labeled forests on $n$ nodes.
More generally, if we set some $y_i$'s to zero, then we deduce that 
the number of lattice points in a graphical zonotope equals
the number of forests in the corresponding graph; 
see Proposition~\ref{prop:vol-Z-G}.
\end{example}

Theorem~\ref{th:gen_ehrhrart} and 
Example~\ref{exam:lattice_points_reg_permut} also imply 
the following statement.

\begin{corollary} Let $\Gamma$ be a connected graph on the vertex set $[n]$.
Let $Z_\Gamma$ be the graphical zonotope, i.e.,
the Minkowski sum of intervals $[e_i,e_j]$, for edges $(i,j)$ of $\Gamma$.
Also consider the Minkowski difference $Z_\Gamma^{-} = Z_\Gamma-\Delta_{[n]}$.
Then the volume of $Z_\Gamma$ equals the number of lattice points in $Z_\Gamma^{-}$:
$$
\Vol Z_\Gamma = \#(Z_\Gamma^-\cap \Z^n), 
$$
and both these numbers are equal to the number of spanning trees in the graph $\Gamma$.
In particular, the number of lattice points
in the permutohedron $P_n(n-2,n-2,n-3,\dots,0)$ equals $n^{n-2}$.
\end{corollary}

\begin{example}
Suppose that $I_i=[n+1-1]$, for $i=1,\dots,m$, where $m=n-1$,
as in Example~\ref{exam:pitman-stanley-volume}.
Theorem~\ref{th:gen_ehrhrart} implies the following expression
for the number of lattice points in the Pitman-Stanley polytope
proved in~\cite{PiSt}:
$$
\#(P_G(y_1,\dots,y_m)\cap \Z^n) =
\sum_{(a_1,\dots,a_m)}
\frac{(y_1+1)_{a_1}}{a_1!}\cdots \frac{(y_m)_{a_m}}{a_m!},
$$
where the sum is over Catalan sequences
$(a_1,\dots,a_m)$ 
as in Example~\ref{exam:pitman-stanley-volume}.
Thus the number of lattice points in $P_G^-(1,\dots,1) = P_G(0,1,\dots,1) =
\Delta_{[2]} + \cdots + \Delta_{[n-1]}$ equals the Catalan number $C_m= C_{n-1}$.
Also the number of lattice points in $P_G(1,\dots,1) = \Delta_{[2]}+\cdots +\Delta_{[n]}$
equals the Catalan number $\sum_{(a_1,\dots,a_m)} (a_1+1)  = C_n$, 
where the sum is over Catalan sequences.
\end{example}

For a bipartite graph $G\subseteq K_{m,n}$, let
$G^*\subseteq K_{n,m}$ be mirror image of $G$ obtained
by switching the left and write components.  In other words,
$G^*$ is the same graph with the relabeled vertices
$1,\dots,m,\bar 1,\dots,\bar n \longrightarrow 
\bar 1,\dots,\bar m,1,\dots,n$.

\begin{lemma}
\label{lem:G-draconain=G*} 
The set of $G$-draconian sequences is exactly 
the set of lattice points of the polytope
$P_{G^*}^-(1,\dots,1)\subset\R^m$.
\end{lemma}

\begin{proof}  In order to prove the lemma we just need to 
check all definitions.  Let $I_1^*,\dots,I_n^*\subseteq[m]$ be the 
collection of subsets associated with the graph $G^*$,
i.e., $j\in I_i^*$ whenever $(i,\bar j)\in G^*$, or, equivalently,
$(j,\bar i)\in G$.  Then $P_{G^*}(1,\dots,1) = 
\Delta_{I_1^*}+\cdots +\Delta_{I_n^*}\subseteq \R^m$.  This is
exactly the polytope $P_m^z(\{z_I\}$, 
where $z_I=\#\{i\mid I_i^*\subseteq I\}$, for nonempty 
$I\subseteq[m]$; see Proposition~\ref{prop:Py=Pz}.  
According to Section~\ref{sec:generalized_permutohedra},
this polytope
is given by the inequalities
$$
P_{G^*}(1,\dots,1)=\{(t_1,\dots,t_m)\in\R^m\mid \sum_{i\in [m]} t_i = n,\ 
\sum_{i\in I} t_i \geq z_I, \textrm{ for }I\subset [m]\}.
$$
Thus the polytope $P_{G^*}^-(1,\dots,1)$, which is the Minkowski
difference of the above polytope and $\Delta_{[m]}$, is given by
$$
P_{G^*}^{-}(1,\dots,1)=\{(t_1,\dots,t_m)\in\R^m\mid \sum_{i\in [m]} t_i = n-1,\ 
\sum_{i\in I} t_i \geq z_I, \textrm{ for }I\subset [m]\}.
$$
We have $z_I = 
\#\{j\in [n]\mid i\in I,\textrm{ for any edge }(i,\bar j)\in G\}
= n -\left|\bigcup_{j\in J} I_j\right|$,
for $I\subseteq [m]$ and $J=[m]\setminus I$.  Thus we 
can rewrite the inequality $\sum_{i\in I} t_i\geq z_I$ 
as $\sum_{j\in J}t_j \leq \left|\bigcup_{j\in J} I_j\right|-1$.
These are exactly the  inequalities from
the definition of $G$-draconian sequence, which proves the claim.
\end{proof}

This shows that Theorem~\ref{th:gen_ehrhrart} gives
a formula for the number of lattice points of the polytope
$P_G^-(y_1,\dots,y_m)$ as a sum over the lattice points
of $P_{G^*}^-(1,\dots,1)$, and vise versa.
In particular, we obtain the following duality 
for trimmed generalized permutohedra.

\begin{corollary}
\label{cor:duality-lattice-points}
 The number of lattice points in the polytope
$P_{G}^-(1,\dots,1)$ equals the number of lattice points 
in the polytope $P_{G^*}^-(1,\dots,1)$:
$$
\#(P_{G}^-(1,\dots,1)\cap \Z^n) = \#(P_{G^*}^-(1,\dots,1)\cap \Z^m ).
$$
\end{corollary}

Notice that the polytopes $P_{G}^-(1,\dots,1)$ and $P_{G^*}^-(1,\dots,1)$ have
different dimensions and they might be very different.   In
Theorem~\ref{th:Left-Right-degrees} we will describe a class of bijections
between lattice points of these polytopes.

\begin{example}
\label{examp:delta_times_delta}
Let $G=K_{m,n}$ be the complete bipartite graph.
Then $P_{K_{m,n}}^-$ is the $(n-1)$-dimensional simplex inflated
$m-1$ times: $P_{K_{m,n}}^- = (m-1)\Delta_{[n]}$.
The polytope for the mirror image of the graph 
is obtained by switching $m$ and $n$:
$P_{K_{m,n}^*}^- = (n-1)\Delta_{[m]}$.
Corollary~\ref{cor:duality-lattice-points} says that these two 
polytopes have the same number of lattice points.
This is a advanced way to say that 
$\binom{m+n-2}{m-1}  = \binom{m+n-2}{n-1}$.
\end{example}

Theorem~\ref{th:Todd-Euler-Maclaurin}(2)
from Appendix~B (Euler-MacLaurin formula for polytopes)
gives the following alternative expression for the generalized Ehrhart 
polynomial, i.e,
for the number of lattice points in $P_n^z(\{z_I\})$.
Without loss of generality, we will assume that $z_{[n]} = 0$.
The volume $\Vol P_n^z(\{z_I\})$ is a homogeneous polynomial 
$\tilde V_n$ in the $z_I$, for all nonempty $I\subsetneq [n]$.

\begin{proposition}
The number of lattice points in the generalized 
permutohedron $P_n^z(\{z_I\})$ is given by the polynomial
obtained from the polynomial $\tilde V_n$ by applying the 
Todd operator
$\Todd_n =
\prod_{I\subsetneq [n]} \Todd\left(-
\frac{\partial}{\partial z_I}\right)$,
where $\Todd(q) = q/(1-e^{-q}) = 1 + \frac{t}{2} + \frac{t^2}{12} - 
\frac{t^4}{720}+\cdots$.
\end{proposition}

\section{Root polytopes and their triangulations}
\label{sec:root_polytope}

\begin{definition}
For a graph $G$ on the vertex set $[n]$, let
$\tilde Q_G\subset \R^{n}$ be the convex hull of the origin $0$ and
the points $e_i-e_j$, for all edges $(i,j)$, $i<j$, of $G$.  
We will call polytopes $\tilde Q_G$ {\it root polytopes.}
In other words, a root polytope is the convex hull of the origin and 
some subset of end-points of positive roots for a root system of type $A_{n-1}$.
Polytopes $\tilde Q_G$ belong to an $n-1$ dimensional hyperplane.
\end{definition}

In the case of the complete graph $G=K_{n}$, the polytope 
$\tilde Q_{K_{n}}$ was studied in~\cite{GGP}.  In particular, we 
constructed a triangulation of this polytope and proved
that its $(n-1)$-dimensional volume equals $\frac{1}{(n-1)!} C_{n-1}$,
where $C_{n-1} = \frac{1}{n}\binom{2(n-1)}{n-1}$ is the $(n-1)$-st 
Catalan number.

In this section we study root polytopes for a bipartite graphs
$G\subseteq K_{m,n}$.  It is convenient to introduce 
related polytopes
$$
Q_G = \mathrm{ConvexHull}(e_i-e_{\bar j}\mid \textrm{for edges }(i,\bar j)
\textrm{ of } G) \subset \R^{m+n},
$$
where $e_1,\dots,e_m,e_{\bar 1},\dots,e_{\bar n}$ are the coordinate vectors
in $\R^{m+n}$.
Since $G$ is a bipartite graph, the polytope $Q_G$ belongs
to an $(m+n-2)$-dimensional affine subspace. 
The polytope $\tilde Q_G$ is the pyramid with the base $Q_G$ and the vertex $0$.
Thus $r! \,\mathrm{Vol}_{r}\, \tilde Q_G =
(r-1)!\, \mathrm{Vol}_{r-1}\, Q_G $, where 
$\mathrm{Vol}_r$ stands for the $r$-dimensional volume.
Slightly abusing notation, we will also refer to polytopes $Q_G$ as 
{\it root polytopes}.

The polytope $Q_{K_{m,n}}$ for the complete bipartite graph 
$K_{m,n}$ is the direct product of two 
simplices $\Delta^{m-1}\times \Delta^{n-1}$ of dimensions $(m-1)$
and $(n-1)$.  (Here $\Delta^{m-1}\simeq \Delta_{[m]}$.)
For other bipartite graphs, the polytope $Q_G$ is 
the convex hull of some subset of vertices of 
$\Delta^{m-1}\times \Delta^{n-1}$.
These polytopes are intimately related to generalized
permutohedra.

Let $I_1,\dots,I_m$ be the sequence of subsets associated
with the graph $G$, i.e., $j\in I_i$ whenever $(i,\bar j)\in G$.
Let $P_G =P_{G}(1,\dots,1) = \Delta_{I_1} + \cdots + \Delta_{I_m}$
and $P_G^-  = P_G - \Delta_{[n]}$.

\begin{theorem} 
\label{th:Vol_Q=N(P)}
For any connected bipartite graph $G\subseteq K_{m,n}$,
the $(m+n-2)$-dimensional volume of the root polytope $Q_{G}$ is 
expressed in terms of the number of lattice points of
the trimmed generalized permutohedron $P_G^-$ as
$$
\Vol Q_{G} = \frac{\# (P_G^{-}\cap \Z^n)}{(m+n-2)!}.
$$
\end{theorem}

We will prove this theorem by constructing a bijection between
simplices in a triangulation of the polytope $Q_G$
and lattice points of the polytope $P_G^-$; see
Theorem~\ref{th:Left-Right-degrees}.

For a bipartite graph $G\subseteq K_{m,n}$, let $G^+\subseteq K_{m+1,n}$
be the bipartite graph obtained from $G$ by adding a new vertex $m+1$
connected by the edges $(m+1,\bar j)$, $j=1,\dots,n$, with all vertices 
of the second part.  Then $P_{G^+}^- = P_G$.

\begin{corollary}
For any bipartite graph $G\subseteq K_{m,n}$ without isolated vertices,
the $(m+n-1)$-dimensional volume of the polytope $Q_{G^+}$ is
related to the number of lattice points in the generalized permutohedron
as
$$
\Vol Q_{G^+} = \frac{\# (P_G\cap \Z^n)}{(m+n-1)!}.
$$
\end{corollary}

\begin{definition}
A {\it polyhedral subdivision} of a polytope $Q$ is a subdivision of $Q$ into a
union of cells of the same dimension as $P$ such that each cell is the convex
hull of some subset of vertices of $Q$ and any two cells intersect properly,
i.e., the intersection of any two cells is their common face.  Polyhedral
subdivisions are partially ordered by refinement.  Minimal elements of this
partial order, i.e., unsubdividable polyhedral subdivisions, are called {\it
triangulations}.  In a triangulation each cell is a simplex.
\end{definition}

Triangulations of the product $\Delta^{m-1}\times \Delta^{n-1}$ were
first discussed by Gelfand-Kapranov-Zelevinsky~\cite[7.3.D]{GKZ}
and then studied by several authors; e.g., see Santos~\cite{San}.
We will analyze triangulations of more general root polytopes $Q_G$.
The following 3 lemmas were originally discovered circa 1992 
by the author in collaboration with Zelevinsky and Kapranov
in the context of triangulations of 
$\Delta^{m-1}\times \Delta^{n-1}$. 

Assume that the graph $G\subseteq K_{m,n}$ is connected.
First, let us describe the simplices inside the polytope $Q_G$.

\begin{lemma}
For a subgraph $H\subseteq G$, the convex hull of the 
collection $\{e_i-e_{\bar j}\mid (i,\bar j)\textrm{ is an edge of }H\}$
of vertices of $Q_G$ is a simplex if and only if 
$H$ is a forest in the graph $G$.
Such a simplex has maximal dimension $m+n-2$
if and only $H$ is a spanning tree of $G$.
All $(m+n-2)$-dimensional simplices of this form have 
the same volume $\frac{1}{(m+n-2)!}$.
\end{lemma}

\begin{proof} If $H$ contains a cycle $(i_1,\bar j_1),
(\bar j_1, i_2),(i_2,\bar j_2),\dots, (\bar j_k,i_1)$,
then the vectors $e_{i_1}- e_{\bar j_1}$, 
$e_{\bar j_1}- e_{i_2}$, \dots , $e_{\bar j_k}-e_{i_1}$ corresponding
to the edges in this cycle are linearly dependent.  (Their sum is zero.)
Thus the end-points of these vectors cannot be vertices of a simplex.
Conversely, for a forest, i.e., a graph without cycles, all vectors
are linearly independent and, thus, form a simplex.
\end{proof}

For a forest  $F\subseteq G$, we will denote the simplex from
this lemma by 
$$
\Delta_F: = \mathrm{ConvexHull}
(e_i-e_{\bar j}\mid (i,\bar j)\textrm{ is an edge of }F).
$$
A triangulation of $Q_G$ as a collection 
of simplices $\{\Delta_{T_1},\dots,\Delta_{T_s}\}$, for some spanning trees 
$T_1,\dots,T_s$ of $G$ such that
$Q_G= \cup \Delta_{T_i}$; and
each intersection $\Delta_{T_i}\cap \Delta_{T_j}$
is the common face of these two simplices.

Let us now describe pairs of simplices that intersect properly.
For two spanning trees $T$ and $T'$ of $G$, let $U(T,T')$ be the {\it
directed} graph with the edge set $\{(i,\bar j)\mid 
(i,\bar j)\in T\}\cup
\{(\bar j,i)\mid (i,\bar j)\in T'\}$, i.e., 
$U(T,T')$ is the union of edges $T$ and $T'$ 
with edges of $T$ oriented from left to right and edges of $T'$ oriented
from right to left.  An directed {\it cycle} is a sequence 
of directed edges $(i_1,i_2),(i_2,i_3),\dots,(i_{k-1},i_k),(i_k,i_1)$
such that all $i_1,\dots,i_k$ are distinct.

\begin{lemma} For two trees $T$ and $T'$, the intersection
$\Delta_T\cap \Delta_{T'}$ is a common face of the simplices
$\Delta_{T}$ and $\Delta_{T'}$  if and only if the directed
graph $U(T,T')$ has no directed cycles of length $\geq 4$.
\end{lemma}

\begin{proof}
Suppose that $U(T,T')$ has a directed cycle of length $\geq 4$.
Then the graphs $T$ and $T'$ have nonempty partial matching (i.e.,
subgraphs with disjoint edges) $M$ and $M'$ such that
(1) $M$ and $M'$ have no common edges;  and
(2) $M$ and $M'$ are matching on the same 
vertex set. Then both $M$ and $M'$ should have $k\geq 2$ edges.
Let $x=\frac {1}{k}\sum_{(i,\bar j)\in M} (e_i - e_{\bar j}) = 
\frac {1}{k} \sum_{(i,\bar j)\in M'} (e_i - e_{\bar j})$.  Thus 
$x\in \Delta_{T}\cap \Delta_{T'}$.
However, the minimal face of the simplex $\Delta_T$ that contains $x$ 
is $\Delta_M$ and the minimal face of $\Delta_{T'}$ that contains $x$ is
$\Delta_{M'}$.  Since $M\ne M'$, we have $\Delta_M \ne \Delta_{M'}$.
Thus the intersection of the simplices $\Delta_{T}$ and $\Delta_{T'}$ 
{\it is not\/} their common face. 

Conversely, assume that $U(T,T')$ has no directed cycles of length $\geq 4$.
Let $F=T\cap T'$ be the forest formed by the common edges of $T$
and $T'$.  Because $U(T,T')$ is acyclic, we can find a function 
$h:\{1,\dots,m,\bar 1,\dots,\bar n\}\to \R$ such that 
(1) $h$ is constant on connected components of the forest $F$;
and (2) for any directed edge $(a,b)\in U(T,T')$ that joins
two different connected components of $F$, we have $h(a)<h(b)$.
The second condition says that if $(a,b)=(i,\bar j)$ is an edge of $T$
then $h(i)<h(\bar j)$, and if $(a,b) = (\bar j, i)$ is an edge of $T'$
then $h(i)>h(\bar j)$. The function $h$ defines a linear form
$f_h$ on the space $\R^{m+n}$ with the coordinates $h(1),\dots,h(m),
h(\bar 1),\dots,h(\bar n)$ in the standard basis.  The above conditions
imply that
(1) for any vertex $x$ in the common face $\Delta_F$ of $\Delta_T$ and
$\Delta_{T'}$, we have $f_h(x) =0$,
(2) for any vertex $x\in \Delta_{T}\setminus \Delta_F$, we have $f_h(x)<0$;
and (3) for any vertex 
$x\in \Delta_{T'}\setminus \Delta_F$, we have $f_h(x)>0$.
In other words, the hyperplane $f_h(x) = 0$ intersects the simplices
$\Delta_T$ and $\Delta_{T'}$ at their common face and separates
the remaining vertices of these simplices.  This implies that 
$\Delta_{T}\cap \Delta_{T'} = \Delta_F$, as needed.
\end{proof}

\begin{definition}
For a spanning tree $T\in K_{m,n}$, let us define the 
{\it left degree vector} $LD=(d_1,\dots,d_m)$ and the 
{\it right degree vector} $RD=(d_{\bar 1},\dots,d_{\bar n})$,
where $d_i=\deg_i(T)-1$ and $d_{\bar j}=\deg_{\bar j}(T)-1$ are the 
degrees of the vertices $i$ and $\bar j$ in $T$ minus 1.
Note that $LD(T)$ and $RD(T)$ are nonnegative integer vectors
because all degrees of vertices in a tree are strictly positive.
\end{definition}

\begin{lemma}
\label{lem:different_RD_LD}
Let $\{\Delta_{T_1},\dots,\Delta_{T_s}\}$
be a triangulation of $Q_G$.
Then, for $i\ne j$, $T_i$ and $T_j$ have different
left degree vectors $LD(T_i)\ne LD(T_j)$ and different right degree
vectors $RD(T_i)\ne RD(T_j)$.
\end{lemma}

\begin{proof}
It is enough to prove that it is impossible to find two different 
spanning trees $T$ and $T'$ have have
same degrees in, say, the left part $\deg_i(T)=\deg_i(T')$, for $i=1,\dots,m$,
and such that the directed graph $U(T,T')$ has no directed cycles of length 
$\geq 4$.  Suppose that we found two such trees.  
Let $F$ be the forest formed by the common edges of $T$ and $T'$.
The directed graph $U(T,T')$ induces an acyclic directed graph
on connected components of $F$.  Because
of the acyclicity of this graph, we can find a minimal connected
component $C$ of $F$ such that no directed edge of 
$U(T,T')$ enters to any vertex of $C$ from outside of this component. 
Since $T\ne T'$, the 
component $C$ cannot include all vertices.
Thus some vertex $i$ of $C$ should be joined by an edge 
$(i,\bar j)\in T\setminus F$ with a vertex in some other component.
Since we assumed that $\deg_i(T) = \deg_i(T')$, there
is an edge $(i,\bar k)\in T'\setminus F$.  But this edge 
should be oriented as $(\bar k, i)$ in the graph $U(T,T')$,
i.e., it enters the vertex $i$ of $C$.  Contradiction.
\end{proof}

An alternative proof of Lemma~\ref{lem:different_RD_LD} follows
from Lemma~\ref{lem:shifts_right_degree} below.

For a bipartite graph $G\in K_{m,n}$, let $G^*\in K_{n,m}$
be the same graph with the left and right components switched, i.e., $G^*$ 
is the mirror image of $G$.  Recall that the trimmed
generalized permutohedron $P_G^-$ is the 
Minkowski difference of the generalized permutohedron $P_G$
and the simplex $\Delta_{[n]}$.

\begin{theorem}
\label{th:Left-Right-degrees} 
For any triangulation 
$\{\Delta_{T_1},\dots,\Delta_{T_s}\}$ of the root polytope $Q_G$,
the set of right degree vectors $\{RD(T_1),\dots,RD(T_s)\}$
is exactly the set of lattice points in the trimmed 
generalized permutohedron $P_{G}^-$
(without repetitions).
Similarly, the set of left degree vectors $\{LD(T_1),\dots,LD(T_s)\}$
is exactly the set of lattice points in the polytope $P_{G^*}^-$
for the mirror image of the graph $G$.
\end{theorem}

We will prove this theorem in Section~\ref{sec:subdivision}.
This theorem says that each triangulation
$\tau = \{\Delta_{T_1},\dots,\Delta_{T_s}\}$ of the root polytope $G_Q$ 
gives a bijection
$$\phi_\tau: \#(P_{G}^-\cap \Z^n) \to \#(P_{G^*}^-\cap \Z^m )
$$
between lattice points of the polytope $P_{G}^-$ and 
the lattice points of the polytope $P_{G^*}^-$ such that
$\phi_\tau:RD(T_i)\mapsto LD(T_i)$, for $i=1,\dots,s$.

It is interesting to investigate which properties of a triangulation
$\tau$ can be recovered from the bijection $\phi_\tau$.
Also it is interesting to intrinsically describe the class of
bijections associated with triangulations of $Q_G$.

\begin{example}   Suppose that $G=K_{m,n}$. 
Theorem~\ref{th:Left-Right-degrees} says that 
each triangulation of the product 
$\Delta^{m-1}\times\Delta^{n-1}$ of two simplices 
gives a bijection between lattice points two inflated
simplices $P_{K_{m,n}}^- = (m-1)\Delta^{n-1}$  and
$P_{K_{n,m}}^- = (n-1)\Delta^{m-1}$;
see Example~\ref{examp:delta_times_delta}.
\end{example}

Another instance of a similar phenomenon related to maximal minors
of matrices was investigated by Bernstein-Zelevinsky~\cite{BZ}.

\section{Root polytopes for non-bipartite graphs}
\label{sec:root_non_bipartite}

Let us show how to extend the above results to root polytopes $\tilde Q_G$
for a more general class of graphs $G$ that may not be bipartite.
Assume that $G$ is a connected graph on the vertex set $[n]$ that satisfies
the following condition:
$$
\textrm{For $i<j<k$, if $(i,j)$ and $(j,k)$ are edges of $G$, 
then $(i,k)$ is also an edge of $G$.}
$$
The polytope $\tilde Q_G$ is has the dimension $(n-1)$.
Let us say that a triangulation of the polytope $\tilde Q_G$ is 
{\it central} if any $(n-1)$-dimensional simplex in this triangulation
contains the origin $0$.

\begin{definition}
Let us say that a tree is {\it alternating} if there are no
$i<j<k$ such that $(i,j)$ and $(j,k)$ are edges in $T$.
Equivalently, labels in any path in an alternating tree $T$ should alternate.
\end{definition}

Alternating trees were first introduced in~\cite{GGP}
in order to describe triangulations of $\tilde Q_{K_n}$.
They also appeared in~\cite{Pos} and~\cite{PS1}.

For a spanning tree $T\subseteq G$, let 
$\tilde \Delta_T=\mathrm{ConvexHull}(0,e_i-e_j\mid (i,j)\in T, i<j)$.

\begin{lemma}
{\rm cf.~\cite{GGP}} \ 
A simplex $\tilde\Delta_T$ may appear in a central
triangulation of $\tilde Q_G$ if and only if $T$ is an alternating tree.
All these simplices have the same volume $\frac{1}{(n-1)!}$.
\end{lemma}

\begin{proof} Suppose that a tree $T$ is not alternating.
Let us find a pair of edges $(i,j)$ and $(j,k)$ in $T$ with $i<j<k$.
Let $T'$ be the tree obtained from $T$ by replacing the edge $(i,j)$
with $(i,k)$ and $T''$ be the tree obtained for $T$ by replacing
the edge $(j,k)$ with $(i,k)$.
Then two simplices $\tilde \Delta_{T'}$ and $\tilde \Delta_{T''}$ 
intersect at their common face.  Their union $\tilde \Delta_{T'}
\cup\tilde \Delta_{T''}$ properly contain the simplex $\tilde \Delta_T$.
Moreover, for neighborhood $B$ of the origin, 
$(\tilde \Delta_{T'} \cup\tilde \Delta_{T''})\cap B = \tilde\Delta_T \cap B$.
If the simplex $\tilde\Delta_T$ belongs to some central triangulation then 
with can replace it by the pair of simplices 
$\tilde\Delta_{T'}$ and $\tilde\Delta_{T'}$ and obtain a ``bigger''
triangulation, which is impossible.
\end{proof}

For an alternating tree $T$, we say that a vertex $i\in[n]$ is 
a {\it left vertex} if, for any edge $(i,j)$ in $T$, we have $i<j$.
Otherwise, if, for any edge $(i,j)$ in $T$, we have $i>j$, we say
that $i$ is a {\it right vertex}.
For a disjoint decomposition $[n]=L\cup R$, let $G_{L,R}$ be
the subgraph of $G$ given by 
$$
G_{L,R} = \{(i,j)\in G\mid i\in L, j\in R, i<j\}.
$$
The graph $G_{L,R}$ is a bipartite graph with the parts
$L$ and $R$.
Spanning trees of the graph $G_{L,R}$ are exactly alternating trees 
of $G$ with fixed sets $L$ and $R$ of left and right vertices.
Note that in general there are $2^{n-2}$ possible choices of the subsets 
$L$ and $R$ because we always have $1\in L$ and $n\in R$ and for any other 
vertex we have 2 options. However, some of these choices may lead
to disconnected graphs $G_{L,R}$ that contain no spanning trees.

Since each alternating tree in $G$ belongs to one of the graphs $G_{L,R}$,
we deduce that each simplex $\tilde \Delta_{T}$ in a central triangulation 
of $\tilde Q_G$ belongs to one of the polytopes $\tilde Q_{G_{L,R}}$.
Thus we obtain the following claim.

\begin{proposition}
\label{prop:LR-decomposes}
The polytope $\tilde Q_G$ decomposes into the union of polytopes
$\tilde Q_G = \bigcup_{L,R} \tilde Q_{G_{L,R}}$
over disjoint decompositions $[n]=L\cap R$ such
that the graph $G_{L,R}$ is connected.
Terms of this decompositions are in a bijection with 
the facets of $\tilde Q_G$ that do not contain
the origin. Each such facet $F$ has the form $F=Q_{G_{L,R}}$ and 
$\tilde Q_{G_{L,R}}$ is the pyramid with the base $F$.
Each central triangulation of $\tilde Q_G$ is obtained by selecting
a triangulation of each part $Q_{G_{L,R}}$.
\end{proposition}

Since each graph $G_{L,R}$ is bipartite, we can apply the 
results of this section and relate the volume of $\tilde Q_{G_{L,R}}$
to the number of lattice points in a certain (trimmed) generalized
permutohedron.  By Proposition~\ref{prop:LR-decomposes},
we can express the volume of the root polytope 
$\tilde Q_G$ as a sum of numbers of lattice points in several 
trimmed generalized permutohedra.

\begin{example}  In~\cite{GGP} we constructed a triangulation of the polytope
$\tilde Q_{K_n}$, for the complete graph $G=K_n$. 
This triangulation is formed by the simplices $\tilde\Delta_{T}$,
for all {\it noncrossing\/} alternating trees $T$, i.e., alternating
trees that contain no pair of crossing edges $(i,k)$ and $(j,l)$,
for $i<j<k<l$.  The number of such trees equals the $(n-1)$-st
Catalan number $C_{n-1}$.

For a disjoint decomposition $[n]=L\cap R$, let $K_{L,R}$ be the bipartite
graph with the edge set $\{(i,j)\mid i\in L, j\in R, i<j\}$.
According to Proposition~\ref{prop:LR-decomposes}, 
we have $\tilde Q_{K_n} = \bigcup_{L,R} \tilde Q_{K_{L,R}}$, where different
terms have no common interior points.  The collection of
simplices $\tilde \Delta_{T}$, for all noncrossing spanning trees
$T$ of the graph $K_{L,R}$, form a triangulation of the polytope 
$\tilde Q_{K_{L,R}}$.
\end{example}

This example and Theorem~\ref{th:Vol_Q=N(P)}
imply the following statement.

\begin{corollary}
For any disjoint decomposition $[n]=L\cap R$ such that $1\in L$
and $n\in R$, the number
of noncrossing spanning trees of the graph $K_{L,R}$ equals
the number of lattice points in the trimmed
generalized permutohedron $P_{K_{L,R}}^-$.
\end{corollary}

For example, if $L=\{1,\dots,l\}$ and $R=\{l+1,\dots,n\}$,
then $K_{L,R} = K_{l,n-l}$ is the complete bipartite graph. 
We deduce that the number of noncrossing trees in 
the complete bipartite graph $K_{l,n-l}$ equals the number of
lattice points in the polytope $P_{K_{l,n-l}}^- = (l-1)\Delta^{n-l-1}$,
which equals $\binom {n-2} {l-1}$.

\section{Mixed subdivisions of generalized permutohedra}
\label{sec:subdivision}

In this section we study mixed subdivisions of generalized 
permutohedra into parts isomorphic to direct products
of simplices.   For this we use the Cayley trick that relates 
mixed subdivisions of the Minkowski sum 
of several polytopes $P_1+\cdots + P_m$ to all polyhedral
subdivision of a certain polytope $\CC(P_1,\dots,P_m)$ of higher dimension.
The Cayley trick was first developed by Sturmfels~\cite{Stu}
for coherent subdivisions and by Humber-Rambau-Santos~\cite{HRS}
for arbitrary subdivisions. Santos~\cite{San} used this trick to study 
triangulations of the product of two simplices.

\begin{definition}
Let $d$ be the dimension of the Minkowski sum $P_1+\cdots + P_m$.
A  {\it Minkowski cell\/} in this Minkowski sum 
is a polytope $B_1+\cdots +B_m$ of the top dimension $d$, 
where each $B_i$ is a convex hull of
some subset of vertices of $P_i$.  A {\it mixed subdivision\/} of
the Minkowski sum is its decomposition into a union of several
Minkowski cells such that the intersection of any two cells is 
their common face.  Mixed subdivisions form a poset with respect
to refinement.  
A {\it fine mixed subdivision\/} is a minimal element in this poset.
\end{definition}

\begin{lemma}
\label{lem:fine_mixed_cells}
A mixed subdivision is fine if and only if, for each mixed
cell $B=B_1+\cdots + B_m$ in this subdivision, all
$B_i$ are simplices and $\sum\dim B_i =\dim B = d$.
\end{lemma}

\begin{proof} 
We leave this claim as an exercise,
or refer to~\cite[Proposition~2.3]{San}. 
\end{proof}

The mixed cells described in this lemma are called 
{\it fine mixed cells.}  The lemma implies that each
fine mixed cell $B_1+\cdots + B_m$ is isomorphic to 
the direct product $B_1\times \cdots\times B_m$
of simplices, i.e., the simplices $B_i$ span 
independent affine subspaces.  In order to emphasize this fact,
we will use the direct product notation
for fine cells.

Let $e_1,\dots,e_m,e_{\bar 1},\dots,e_{\bar n}$ be the standard
basis of $\R^{m+n}=\R^m\times \R^n$.  Embed the space $\R^n$, where
the polytopes $P_1,\dots,P_n$ live, into $\R^{m+n}$  
as the subspace with the basis $e_{\bar 1}, \dots,e_{\bar n}$.

\begin{definition}
Following Sturmfels~\cite{Stu} and Huber-Rambau-Santos~\cite{HRS}, 
we define the 
{\it Cayley embedding} of $P_1,\dots,P_m$
as the polytope $\CC(P_1,\dots,P_m)$ given by
$$
\CC(P_1,\dots,P_m) = \mathrm{ConvexHull}(e_i + P_i\mid i=1,\dots,m).
$$
\end{definition}

Let $(y_1,\dots,y_m)\times \R^n$ denote the $n$-dimensional affine 
subspace in $\R^{m+n}$ such that the first $m$ coordinates are
equal to some fixed parameters $y_1,\dots,y_m$.
(Here we think of the $y_i$ not as coordinates but as fixed parameters.)

\begin{lemma} {\rm \cite{Stu, HRS}} \
For any choice of parameters $y_1,\dots,y_m\geq 0$ such that $\sum y_i = 1$,
the intersection of $\CC(P_1,\dots,P_m)$ with the affine 
subspace $(y_1,\dots,y_m)\times \R^n$
is exactly the weighted Minkowski sum $y_1P_1+\cdots + y_m P_m$
(shifted into this affine subspace).
\end{lemma}

\begin{proof}
Indeed, by the definition, the polytope $\CC(P_1,\dots,P_m)$ 
is the locus of points
of the form $\sum_{i=1}^m \lambda_i(e_i + p_i)$, where $p_i\in P_i$,
$\lambda_i\geq 0$ and $\sum \lambda_i = 1$.
Intersecting a point of this form with $(y_1,\dots,y_n)\times \R^n$
means that we fix $\lambda_i = y_i$, for $i=1,\dots,m$.  This gives
the needed Minkowski sum.
\end{proof}

The next proposition expresses the {\it Cayley trick}.

\begin{proposition}
\label{prop:cayley_trick}
{\rm \cite{HRS}} \ 
Fix strictly positive parameters $y_1,\dots,y_m>0$
such that $\sum y_i = 1$.
For a polyhedral subdivision of $\CC(P_1,\dots,P_m)$,
intersecting its cells with $(y_1,\dots,y_n)\times R^n$
we obtain a mixed subdivision of $y_1 P_1 + \cdots + y_m P_m$.
This gives a poset isomorphism between polyhedral subdivisions of
$\CC(P_1,\dots,P_m)$ and mixed subdivisions of $y_1P_1+\cdots + y_m P_m$.
\end{proposition}

\begin{proof}
The first claim that a polyhedral subdivision of $\CC(P_1,\dots,P_m)$
gives a mixed subdivision of $y_1 P_1 + \cdots + y_m P_m$ is immediate.
On the other hand, we can recover a polyhedral subdivision
of $\CC(P_1,\dots,P_m)$ from a mixed subdivision 
of $y_1 P_1 + \cdots + y_m P_m$. 
We can always rescale cells of the mixed subdivision by changing values
of $y_1,\dots,y_m$ and obtain a mixed subdivision of
$y_1' P_1 + \cdots + y_m' P_m$, for any nonnegative $y_1',\dots,y_m'$.
As we vary $y=(y_1,\dots,y_m)$ over all points of the simplex 
$y_1,\dots,y_m\geq 0$, $y_1+\dots+y_m=1$, the unions 
$\cup_{y\in\Delta^{m-1}} yB$, for each mixed  cell $B$, 
form cells of the polyhedral subdivision of $\CC(P_1,\dots,P_m)$;
see~\cite{HRS} for details.
\end{proof}

Let $G\subseteq K_{m,n}$ be a connected bipartite graph.
Let $I_1,\dots,I_m\subseteq [n]$ be the associated collection of 
nonempty subsets: $I_i = \{j\mid (i,\bar j)\in G\}$, for $i=1,\dots,m$.
Then the Cayley embedding of the simplices $\Delta_{I_1},\dots,
\Delta_{I_m}$ is exactly
the root polytope $Q_G$ from Section~\ref{sec:root_polytope}:
$$
Q_G = \CC(\Delta_{I_1},\dots,\Delta_{I_m}).
$$
Recall that the generalized permutohedron $P_G(y_1,\dots,y_m)$ 
$$
P_G(y_1,\dots,y_m) = y_1\Delta_{I_1} + \cdots + y_m\Delta_{I_m},
$$
for the nonnegative $y_i$.
Proposition~\ref{prop:cayley_trick} specializes to 
the following claim.

\begin{corollary}  
\label{cor:cayley_trick_permutohedron}
For any strictly positive $y_1,\dots,y_m$, 
mixed subdivisions of the generalized permutohedron 
$P_G(y_1,\dots,y_m)$ are in one-to-one correspondence
with polyhedral subdivisions of the root polytope $Q_G$.
In particular, fine mixed subdivisions of $P_G(y_1,\dots,y_m)$
are in one-to-one correspondence with triangulations of $Q_G$.
This correspondence is given by intersecting a polyhedral 
subdivision of $Q_G$ with the subspace 
$(\frac{y_1}{s},\dots, \frac{y_m}{s})\times \R^n$, where $s=\sum y_i$,
and then inflating the intersection by the factor $s$.
\end{corollary}

In particular, this implies that the number of cells 
in a fine mixed subdivision of $P_G$ equals $(m+n-2)!\,\Vol Q_G$.

Let us describe fine mixed cells that appear in subdivisions
of $P_G(y_1,\dots,y_m)$.  For a sequence of nonempty subsets
$\J=(J_1,\dots,J_m)$,
let $G_{\J}$ be the graph with the edges
$(i,\bar j)$, for $j\in J_i$.

\begin{lemma}
\label{lem:fine_mixed_cells=trees}
Each fine mixed cell in a mixed subdivision of 
$P_G(y_1,\dots,y_m)$
has the form $y_1\Delta_{J_1}\times\cdots \times y_m\Delta_{J_m}$,
for some sequence of nonempty subsets $\J=(J_1,\dots,J_m)$ in $[n]$,
such that $G_{\J}$ is a spanning tree of $G$.
\end{lemma}

\begin{proof} By Lemma~\ref{lem:fine_mixed_cells}, each 
fine cell has the form $y_1\Delta_{J_1}\times \dots \times y_m\Delta_{J_m}$ 
where $J_i\subseteq I_i$, for $i=1,\dots,m$, i.e., $G_{\J}$
is a subgraph of $G$, the simplices $\Delta_{J_i}$ span independent 
affine subspaces, and $\sum \Delta_{J_i} = \sum (|J_i|-1) = n-1$.
This is equivalent to the condition that $G_{\J}$ is 
a tree.
\end{proof}

Let us denote the fine cell associated with a spanning tree 
$T\subseteq G$, as described in the above lemma, by
$$
\Pi_{T}:= y_{1} \Delta_{J_1}\times \cdots \times y_{m} \Delta_{J_m},
$$
where $J_i = \{j\mid (i,\bar j)\in T\}$, for $i=1,\dots,m$.
These fine cells $\Pi_T$ are exactly the cells
associated with the simplices $\Delta_T\subset Q_G$
from Section~\ref{sec:root_polytope} via the Cayley trick:
$$
\Pi_T =  s \left(\Delta_T\cap \left(\frac{y_1}{s},\dots,\frac{y_m}{s}\right)\times \R^n\right),
$$
where $s=\sum y_i$.  So it is not surprising that the fine cells
$\Pi_T$ are labeled by the same objects---spanning trees of $G$.

Let us explain the meaning of the left degree vector
$LD(T)=(d_1,\dots,d_m)$ and the right degree vector 
$RD(T)=(d_{\bar 1},\dots,d_{\bar n})$ of a tree $T\subseteq G$ 
in terms of the fine cell $\Pi_T$.   

\begin{lemma}
\label{lem:left_degree}
Let $LD(T)=(d_{1},\dots,d_{m})$ be the left
degree vector of a tree $T$, then
$$
\Vol \Pi_T = \frac{y_1^{d_{1}}}{d_{1}!}\cdots
\frac{y_m^{d_{m}}}{d_{m}!}.
$$
\end{lemma}

\begin{proof}
Indeed, $d_{\bar i} = |J_i|-1 = \dim \Delta_{J_i}$, for $i=1,\dots,m$.
\end{proof}

\begin{lemma}  
\label{lem:shifts_right_degree}
Let us specialize $y_1=\cdots = y_m = 1$.
For a spanning tree $T\subseteq G$, the fine cell 
$\Pi_T$ contains the shift 
$(a_1,\dots,a_n) + \Delta_{[n]}$ of the simplex $\Delta_{[n]}$
by an integer vector $(a_1,\dots,a_n)\in\Z^n$ if and only if
$(a_1,\dots,a_n)$ is the right degree vector $RD(T)$ of the tree $T$.
Moreover, if $(a_1,\dots,a_n)\in\Z^n$ is not the right degree vector of $T$,
then the shift $(a_1,\dots,a_n) + \Delta_{[n]}$ has no common
interior points with the cell $\Pi_T$.
\end{lemma}

\begin{proof}
Notice that, for two subsets $I,J\subseteq[n]$ with a nonempty 
intersection, we have the following inclusion of Minkowski sums:
$$
\Delta_I+\Delta_J\supseteq 
\Delta_{I\cup J} + \Delta_{I\cap J}.
$$
Indeed, the polytope $\Delta_{I\cup J} + \Delta_{I\cap J}$ is the convex hull
of all possible sums $e_i+e_j$, where $e_i$ is a vertex of 
$\Delta_{I\cup J}$ and 
$e_j$ a vertex of $\Delta_{I\cap J}$, i.e., $i\in I\cup J$ and $j\in I\cap J$.
We have either ($i\in I$ and $j\in J$), or ($i\in J$ and $j\in I$), or both.
In all cases, we have $e_i+e_j\in \Delta_I+\Delta_J$.

For the fine cell $\Pi_T = \Delta_{J_1}\times \cdots \times \Delta_{J_m}
=\Delta_{J_1}+\cdots +\Delta_{J_m}$,
pick two summands $\Delta_{J_i}$ and $\Delta_{J_j}$ with a nonempty 
intersection $J_i\cap J_j$ (what should contain exactly one element $k$)
and replace them by $\Delta_{J_i\cup J_j}$ and 
$\Delta_{J_i\cap J_j}$.  We obtain another cell $\Pi_{T'}\subseteq \Pi_T$,
where the tree $T'$ is obtained from $T$ by replacing all edges 
$(j,\bar l)\in T$, for $l\ne k$, with the edges $(i,\bar l)$.
Notice that the tree $T'$ has exactly the same right degree vector 
$RD(T')=RD(T)$.  Let us keep repeating this operation until we
obtain a cell of the form $\Pi_{T''} = 
\Delta_{\{i_1\}} + \dots + \Delta_{\{i_m\}} +  \Delta_{[n]}  \subseteq \Pi_T$,
i.e., all summands are single vertices except for one simplex $\Delta_{[n]}$.
Since the tree $T''$ has the same right degree vector 
$(d_{\bar 1},\dots,d_{\bar n})=RD(T'')=RD(T)$ as the tree $T$,
we deduce that $\#\{j\mid i_j=i\} = d_{\bar i}$, for $i=1,\dots,n$.
In other words, $\Pi_{T''} = 
(d_{\bar 1},\dots,d_{\bar n}) + \Delta_{[n]}\subseteq \Pi_F$.

It remains to show that any other shift 
$(a_1,\dots,a_n)+\Delta_{[n]}$, for an integer vector 
$(a_1,\dots,a_n)\ne 
(d_{\bar 1},\dots,d_{\bar n})$, has no common interior points with
the cell $\Pi_T$.  Suppose that there exists such a shift
with a common interior point $b\in \Pi_T\cap ((a_1,\dots,a_n)+\Delta_{[n]})$.
Let $r=(d_{\bar 1}-a_1,\dots,d_{\bar n}-a_n)\in
\Z^{n}\setminus \{(0,\dots,0)\}$.  Then the point $b+r$ is an interior
point of $(d_{\bar 1}, \dots, d_{\bar n}) + \Delta_{[n]}\subseteq \Pi_T$.
Thus the whole line interval $[b,b+r]$ belong to the interior
of the fine cell $\Pi_F=\Delta_{J_1}\times\cdots \times \Delta_{J_m}$.  
Here $b\in \R^n$ and $r$ is a nonzero integer vector.
Thus at least one projection $[b',b'+r']$ of the interval $[b,b+r]$
to some component $\Delta_{J_i}$ of the direct product has a nonzero length.
Here $r'$ is should be a nonzero integer vector and $[b',b'+r']$ 
should belong to the interior of the simplex $\Delta_{J_i}$.
But this is impossible.  No coordinate simplex can contain such an 
interval strictly in its interior.  Indeed, the diameter of a
coordinate simplex in the usual Euclidean metric on $\R^n$ is 
$2^{\frac{1}{n}}$.  The only integer vectors
that have smaller length are the coordinate vectors $\pm e_j$.
If $b'$ belongs to a coordinate simplex then $b'\pm e_j$
does not belong to it,  because the vector $\pm e_j$ 
does not lie in the hyperplane
where all coordinate simplices live.  We obtain a contradiction.
\end{proof}

Let us now prove Theorems~\ref{th:Left-Right-degrees}
and~\ref{th:gen_ehrhrart}.

\begin{proof}[Proof of Theorem~\ref{th:Left-Right-degrees}] 
It is enough to prove the statement about right degree vectors
and deduce the statement about left degree vectors by symmetry.
By Corollary~\ref{cor:cayley_trick_permutohedron},
simplices in a triangulation $\{\Delta_{T_1},\dots,\Delta_{T_s}\}$
of the root polytope $Q_G$ are in one-to-one correspondence with cells 
in the corresponding fine mixed subdivision $\{\Pi_{T_1},\dots,\Pi_{T_s}\}$
of the generalized permutohedron $P_G$.
By Lemma~\ref{lem:shifts_right_degree}, each cell $\Pi_{T_i}$
contains the shifted simplex $a+\Delta_{[n]}$,
where $a=RD(T_i)$, and each integer shift 
$a+\Delta_{[n]}\subseteq P_G$
belongs to one of the cells $\Pi_{T_i}$.   Notice that
the set of integer vectors $a\in\Z^n$ such that 
$a+\Delta_{[n]}\subseteq P_G$ is exactly the set of lattice points
of the trimmed generalized permutohedron $P_G^-$.
This proves that the map $\Delta_{T_i}\mapsto RD(T_i)$ is a bijection
between simplices in the triangulations and lattice points of $P_G^{-}$,
as needed.
\end{proof}

\begin{proof}[Proof of Theorem~\ref{th:gen_ehrhrart}] 
Let us fix a fine mixed subdivision
$\{\Pi_{T_1},\dots,\Pi_{T_m}\}$ of the polytope $P_G(y_1,\dots,y_m)$.
According to Lemma~\ref{lem:left_degree}, the volume of 
$P_G(y_1,\dots,y_m)$ can be written as
$$
\Vol P_G(y_1,\dots,y_m) = \sum_{i=1}^m \frac {y_1^{d_1(T)}}{d_1!}\cdots
\frac{y_m^{d_m(T)}}{d_m!}.
$$
Let us compare this expression with the expression for 
$\Vol P_G(y_1,\dots,y_m)$ given by 
Theorem~\ref{th:second-formula-generalized}.
We deduce that the map $\Pi(T_i)\mapsto LD(T_i)$ is a bijection
between fine cells $\Pi_{T_i}$ in this subdivision and $G$-draconian sequences.
According to the Cayley trick and Theorem~\ref{th:Left-Right-degrees}, 
the number of fine cells in this subdivision equals the number of simplices 
in a triangulation for $Q_G$ equals the number of lattice points
in $P_{G^{-}}(1,\dots,1)$.
We deduce that the number of $G$-draconian sequences equals 
the number of lattice points of $P_{G^-}(1,\dots,1)$.
This is exactly the claim of Theorem~\ref{th:gen_ehrhrart}
in the case when $y_1=\dots=y_m=1$.

The case of general $y_1,\dots,y_m$ follows from this special
case.  Indeed, we can write any weighted Minkowski sum $y_1\Delta_{I_1} +
\cdots + y_m\Delta_{I_m}$, for nonnegative integers $y_1,\dots,y_m$, as the
Minkowski sum of $y_1$ copies of $\Delta_{I_1}$, $y_2$ copies of
$\Delta_{I_2}$, etc.  When we do this transformation 
the right-hand sides of expressions given by Theorem~\ref{th:gen_ehrhrart} 
agree.  For example, if we replace the term $y_1 \Delta_{I_1}$ 
in the Minkowski sum with the sum $z_1\Delta_{I_1} + z_2\Delta_{I_1}$,
where $y_1=z_1+z_2$,
then the can correspondingly modify the right-hand side using the identity
$\frac {(y_{1})_{a_1}}{a_1!}
= \binom{y_{1}+a_1-1}{a_1} = 
\sum_{b_{1}+b_2=a_1} \frac {(z_{1})_{b_1}}{b_1!}\,
\frac {(z_2)_{b_2}}{b_2!}$.
\end{proof}

\begin{remark}  We can also deduce that the number of
$G$-draconian sequences equals $(m+n-2)!\,\Vol Q_G$,
i.e., the number of simplices in a triangulation of $Q_G$,
using integration.
Let us calculate the volume $\Vol Q_G$ by integrating
the volume of its slice $P_G(y_1,\dots,y_m)$ given by
Theorem~\ref{th:second-formula-generalized} over all points
of the $(m-1)$-dimensional simplex $\Delta_{[m]}$:
$$
\Vol Q_G = \int_{(y_1,\dots,y_m)\in\Delta_{[m]}}\,\,
\Vol P_G(y_1,\dots,y_m) dy_1\cdots dy_{m-1}.
$$
Now we can use the fact that the integral 
of a monomial $\frac{y_1^{a_1}}{a_1!}\cdots \frac{y_m^{a_m}}{a_m!}$
over the simplex $\Delta_{[m]}$ equals $((m-1+\sum a_i)!)^{-1}$.
\end{remark}

Also remark that the first part of the above proof 
and Theorem~\ref{th:Left-Right-degrees} gives 
an alternative proof of Lemma~\ref{lem:G-draconain=G*} 
saying that the set of $G$-draconian sequences is
the set of lattice points in $P_{G^*}^-(1,\dots,1)$.

\begin{example}
Let us assume that $I_1,\dots,I_m$, $m=2^n-1$, are all 
nonempty subsets of $[n]$ and $G$ is the associated
bipartite graph.
The $G$-draconian sequences of integers are in one-to-one 
correspondence with all {\it unordered\/} collections of subsets 
in $[n]$ satisfying the dragon marriage condition.  
For a draconian sequence $(a_1,\dots,a_m)$ there
are $\binom{n-1}{a_1,\dots,a_m}$ associated ordered
sequences of subsets.
In this case, $P_G = P_n(2^{n-1},2^{n-2}\dots,2,1)$
and $P_G^- = P_n(2^{n-1}-1,2^{n-2}\dots,2,1)$ (both are
usual permutohedra).
The number of draconian sequences is exactly the
number of lattice points in the permutohedron
$P_n(2^{n-1}-1,2^{n-2}\dots,2,1)$.
\end{example}

Another approach to counting lattice points in generalized
permutohedra is based on constructing its fine mixed subdivision
and paying  a special attention to lower dimensional cells.
Let us say that a {\it semi-polytope} is a bounded subset of points 
in a real vector space given by a finite collection of affine weak 
and strict equalities.  Define coordinate {\it semi-simplices} as
$$
\Delta_{I,j}^{semi} =\Delta_I\setminus \Delta_{I\setminus\{j\}} = 
\left\{\sum_{i\in I} x_{i}\,e_{i} \mid \sum_{i\in I} x_{i} =1;\
x_{i}\geq 0,\textrm{ for $i\in I$; and } x_j>0\right\},
$$
for $j\in I\subseteq[n]$.

\begin{proof} [Alternative semiproof of Theorem~\ref{th:gen_ehrhrart}] 
Let $P_G(y_1,\dots,y_m)=y_1\Delta_{I_1}+\cdots +y_m\Delta_{I_m}$.
Assume that $I_1=[n]$.
It seems feasible that there exists a {\it disjoint}
decomposition of the polytope $P_G(y_1,\dots,y_m)$ into 
semipolytopes of the form
\begin{equation}
\label{eq:disjoint_decom_PG}
P_G(y_1,\dots,y_m) = \bigcup_{(J_1,\dots,J_m)} y_{1} \Delta_{J_1}\times 
y_{I_2} \Delta_{J_2,j_2}^{semi} \cdots \times y_{m} \Delta_{J_m,j_m}^{semi},
\end{equation}
where the sum is over sequences of subsets $(J_1,\dots,J_m)$ and
$j_2,\dots,j_m$ such that $j_i\in J_i\subseteq I_i$, 
and bipartite graphs associated with $(J_1,\dots,J_m)$ are
spanning trees $T$ of $G$. In particular, the closure of each term
is a fine mixed cell $\Pi_T$ of top dimension.

Here is a not quite rigorous reason why this should be true.
Let us start with the top dimensional simplex $y_1\Delta_{I_1}$,
$I_1=[n]$.  When we add the simplex $y_2\Delta_{I_2}$, we
create several new fine cells.  Each of these cells is the direct product 
$y_1\Delta_{J_1}\times y_2 \Delta_{J_2}$ of a face of $y_1\Delta_{I_1}$ and
a face of $y_2\Delta_{I_2}$ glued to $y_1\Delta_{I_1}$ 
by one if its facets $y_1\Delta_{J_1}\times y_2\Delta_{J_2\setminus\{j_2\}}$.
This is why we exclude elements of this facet.  When we add 
$y_3\Delta_{I_3}$ we again create several new fine cells.   
Again each of these new 
cells is a direct product of one of the faces of the polytope
created on the earlier stage and a face $y_3\Delta_{J_3}$ of $y_3\Delta_{I_3}$.
Again each of these cells should be glued by a facet of
$y_3\Delta_{J_3}$, etc.

Let us show that just an existence of a decomposition for the 
form~(\ref{eq:disjoint_decom_PG}) already implies 
Theorem~\ref{th:gen_ehrhrart}.  Indeed, the number of lattice points
in one of the terms of this decomposition equals
$\frac{(y_1+1)_{a_1}}{a_1!} \frac{(y_2)_{a_2}}{a_2!}\cdots
\frac{(y_m)_{a_m}}{a_m!}$ and its volume is 
$\frac{{(y_1+1)}^{a_1}}{a_1!} \frac{y_2^{a_2}}{a_2!}\cdots
\frac{y_m^{a_m}}{a_m!}$, where $a_i = \dim \Delta_{J_i} = |J_i|-1$.
Thus the formula for the number of lattice points in
$P_G(y_1,\dots,y_m)$ is obtained from the formula for 
the volume given by
Theorem~\ref{th:second-formula-generalized} by
replacing usual powers with raising powers, as needed.
\end{proof}

In order to make this proof more rigorous, we need to 
carefully analyze all possible cases.  Preferably one would
like to have an explicit construction for a decomposition of the 
form~(\ref{eq:disjoint_decom_PG}).

In Section~\ref{sec:shifted_tableaux}, we will need the following
statement.

\begin{proposition}
\label{prop:lattice=sum_vertices}
Any integer lattice point of the generalized
permutohedron $P_G= \Delta_{I_1}+\cdots +\Delta_{I_m}$
has the form $e_{j_1}+\cdots + e_{j_m}$, where $j_k\in I_k$,
for $k=1,\dots,m$.
\end{proposition}

\begin{remark} Proposition~\ref{prop:lattice=sum_vertices}
says that any lattice point of the generalized permutohedron is the sum
of vertices of its Minkowski summand.
Note that of a similar clam is not true for an arbitrary Minkowski sum.
For example, the Minkowski sum of two line intervals $[(0,1),(1,0)]$
and $[(0,0),(1,1)]$ contains the lattice point $(1,1)$ which cannot
be presented as a sum of the vertices.
\end{remark}

\begin{proof}[Proof of Proposition~\ref{prop:lattice=sum_vertices}]
Each lattice point of $P_G$ belongs to a fine mixed cell in
a fine mixed subdivision of $P_G$; see Section~\ref{sec:subdivision}.
According to Lemma~\ref{lem:fine_mixed_cells=trees}, each fine mixed
cell is a direct product $\Delta_{J_1}\times \cdots \times \Delta_{J_m}$
of simplices, where $J_i\subseteq I_i$, for $i=1,\dots,m$,
and the graph $T=G_{(J_1,\dots,J_m)}\subseteq K_{m,n}$ is a bipartite tree.
Any lattice point $(b_1,\dots,b_n)$ of
$\Delta_{J_1}\times \cdots \times \Delta_{J_m}$ comes from a function
$f:\{(i,\bar j)\}\to\R_{\geq 0}$ defined on edges of the tree $T$ such that
(1) $f(i,\bar j)\geq 0$, (2) $\sum_j f(i,\bar j) = 1$, and
(3) $\sum_i f(i,\bar j) = b_j$, for any $i=1,\dots,m$ and $j=1,\dots,n$.
Since $T$ is a tree and the sum of values of $f$ over edges at any node
of $T$ is
integer, we deduce that $f$ has all nonnegative integer values.
(First, we prove this for leaves of $T$, then for leaves of the tree
obtained by removing the leaves of $T$, etc.)
Thus, for any $i=1,\dots,m$, we have $f(i,\bar j_i)=1$, for some $j_i$,
and $f(i,\bar j)=0$, for $j\ne j_i$.
Thus $(b_1,\dots,b_n) = e_{j_1}+\cdots + e_{j_n}$, as needed.

\end{proof}

\section{Application: diagonals of shifted Young tableaux}
\label{sec:shifted_tableaux}

A {\it standard shifted Young tableaux\/} of the triangular shape
$(n,n-1,\dots,1)$ is a bijective map
$T:\{(i,j)\mid 1\leq i\leq j\leq n\} \to\{1,\dots, \binom{n+1}{2}\}$
increasing in the rows and the columns, i.e., $T((i,j))<T((i+1,j))$ and
$T((i,j))<T((i,j+1))$, whenever the entries are defined.
Let us say that the {\it diagonal vector\/} of such a tableau $T$
is the vector $\diag(T)=(d_1,\dots,d_n):=(T(1,1),T(2,2),\dots,T(n,n))$;
see Example~\ref{exam:shifted-tableau} below.
Is clear, that $d_1=1$, $d_n=\binom{n+1}{2}$, and $d_{i+1}>d_i$.
In this section we describe all possible diagonal vectors.

For a nonnegative integer $(n-1)$-vector $(a_1,\dots,a_{n-1})$,
let $N(a_1,\dots,a_{n-1})$ be the number of standard shifted
Young tableaux $T$ of the triangular shape with the diagonal vector
$\diag(T)=(1,a_1+2,a_1+a_2+3,\dots,a_1+\cdots+a_{n-1} + n)$,
or, equivalently, $a_i = d_{i+1}-d_i-1$, for $i=1,\dots,n-1$.

\begin{theorem}
\label{th:shifted=productij}
We have the following identity:
$$
\sum_{a_1,\dots,a_{n-1}\geq 0}
N(a_1,\dots,a_n)\, \frac{t_1 ^{a_1}} {a_1!}\cdots
\frac{t_{n-1} ^{a_{n-1}}} {a_{n-1}!} =
\prod_{1\leq i<j\leq n} \frac{t_i+t_{i+1}+\cdots + t_{j-1}}{j-i}.
$$
\end{theorem}

\begin{proof}
Let $\lambda=(\lambda_1\geq \dots\geq \lambda_n)$ be a partition.
The {\it Gelfand-Tsetlin polytope\/} $GT(\lambda)$ is defined
as the set of triangular arrays
$(p_{ij})_{i,j\geq 1, i+j\leq n}\in \R^{\binom{n+1}{2}}$
such that the first row is $(p_{11},p_{12},\dots,p_{1n}) = \lambda$
and entries in consecutive rows are interlaced
$p_{i1}\geq p_{i+1\, 1}\geq p_{i2}\geq p_{i+1\,2}\geq \cdots$,
for $i=1,\dots,n-1$.

Let us calculate the volume of the polytope $GT(\lambda)$ in
two different ways.   First, recall that lattice points of $GT(\lambda)$
correspond to elements of the Gelfand-Tsetlin basis of the irreducible
representation $V_\lambda$ of $GL(n)$ with the highest weight $\lambda$.
Thus the number of the lattice points is given by the Weyl dimension
formula:
$\#(GT(\lambda)\cap \Z^{\binom{n+1}{2}}) = \prod_{1\leq i<j\leq n}
\frac{\lambda_i-\lambda_j + j-i}{j-i}$.
We deduce that the volume of $GT(\lambda)$ is given by the top homogeneous
component of this polynomial in $\lambda_1,\dots,\lambda_n$:
$$
\Vol GT(\lambda) = \prod_{1\leq i<j\leq n} \frac{\lambda_i-\lambda_j}{j-i}.
$$
On the other hand, note that the shape of an array
$(p_{ij})\in GT(\lambda)$ is equivalent to the shape of
a shifted tableau.
Let us subdivide $GT(\lambda)$ into parts by the hyperplanes
$p_{ij}=p_{kl}$, for all $i,j,k,l$. A region of this
subdivision of the Gelfand-Tsetlin polytopes $GT(\lambda)$
correspond to a choice of a total ordering of the $p_{ij}$ compatible
with all inequalities.  Such ordering are in one-to-one correspondence
with standard shifted Young tableaux of the triangular shape
$(n,n-1,\dots,1)$.
For a tableau $T$ with the diagonal vector $\diag(T)=(d_1,\dots,d_n)$,
the associated region of $GT(\lambda)$ is isomorphic to
$\{(y_1<\dots < y_{\binom{n+1}{2}})\mid y_{d_i} = \lambda_{i},\textrm{ for }
i=1,\dots,n\}$, 
that is, to the direct product of simplices
$(\lambda_1-\lambda_2)\Delta^{d_{2}-d_1-1}\times \cdots \times
(\lambda_{n-1}-\lambda_n)\Delta^{d_{n}-d_{n-1}-1}$.
The volume of this direct product equals
$$
\prod_{i=1}^{n-1} \frac {(\lambda_i-\lambda_{i+1})^{d_{i+1}-d_i-1}}
{(d_{i+1}-d_i-1)!}.
$$
Thus the volume of $GT(\lambda)$ can be written as the sum of these
expressions over standard shifted tableaux.
Comparing these two expressions for $\Vol GT(\lambda)$ and
writing them in the coordinates $t_i = \lambda_i-\lambda_{i+1}$,
we obtain the needed identity.
\end{proof}

Theorem~\ref{th:shifted=productij} implies that $N(a_1,\dots,a_n)$ can 
be nonzero only if $(a_1,\dots,a_n)$ is a lattice point 
of the Newton polytope 
$$
\Ass_{n-1}:=\Newton\left(\prod_{1\leq i<j\leq n}
(t_i+t_{i+1}+\cdots +t_{j-1})\right) = 
\sum_{1\leq i<j\leq n}\Delta_{[i,j-1]}.
$$
This Newton polytope 
is exactly the associahedron in the Loday realization,
for $n-1$; see Subsection~\ref{ssec:associahedron}.
Using Proposition~\ref{prop:lattice=sum_vertices}, we
obtain the following statement.

\begin{corollary}  The number of different diagonal vectors
in standard shifted Young tableaux of the shape $(n,n-1,\dots,1)$
is exactly the number of integer lattice points in the associahedron
$\Ass_{n-1}$.  More precisely, $N(a_1,\dots,a_{n-1})$ is nonzero
if and only if $(a_1,\dots,a_{n-1})$ is an integer lattice point
of $\Ass_{n-1}$.
\end{corollary}

It would be interesting to extend this claim to other 
shifted shapes.

\begin{example}
\label{exam:lattice_points_ass}
Let $D_n$ be the number of different diagonal vectors, or, equivalently,
the number integer lattice points in $\Ass_{n-1}$, 
or, equivalently, the number of nonzero monomials
in the expansion of the product $\prod_{1\leq i< j\leq n} 
\sum_{k=i}^{j-1} t_k$.
Several numbers $D_n$ are given below.

\smallskip
\begin{center}
\begin{tabular}{|c||l|l|l|l|l|l|l|l|l|}
\hline
$n$   & 1 & 2 & 3 & 4 & 5   & 6    & 7      & 8        & 9          \\
\hline 
$D_n$ & 1 & 1 & 2 & 8 & 55  & 567  & 7958   & 142396   & 3104160     \\ 
\hline
\end{tabular}
\end{center}
\end{example}

Theorem~\ref{th:shifted=productij} also implies that $N(a_1,\dots,a_n)$  
equals $\prod_{i=1}^{n-1} (a_i)!/(1!2!\cdots (n-1)!)$ times the number
of ways to write the point $(a_1,\dots,a_{n-1})$ as a sum of vertices
of the simplices $\Delta_{[i,j-1]}$.  In particular, if
$(a_1,\dots,a_{n-1})$ is a vertex of the associahedron 
$\Ass_{n-1}$ then the second factor is $1$.

Recall that vertices of $\Ass_{n-1}$ correspond to plane 
binary trees on $n-1$ nodes; see Subsection~\ref{ssec:associahedron}.
For a plane binary tree on $n-1$ nodes, let $L_i, R_i$, $i=1,\dots,n-1$,
be the left and right branches of the nodes arranged in the binary
search order; see Subsection~\ref{ssec:associahedron}.
Also let $l_i = |L_i|+1$ and $r_i = |R_i|+1$.

\begin{corollary}
The numbers of standard shifted Young tableaux with diagonal vectors
corresponding to the vertices of the associahedron are given by 
$$
T(l_1\cdot r_1,\dots,l_{n-1}\cdot r_{n-1}) = \frac{(l_1\cdot r_1)!\cdots
(l_{n-1}\cdot r_{n-1})!}{1!\,2!\cdots (n-1)!} = 
f_{l_1\times r_1}\cdots f_{l_{n-1}\times r_{n-1}},
$$
where $f_{k\times l}$ is the number of standard Young tableaux
of the rectangular shape $k\times l$.
\end{corollary}

The second expression can be obtained from the first using
the hook-length formula for the number of standard Young tableaux.
We can also deduce it directly, as follows.
Recall that binary trees on $n-1$ nodes are associated with 
subdivisions of the shifted shape $(n-1,n-2,\dots,1)$
into $n-1$ rectangles of sizes $l_1\times r_1$, \dots,
$l_{n-1}\times r_{n-1}$; see Subsection~\ref{ssec:associahedron}.
Each shifted tableaux with the diagonal vector
$(d_1,\dots,d_n) = 
(1,\ 2 + l_1\cdot r_1,\ 3 + l_1\cdot r_1 + l_2 \cdot r_2,\ \cdots)$
is obtained from such a subdivision by adding $n$ diagonal boxes filled
with the numbers $d_1,\dots,d_n$ and filling the $i$-th rectangle
$l_i\times r_i$ with the numbers $d_{i}+1,d_{i}+2,\dots,d_{i+1}-1$
so that they from a rectangular standard tableau,
for $i=1,\dots,n-1$.

\begin{example}
\label{exam:shifted-tableau} 
The diagonal vector
$(1,3,10,12,15,36,40,43,45)$ is associated with 
the plane binary tree and the subdivision into rectangles
from Example~\ref{ex:plane_binary_trees}.
Here is a shifted tableau with the this diagonal vector 
obtained by filling the rectangles of this subdivision: 
\medskip

\begin{center}
{\tiny
\input{fig15-1.pstex_t}
}
\end{center}
\end{example}

\section{Mixed Eulerian numbers}
\label{sec:Mixed_Eulerian}

Let us return to the usual 
permutohedron $P_{n+1} = P_{n+1}(x_1,\dots,x_{n+1})$.
Let us use the coordinates $u_1,\dots,u_n$ related to $x_1,\dots,x_{n+1}$ by
$$
u_1= x_1-x_2,\
u_2 = x_2 - x_3,\
\cdots,\ 
u_n = x_n - x_{n+1}
$$
This coordinate system is canonically defined for an arbitrary Weyl group
as the coordinate system in the weight space given by the fundamental weights.

The permutohedron $P_{n+1}$ can be written as the Minkowski sum
$$
P_{n+1} = u_1\,\Delta_{1,n+1} + u_2\,\Delta_{2,n}+ \cdots + u_n \,\Delta_{n,n+1}
$$
of the {\it hypersimplices} $\Delta_{k,n+1} := P_{n+1}(1,\dots,1,0,\dots,0)$ 
with $k$ `$1$'s.
For example, the hexagon can be expressed as the Minkowski sum 
of the hypersimplices $\Delta_{1,3}$ and $\Delta_{2,3}$, 
which are two triangles with opposite orientations:
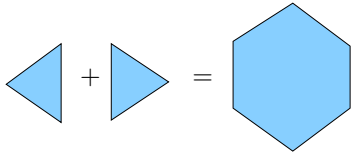
\begin{center}
\input{fig5-1.pstex_t}
\end{center}

According to Proposition~\ref{prop:mixed_volume},
the volume of $P_{n+1}$ can be written as
$$
\Vol P_{n+1} = 
\sum_{c_1,\dots,c_n}
{A_{c_1,\dots,c_n}}\, \frac{u_1^{c_1}}{c_1!}\cdots \frac{u_n^{c_n}}{c_n!},
$$
where the sum is over $c_1,\dots,c_n\geq 0$, $c_1+\cdots+c_n=n$, and
$$
{A_{c_1,\dots,c_n}} = n!\,V(\Delta_{1,n+1}^{c_1},\dots,
\Delta_{n,n+1}^{c_n})\in \Z_{> 0}
$$
is the mixed volume of hypersimplices multiplied by $n!$.
Here $P^l$ means the polytope $P$ repeated $l$ times.

\begin{definition}
Let us call the integers $A_{c_1,\dots,c_n}$
the {\it mixed Eulerian numbers}.
\end{definition}

The mixed Eulerian numbers are nonnegative integers because
hypersimplices are integer polytopes.  
In particular,
$n!\, \Vol P_{n+1}$ is a polynomial in $u_1,\dots,u_n$ with
positive integer coefficients.

\begin{example}
We have 
$$
\begin{array}{l}
\Vol P_2 = {\bf 1}\,u_1;\\
\Vol P_3 = {\bf 1}\, \frac{u_1^2}{2}  + {\bf 2}\, u_1 u_2 + 
{\bf 1} \,\frac{u_2^2}{2};\\
\Vol P_4 = {\bf 1}\,\frac{u_1^3}{3!} + {\bf 2}\,\frac{u_1^2}{2} u_2 +
{\bf 4}\,{u_1}\frac{u_2^2}{2} + {\bf 4}\, \frac{u_2^3}{3!} + 
{\bf 3}\,\frac{u_1^2}{2}u_3 + {\bf 6} \,{u_1u_2u_3} + \\
 \qquad \qquad \qquad 
 \qquad \qquad \qquad 
 \qquad 
+ {\bf 4}\, \frac{u_2^2}{2}u_3 + {\bf 3}\, u_1\frac{u_3^2}{2} + 
{\bf 2}\,{u_2}\frac{u_3^2}{2} + {\bf 1}\frac{u_3^3}{3!}.
\end{array}
$$
Here the mixed Eulerian numbers are marked in bold.
\end{example}

Recall that the usual {\it Eulerian number\/} $A(n,k)$ is 
defined as the number of permutations in $S_n$ with exactly $k-1$ descents.
It is well-known that $n!\,\Vol \Delta_{k,n+1}=A(n,k)$;
see Laplace~\cite[p.~257ff]{Lap}.

\begin{theorem} 
\label{th:mixed_eul_properties}
The mixed Eulerian numbers have the following
properties:

\begin{enumerate}
\item The numbers $A_{c_1,\dots,c_n}$ are positive integers defined
for $c_1,\dots,c_n\geq 0$, $c_1+\cdots+c_n=n$.

\item We have $A_{c_1,\dots,c_n} = A_{c_n,\dots,c_1}$.

\item For $1\leq k\leq n$, the number $A_{0^{k-1},n,0^{n-k}}$  
is the usual Eulerian number $A(n,k)$.
Here and below $0^l$ denotes the sequence of $l$ zeros.

\item We have $\sum 
\frac{1}{c_1!\cdots c_n!}\,A_{c_1,\dots,c_n} = (n+1)^{n-1}$,
where the sum is over $c_1,\dots,c_n\geq 0$ with $c_1+\cdots + c_n = n$.

\item We have $\sum A_{c_1,\dots,c_n} = n! \, C_n$,
where again the sum is over all $c_1,\dots,c_n\geq 0$ with $c_1+\cdots + c_n = n$
and $C_n= \frac{1}{n+1}\,\binom{2n}{n}$ is the Catalan number.

\item 
For $1\leq k\leq n$ and $i=0,\dots,n$,
the number $A_{0^{k-1},n-i,i,0^{n-k-1}}$ is equal to the number of permutations 
$w\in S_{n+1}$ with $k$ descents and $w(n+1)=i+1$.

\item We have $A_{1,\dots,1} = n!$.

\item We have $A_{k,0,\dots,0,n-k} = \binom{n}{k}$.

\item We have $A_{c_1,\dots,c_n} = 1^{c_1} 2^{c_2} \cdots n^{c_n}$
if $c_1+\cdots + c_i \geq i$, for $i=1,\dots,n-1$,
and $c_1+\cdots + c_n = n$.
There are exactly $C_n$ such sequences $(c_1,\dots,c_n)$.

\end{enumerate}
\end{theorem}

\begin{proof}
Properties (1) and (2) follow 
from the definition of the mixed Eulerian numbers.
Property (3) follows from the fact that $n!\,\Vol \Delta_{k,n+1}=A(n,k)$.
Property (4) follows from the fact that the volume of the regular permutohedron
$P_{n+1}(n,n-1,\dots,0)$, which corresponds to $u_1=\dots=u_n=1$,
equals $(n+1)^{n-1}$; see Proposition~\ref{prop:vol-Z-G}.
Property (5) follows from Theorem~\ref{th:equivalence-classes} below.
It was conjectured by R.~Stanley.
Property (6) is equivalent to the result by Ehrenborg, Readdy, and Steingr\'imsson
\cite[Theorem~1]{ERS} about mixed volumes of two adjacent hypersimplices.
Property (7) is a special case of Property~(9).

(8)  According to Theorem~\ref{th:vol=descent_number},
we have 
$$
\Vol P_{n+1}(x_1,0,\dots,0,x_{n+1}) = 
\sum_{k=0}^n (-1)^{n-k}\,D_{n+1}([k+1,n])\,
\frac {x_1^{k}}{k!}\, \frac {x_{n+1}^{n-k}}{(n-k)!},
$$
where $D_{n+1}([k+1,n])=\binom{n}{k}$ is the number of permutations
$w\in S_{n+1}$ such that $w_1<\cdots < w_{k+1}>w_{k+2}>\cdots >w_{n+1}$.
This permutohedron corresponds to
$u_1  = x_1$, $u_2 = \cdots u_{n-1}=0$, $u_n=-x_{n+1}$, which implies
that $A_{k,0,\dots,0,n-k} = \binom nk$.

(9) Let us use Theorem~\ref{th:second-formula}.
The $y$-variables are related to the $u$-variables as 
$$
\left\{
\begin{array}{l}
y_2 = u_1,\\
y_3 = u_2 - u_1,\\
y_4= u_3-2u_2 + u_1,\\
\vdots\\
\ds y_{n+1} = \sum_{i=0}^{n-1} (-1)^i \binom {n-1}{i}\, u_{n-i}
\end{array}
\right.
$$

Using these relations, we can express
any coefficient $[u_n^{c_1}\cdots u_1^{c_n}]\,V_{n+1}$ 
of the polynomial $V_{n+1} = \Vol P_{n+1}$ written in 
the $u$-coordinates as a combination
of coefficients $[y_{n+1}^{c_1'}\cdots y_2^{c_n'}]\,V_{n+1}$ 
of this polynomial written in the $y$-coordinates.
Let us assume that $(c_1,\dots,c_{n})$ satisfies
$c_1+\cdots + c_{i} \geq i$, for $i=1,\dots,n-1$, and $c_1+\cdots+ c_n = n$. 
Then any sequence $(c_1',\dots,c_n')$ that appears in this expression
satisfies the same conditions.
For such a sequence, we have
$$
[y_{n+1}^{c_1'}\cdots y_2^{c_n'}] \,V_{n+1}= 
\frac 1{c_1'!\cdots c_n'!} \, {\binom {n+1}{n+1}}^{c_1'}
{\binom{n+1}{n}}^{c_2'} \cdots {\binom{n+1}{2}}^{c_n'}.
$$
Indeed, any collection of subsets $J_1,\dots,J_n\subseteq[n+1]$
such that $c_i'$ of them have the cardinality $n+2-i$,
for $i=1,\dots,n$, automatically satisfies the dragon 
marriage condition; see Theorem~\ref{th:second-formula}.
Thus we have
$$
\begin{array}{l}
A_{c_1,\dots,c_n} = 
\left(\frac {\partial}{\partial u_n}\right)^{c_1}
\cdots
\left(\frac {\partial}{\partial u_1}\right)^{c_n} \,V_{n+1}
= \left(
\left(\frac {\partial}{\partial y_{n+1}}\right)^{c_1}
\left(
\frac {\partial}{\partial y_n}
- \binom{n-1}{1} \frac {\partial}{\partial y_{n+1}}
\right)^{c_2} \right.\times
\\[.1in]
\quad\times\left.
\left(
\frac {\partial}{\partial y_{n-1}} 
- \binom {n-2}{1} \frac {\partial}{\partial y_{n}}
+ \binom{n-1}{2} \frac {\partial}{\partial y_{n+1}}
\right)^{c_3}
\cdots \right)
\,V_{n+1} =\\[.1in]
\quad={\binom{n+1}{n+1}}^{c_1}
\left(\binom{n+1}{n} - \binom{n-1}{1}\binom{n+1}{n+1}\right)^{c_2}
\left(\binom{n+1}{n-1} - \binom{n-2}{1}\binom{n+1}{n} + 
\binom{n-1}{2}\binom{n+1}{n+1}\right)^{c_3}\cdots = \\[.1in]
\quad = 1^{c_1} 2^{c_2} \cdots n^{c_n}.
\end{array}
$$
In the last equality we used the binomial identity
$$
\sum_{i=0}^{k-1} (-1)^i\binom{n-k+i}{i} \binom{n+1}{n+2-k+i} = k,
\quad\textrm{for } 1\leq k\leq n,
$$
which we leave as an exercise.
\end{proof}

Let ``$\sim$'' be the equivalence relation of the set of
nonnegative integer sequences $(c_1,\dots,c_n)$ with $c_1+\cdots +c_n =n$
given by $(c_1,\dots,c_n)\sim (c_1',\dots,c_n')$
whenever $(c_1,\dots,c_n,0)$ is a cyclic shift of
$(c_1',\dots,c_n',0)$.

\begin{theorem}
\label{th:equivalence-classes}
For a fixed $(c_1,\dots,c_n)$, we have
$$
\sum_{(c_1',\dots,c_n')\sim (c_1,\dots,c_n)} A_{c_1',\dots,c_n'} = n!
$$
In other words, the sum of mixed Eulerian numbers in each
equivalence class is $n!$.

There are exactly the Catalan number $C_n=\frac{1}{n+1}\,\binom{2n}{n}$  
equivalence classes.
\end{theorem}

This claim was conjectured by R.~Stanley.
For example, it says that $A_{1,\dots,1} = n!$ and that
$A_{n,0,\dots,0} + A_{0,n,0,\dots,0}+A_{0,0,n,\dots,0}+\cdots 
+A_{0,\dots,0,n}= n!$, i.e., 
the sum of usual Eulerian numbers $\sum_k A(n,k)$ is $n!$.

\begin{remark}
The claim that there are $C_n$ equivalence classes 
is well-known.
Every equivalence class contains exactly one sequence $(c_1,\dots,c_n)$
such that $c_1+\cdots+c_i\geq i$, for $i=1,\dots,n$.
For this special sequence, the 
mixed Eulerian number is given by the simple product 
$A_{c_1,\dots,c_n} = 1^{c_1} \cdots n^{c_n}$; see
Theorem~\ref{th:mixed_eul_properties}.(9).
\end{remark}

Theorem~\ref{th:equivalence-classes}
follows from the following claim.

\begin{proposition} 
\label{prop:ciclic_symmetrization}
Let us write $\Vol P_{n+1}$ as a polynomial
$\hat V_{n+1}(u_1,\dots,u_{n+1})$ in $u_1,\dots, u_{n+1}$.
(This polynomial does not depend on $u_{n+1}$.)
Then the sum of cyclic shifts of this polynomial
equals
$$
\begin{array}{r}
\hat V_{n+1}(u_1,\dots,u_{n+1}) + \hat V_{n+1}
(u_{n+1},u_1,\dots,u_n) + \cdots +
\hat V_{n+1}(u_2,\dots,u_{n+1},u_1) =\quad \\ 
= (u_1+\cdots+u_{n+1})^n
\end{array}
$$
\end{proposition}

This claim has a simple geometric explanation
in terms of alcoves of the affine Weyl group.  
Cyclic shifts come from symmetries of the type $A_n$ extended 
Dynkin diagram.

\begin{proof}
Let $W=S_{n+1}$ be the type $A_n$ Weyl group.
The associated {\it affine Coxeter arrangement\/} 
is the hyperplane arrangement in the vector
space $\R^{n+1}/(1,\dots,1)\R \simeq \R^n$ given by 
$t_i -t_j =k$, for $1\leq i<j\leq n+1$ and $k\in \Z$.
Here and below in this proof the coordinates
$t_1,\dots,t_{n+1}$ in $\R^{n+1}$ are understood modulo
$(1,\dots,1)\R$.
These hyperplanes subdivide the vector space into simplices,
which are called the {\it alcoves.}  
The reflections with respect to these hyperplanes
generate the {\it affine Weyl group} $\Waff$ that acts 
simply transitively on the alcoves.

The {\it fundamental alcove\/} $A_\circ$ is given by the inequalities
$t_1>t_2>\cdots > t_{n+1}> t_1-1$.  It is the $n$-dimensional simplex
with the vertices
$v_0= (0,\dots,0)$, $v_1=(1,0,\dots,0)$, $v_2=(1,1,0,\dots,0)$, \dots,
$v_n=(1,\dots,1,0)$.
For $i=1,\dots,n$, the map 
$$
\phi_i:(t_1,\dots,t_{n+1})\mapsto 
(t_{i+1},\dots,t_{n+1},t_1-1,\dots,t_i-1)
$$ 
preserves the fundamental
alcove and sends the vertex $v_i$ to the origin $v_0$. 
We have $\Vol A_\circ = \frac{1}{|W|} = \frac{1}{(n+1)!}$,
assuming that we normalize the volume as in 
Section~\ref{sec:weight-polytopes}.

Let up pick a point
$x=(x_1,\dots,x_{n+1})$ in $A_\circ$.
The $\Waff$-orbit of $x$ has a unique representative in each alcove.
For any vertex $v$ of the affine Coxeter arrangement, i.e., for a
0-dimensional intersection of its hyperplanes, the convex hull of elements  
the orbit $\Waff \cdot x$ contained in the alcoves adjacent to $v$ is 
a (parallel translation) of a permutohedron.  
This collection of permutohedra 
associated with vertices of the arrangement forms a subdivision
of the linear space.

For the origin $v=v_0$, we obtain the permutohedron
$P_{(0)}=P_{n+1}(x_1,\dots,x_{n+1})$, and, for the vertex $v_i$, 
$i=1,\dots,n$, we obtain the permutohedron
$$
P_{(i)}=\phi_i^{-1}P_{n+1}(\phi_i(x)) =
\phi_i^{-1} P_{n+1}(x_{i+1},\dots,x_{n+1}, x_1-1,\dots,x_{i}-1).
$$
Note that, for $i=0,\dots,n$, we have
$\Vol P_{(i)} \cap A_\circ = \frac 1{|W|}\, \Vol P_{(i)}$.
Indeed, each permutohedron
$P_{(i)}$ is composed of $|W|$ isomorphic parts obtained by reflections
of $\Vol P_{(i)} \cap A_\circ$.

Thus the volume of the fundamental alcove times $|W|$ 
equals the sum of volumes of $n+1$ adjacent permutohedra,
For example, the 6 areas of the blue triangle on the following picture
is the sum of the areas of three hexagons.
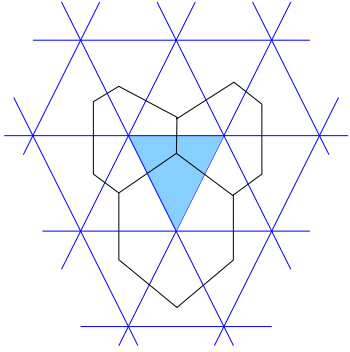
\begin{center}
\begin{picture}(0,0)%
\includegraphics{fig6-1.pstex}%
\end{picture}%
\setlength{\unitlength}{1579sp}%
\begingroup\makeatletter\ifx\SetFigFont\undefined%
\gdef\SetFigFont#1#2#3#4#5{%
  \reset@font\fontsize{#1}{#2pt}%
  \fontfamily{#3}\fontseries{#4}\fontshape{#5}%
  \selectfont}%
\fi\endgroup%
\begin{picture}(5424,5424)(2689,-5473)
\end{picture}

\end{center}
\nopagebreak
In other words,  we have
$1=|W|\cdot \Vol A_\circ = \sum_{i=0}^n\Vol P_{(i)}$.
The last expression can be written in the $u$-coordinates as
$$
\hat V_{n+1}(u_1,\dots,u_{n+1}) 
+ \hat V_{n+1}(u_2,\dots,u_{n+1},u_1) + \cdots
+ \hat V_{n+1} (u_{n+1},u_1,\dots,u_n),
$$
assuming that $u_1+\cdots + u_n = 1$.  The case of arbitrary
$u_1,\dots,u_n$ is obtained by multiplying all $u_i$'s by the same
factor $\alpha$ which corresponds to multiplying the volume by 
$\alpha^n$.
\end{proof}

\begin{proof}[Proof of Theorem~\ref{th:equivalence-classes}]
We obtain the required equality when we extract the coefficient
of $u_1^{c_1}\cdots u_n^{c_n} u_{n+1}^0$ in the both sides
of the identity in 
Proposition~\ref{prop:ciclic_symmetrization}.
\end{proof}

Proposition~\ref{prop:ciclic_symmetrization} together with Theorem~\ref{th:f1}
implies the following identity.  It would be interesting to find a direct proof
of this claim.

\begin{corollary}  
The symmetrization of the expression
$$
\frac{1}{n!} \,\frac{(\lambda_1 u_1 + (\lambda_1+\lambda_2)u_2 + \cdots
(\lambda_1+\cdots+\lambda_{n+1})u_{n+1})^n}{
(\lambda_1-\lambda_2)\cdots (\lambda_n-\lambda_{n+1})}
$$
with respect to $(n+1)!$ permutations of $\lambda_1,\dots,\lambda_{n+1}$
and $(n+1)$ {\it cyclic permutations} of $u_1,\dots u_{n+1}$ equals
$(u_1+\cdots+u_{n+1})^{n}$.
\end{corollary}

\section{Weighted binary trees}
\label{sec:weighted_binary_trees}

Let us give a combinatorial interpretation for 
the mixed Eulerian numbers based on plane binary trees.

Let $T$ be a plane binary tree on $[n]$ with 
the binary search labeling of the nodes;
see Subsection~\ref{ssec:associahedron}. 
There are the Catalan number $C_n$ of such trees.
For any node $i=1,\dots,n$, the set $\desc(i,T)$ of descendants of $i$ 
(including the node $i$ itself) is a  consecutive interval $\desc(i,T) = [l_i,r_i]$ 
of integers.  In particular, we have $l_i\leq i\leq r_i$.
For a pair nodes $i$ and $j$ in $T$ such that  $i\in \desc(j,T)$, i.e.,
$l_j\leq i\leq r_j$, define the weight
\begin{equation}
\label{eq:wij}
wt(i,j) = 
\min\left(
\frac{i - l_{j} +1}{j - l_{j} + 1},
\frac{r_{j} - i +1}{r_{j} - j + 1}\right) =
\left\{
\begin{array}{cl}
\frac{i - l_{j} +1}{j - l_{j} + 1} & \textrm{if } i \leq j,\\[.1in]
\frac{r - i +1}{r_{j} - j + 1} & \textrm{if } i > j.
\end{array}        
\right.
\end{equation}

Let $h(j,T):=|\desc(j,T)|$ be the ``hook-length'' of a node $j$ in a rooted tree $T$.

\begin{theorem}
\label{th:vol_perm_binary_trees}
The volume of the permutohedron $P_{n+1}$ is given by the following 
polynomial in the variables $u_1,\dots,u_n$:
$$
\Vol P_{n+1} = 
\sum_{T} 
\frac{n!}{\prod_{j=1}^n h(j,T)}\,
\prod_{j=1}^n \left(\sum_{i\in \desc(j,T)} wt(i,j) \, u_i\right),
$$
where the sum is over $C_n$ plane binary trees $T$ with $n$ nodes.
\end{theorem}

\begin{example} 
\label{ex:n=3_binary_trees}
For $n=3$, we have the following five binary trees, where we indicated the binary search 
labeling inside the nodes and also indicated the hook-lengths of the nodes:
\begin{center}
\input{fig8-1.pstex_t}
\end{center}
Theorem~\ref{th:vol_perm_binary_trees} says that
$$
\begin{array}{l}
\Vol P_4 = (u_1)(\frac{1}{2}\,u_1 + u_2) (\frac{1}{3}\,u_1+\frac{2}{3}\,u_2 + u_3)
+ (u_1 + \frac{1}{2}\,u_2) (u_2) 
(\frac{1}{3}\,u_1+\frac{2}{3}\,u_2+u_3)  \\[.1in]
\qquad\quad {}+ (u_1 + \frac{2}{3}\,u_2 + \frac{1}{3}\,u_3) (u_2)(\frac{1}{2}u_2 + u_3)
+ (u_1 + \frac{2}{3}\,u_2 + \frac{1}{3}\,u_3)
(u_2 + \frac{1}{2} u_3) (u_3)  \\[.1in]
\qquad\quad {}+ 2\cdot (u_1) (\frac{1}{2}\,u_1 + u_2 + \frac{1}{2}\,u_3)(u_3).
\end{array}
$$
\end{example}

\begin{corollary}
\label{cor:hook_lengths_binary_trees}  
We have
$$
(n+1)^{n-1} = \sum_T \frac{n!}{2^{n}} \, \prod_{j\in T} 
\left(1+\frac{1}{h(j,T)}\right),
$$
where is sum is over $C_n$ plane binary trees $T$ with $n$ nodes.
\end{corollary}

For $n=3$, the corollary says that
$(3+1)^2 = 3 + 3 + 3 +3 + 4$; see figure in Example~\ref{ex:n=3_binary_trees}.

\begin{proof} 
Let us specialize Theorem~\ref{th:vol_perm_binary_trees} for
$u_1=\dots=u_n=1$.  In this case, $P_{n+1}$ is the regular permutohedron
with volume $(n+1)^{n-1}$, see Proposition~\ref{prop:vol-Z-G}.
Easy calculation shows that $\sum_{i\in\desc(j,T)}wt(i,j) = \frac{h(j,T) +1}{2}$.
Thus the right-hand side of Theorem~\ref{th:vol_perm_binary_trees} gives the 
needed expression.
\end{proof}

Various combinatorial proofs and generalizations of
Corollary~\ref{cor:hook_lengths_binary_trees} were given by
Seo~\cite{Seo}, Du-Liu~\cite{DL}, and Chen-Yang~\cite{CY}.

An {\it increasing labeling\/} of nodes in a rooted tree $T$ on $[n]$
is a permutation $v\in S_n$ such that, whenever $i\in\desc(j,T)$,
i.e., the node $i$ is a descendant of the node $j$,
we have $v(i)\geq v(j)$.  
It is well-known that the number of increasing labelings 
is given by the following ``hook-length formula;''
see Knuth~\cite[Exer.~5.1.4.(20)]{Knu} and Stanley~\cite[Prop.~22.1]{St-OS}.
It can be easily proved by induction.

\begin{lemma}
\label{lem:num_incr_lab}
The number of increasing labeling of a tree $T$ 
equals $\frac{n!}{\prod_{j=1}^n h(j,T)}$.
\end{lemma}

Let us say that an {\it increasing binary tree\/} 
$(T,v)$ is a plane binary tree $T$ with the binary search labeling
as above and a choice of an increasing labeling $v$ of its nodes.
It is well-known that there are $n!$ increasing binary trees.
The map $(T,v)\mapsto v$ is a bijection between increasing binary trees
and permutations $v\in S_n$; cf.~\cite[1.3.13]{EC1}.

Let $\mathbf{i}=(i_1,\dots,i_n)\in [n]^n$ be a sequence of integers.
Let us say that an increasing binary tree $(T,v)$ 
is {\it $\mathbf{i}$-compatible\/}
if $i_{v(j)}\in [l_{j},r_{j}]$, for $j=1,\dots,n$.  
Define the {\it $\bf i$-weight} of an $\bf i$-compatible increasing binary
tree $(T,v)$ as
$$
wt(\mathbf{i},T,v) = \prod_{j=1}^n wt(i_{v(j)},j).
$$
where $wt(i_{v(j)},j)$ is given by~(\ref{eq:wij}).
The number $n!\,wt(\mathbf{i}, T, v)$ is always a positive integer.  The following lemma 
can be easily proved by induction, cf.\ Lemma~\ref{lem:num_incr_lab}.
We leave it as an exercise.

\begin{lemma}
We have, $n!$ divided by all denominators in $wt(\mathbf{i},T,v)$
equals the number labelings 
of the nodes of $T$ by permutations $w\in S_n$  such that,
for any node $j$, for which we pick the first (respectively, second) case in 
the definition of $wt(i_{v(j)}, j)$, the label $w(j)$ is less than labels $w(k)$
of all nodes $k$ in the left (respectively, right) branch of the node $j$.
\end{lemma}

\begin{example}
The following figure shows an $\bf i$-compatible increasing binary tree, for 
${\bf i} = (3,4,8,7,1,7,4,3)$.  The labels for the binary search labeling
are shown inside the nodes.  The increasing labeling is
$v = 5,2,8,7,1,3,6,4$ (shown in blue color).
The intervals $[l_j,r_j]$ are 
$[1,1]$, $[1,2]$, $[3,3]$, $[3,4]$, $[1,8]$, $[6,8]$, $[7,7]$, $[7,8]$.  
We also marked each node $j$  by the variable $u_{i_{v(j)}}$ (shown in red color).
The $\bf i$-weight of this tree is $wt({\bf i},T,v)=
\frac{3}{5}\cdot
\frac{1}{3}\cdot
\frac{1}{3}\cdot
\frac{1}{2}\cdot
\frac{1}{1}\cdot
\frac{1}{1}\cdot
\frac{2}{2}\cdot
\frac{1}{1}$.

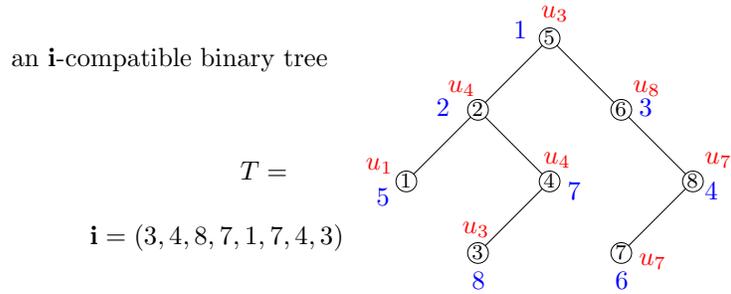
\begin{center}
\input{fig7-1.pstex_t}
\end{center}
\end{example}

Let us give a combinatorial interpretation for the mixed Eulerian numbers.

\begin{theorem}
\label{th:mixed_eulerian_binary_trees}
Let $(i_1,\dots,i_n)$ be any sequence such that
$u_{i_1}\cdots u_{i_n} = u_1^{c_1}\cdots u_n^{c_n}$.
Then
$$
A_{c_1,\dots,c_n} = \sum_{(T,v)}
n! \,wt(\mathbf{i}, T,v),
$$
where the sum is over $\bf i$-compatible increasing binary
trees $(T,v)$ with $n$ nodes.
\end{theorem}

Note that all terms $n!\,wt(\mathbf{i},T,v)$ 
in this formula are positive integers.
Actually, this theorem gives not just one but $\binom{n}{c_1,\dots, c_n}$
different combinatorial interpretations of the mixed Eulerian  numbers
$A_{c_1,\dots,c_n}$ for each way to write 
$u_1^{c_1}\cdots u_n^{c_n}$ as $u_{i_1}\cdots u_{i_n}$.
We will extend and prove Theorem~\ref{th:vol_perm_binary_trees}
in Section~\ref{sec:vol_weight_Phi}.
Let us now derive Theorem~\ref{th:mixed_eulerian_binary_trees} from it.

\begin{proof}[Proof of Theorem~\ref{th:mixed_eulerian_binary_trees}]
The volume of the permutohedron is obtained by multiplying 
the right-hand side of Theorem~\ref{th:vol_perm_binary_trees}
by  $\frac{1}{n!}\,u_{i_1}\cdots u_{i_n}$ and summing 
over all sequences ${\bf i} = (i_1,\dots,i_n)\in [n]^n$:
$$
\Vol P_{n+1} = \sum_{{\bf i}\in[n]^n} u_{i_1}\cdots u_{i_n}
\sum_{(T,v)} wt(\mathbf{i}, T,v),
$$
where the second sum is over $\bf i$-compatible increasing binary
trees $(T,v)$ with $n$ nodes.  This formula together with Lemma~\ref{lem:num_incr_lab}
implies the needed expression.
\end{proof}

\section{Volumes of weight polytopes via $\Phi$-trees}
\label{sec:vol_weight_Phi}

In this section we extend the results of the previous section
to weight polytopes for an arbitrary root system.
\medskip

Let $\Phi$ be an irreducible root system of rank $n$ with a choice of simple
roots $\alpha_1,\dots,\alpha_n$, and let $W$ be the associated Weyl group.
Let $(x,y)$ be a $W$-invariant inner product.
Let $\omega_1,\dots,\omega_n$ be the fundamental weights.
They form the dual basis to the basis of simple coroots 
$\alpha_i^\vee = \frac{2\alpha_i}{(\alpha_i,\alpha_i)}$.
Let $P_W(x)$ be the associated weight polytope, 
where $x = u_1\omega_1+\dots+u_n\omega_n$;
see Definition~\ref{def:weight_polytope}.
Its volume it a homogeneous polynomial $V_\Phi$ of degree $n$ 
in the variables $u_1,\dots,u_n$:
$$
V_\Phi(u_1,\dots,u_n) :=\Vol P_W(u_1\omega_1+\dots+u_n\omega_n).
$$

Recall the definition of $B(\Gamma)$-trees;
cf.\ Definition~\ref{def:B_forests} and Subsection~\ref{ssec:graph_ass}.

\begin{definition}
For a connected graph $\Gamma$,  a {\it $B(\Gamma)$-tree\/} is a rooted tree $T$ on the same 
vertex set such that
\begin{enumerate}
\item[(T1)] For any node $i$ and the set $I=\desc(i,T)$ of all descendants of $i$ in $T$, 
the induced graph $\Gamma|_I$ is connected.
\item[(T2)]  There are no two nodes $i\ne j$ such that the sets $I=\desc(i,T)$
and $J=\desc(j,T)$ are disjoint and the induced graph $\Gamma|_{I\cup J}$ is connected.
\end{enumerate}
\end{definition}

An {\it increasing $B(\Gamma)$-tree\/} $(T,v)$ is $B(\Gamma)$-tree $T$ together
with an increasing labeling $v$ of its nodes, defined as in 
Section~\ref{sec:weighted_binary_trees}.
In the case when $\Gamma$ is the Dynkin diagram of the root system $\Phi$,
we will call these objects {\it $\Phi$-trees\/} and {\it increasing $\Phi$-trees}.

The next proposition extends the well-known claim that there are
$n!$ increasing binary trees on $n$ nodes.

\begin{proposition}  For any connected graph $\Gamma$ on $n$ nodes,
the number of increasing $B(\Gamma)$-trees equals $n!$.
\end{proposition}

\begin{proof}
The map $(T,v)\mapsto v$ is a
bijection between increasing $B(\Gamma)$-trees and permutations $v\in S_n$.
\end{proof}

For a subset $I\subseteq [n]$, let $\Phi_I$ be the root system
with simple roots $\{\alpha_i\mid i\in I\}$, and let $W_I\subset W$ 
be the associated parabolic subgroup.
Let $\omega_i^I$, $i\in I$ be the fundamental weights for the 
root system $\Phi_I$.
For $j\in I\subseteq[n]$, let us define the linear form
$f_{I,j}(u) := \frac{1}{|I|} \sum_{i\in I} u_i\, (\omega_i^I,\omega_j^I)$
in the variables $u_i$.

\begin{theorem}
\label{th:Phi_volume_trees}
The volume of the weight polytope $P_W(x)$ is given by
$$
V_\Phi(u_1,\dots,u_n) = \frac{2^n\cdot |W|}
{\prod_{i=1}^n (\alpha_i,\alpha_i)}\,
\sum_{T} \prod_{j=1}^n f_{\desc(j,T),j}(u),
$$
where the sum is over all $\Phi$-trees $T$.
\end{theorem}

\begin{definition}
The {\it mixed $\Phi$-Eulerian numbers\/}
$A_{c_1,\dots,c_n}^\Phi$, for
$c_1,\dots,c_n\geq 0$, $c_1+\cdots + c_n =n$,
are defined as the coefficients of the polynomial expressing
the volume of the weight polytope:
$$
V_\Phi(u_1,\dots,u_n) = \sum_{c_1,\dots,c_n} A_{c_1,\dots,c_n}^\Phi\,
\frac{u_1^{c_1}}{c_1!}\cdots \frac{u_n^{c_n}}{c_n!}.
$$
Equivalently, 
the mixed $\Phi$-Eulerian numbers are the mixed volumes of
the {\it $\Phi$-hypersimplices,} which are the weight polytopes for the
fundamental weights.
\end{definition}

For a sequence $\mathbf{i}=(i_1,\dots,i_n)\in[n]^n$,
let us say that an increasing $\Phi$-tree $(T,v)$ is 
{\it $\mathbf{i}$-compatible\/}
if $i_{v(j)}\in \desc(j,T)$, for $j=1,\dots,n$.

\begin{theorem}
\label{th:Phi_mixed_eulerian_binary_trees}
Let $(i_1,\dots,i_n)$ be any sequence such that
$u_{i_1}\cdots u_{i_n} = u_1^{c_1}\cdots u_n^{c_n}$.
Then
$$
A_{c_1,\dots,c_n}^\Phi = 
\frac{2^n\cdot |W|}
{\prod_{i=1}^n (\alpha_i,\alpha_i)}\,
\sum_{(T,v)} \prod_{j=1}^n \left(\omega_{i_{v(j)}}^{\desc(j,T)},
\omega_{j}^{\desc(j,T)}\right),
$$
where the sum is over $\bf i$-compatible increasing 
$\Phi$-trees $(T,v)$.
\end{theorem}

The proof of these results is based on the following
recurrence relation for volumes of weight polytopes.
Let $\Phi_{(j)} := \Phi_{[n]\setminus \{j\}}$ be the root system
whose Dynkin diagram is obtained by removing the $j$th node,
and let $W_{(j)} := W_{[n]\setminus \{j\}}$ be the corresponding Weyl group,
for $j=1,\dots,n$.

\begin{proposition}
\label{prop:general_recur_weight_poly} 
For $i=1,\dots,n$, we have
$$
\frac{\partial }{\partial u_i} V_{\Phi}(u_1,\dots,u_n) = \sum_{j=1}^n
\frac{|W|}{|W_{(j)}|}\,
\frac {(\omega_i,\omega_j)} {(\alpha_j,\omega_j)}\,
V_{\Phi_{(j)}}(u_1,\dots,u_{j-1},u_{j+1},\dots,u_n).
$$
\end{proposition}

Note that $(\alpha_j,\omega_j)  = \frac{1}{2}\, (\alpha_j,\alpha_j)\,
(\alpha_j^\vee,\omega_j) = \frac{1}{2}\, (\alpha_j,\alpha_j)$.

\begin{proof}
The derivative $\partial V_{\Phi}/\partial u_i$ is the rate of change of the 
volume of the weight polytope as we move its generating vertex $x$
in the direction of the $i$th fundamental weight $\omega_i$.
It can be written as the sum of $(n-1)$-dimensional volumes of facets of 
$P_W(x)$ scaled by some factors,  which tell how fast the facets move.
Facets of $P_W(x)$ have the form $w( P_{W_{(j)}}(x))$, 
where $j\in [n]$ and $w\in W/W_{(j)}$.
In other words, $P_W(x)$ has $\frac{|W|}{|W_{(j)}|}$ 
facets isomorphic to $P_{W_{(j)}}(x)$.

The facet $P_{W_{(j)}}(x)$ 
is perpendicular to the fundamental weight $\omega_j$.
Note that this facet $P_{W_{(j)}}(x)$ is a parallel translate of
$P_{W_{(j)}}(x')$, where $x' =
 u_1 \omega_1^{(j)} + \cdots + u_{j-1}\omega_{j-1}^{(j)} +
 u_{j+1} \omega_{j+1}^{(k)} + \cdots + u_{n}\omega_{n}^{(j)}$
and $\omega_i^{(j)} := \omega_i^{[n]\setminus \{j\}}$.
Indeed, the fundamental weights $\omega_i^{(j)}$ for the root system 
$\Phi_{(j)}$ are projections of the fundamental weights $\omega_i$, $i\ne j$, 
for $\Phi$ to the hyperplane perpendicular to $\omega_j$.
Thus the $(n-1)$-dimensional volume of this facet is 
$\Vol P_{W_{(j)}}(x) = V_{\Phi_{(j)}}(u_1,\dots,u_{j-1},u_{j+1},\dots,u_n)$.

If we move $x$ in the direction of a vector $v$, then the facet 
$P_{W_{(j)}}(x)$
moves with the velocity proportional to $(v,\omega_j)$.
Recall that we normalize the volume so that the volume of the parallelepiped
generated by the simple roots $\alpha_1,\dots,\alpha_n$ is $1$;
see Section~\ref{sec:weight-polytopes}.  Thus the scaling factor for 
$v=\alpha_j$ is $1$, and, in general, the scaling factor 
is $\frac{(v,\omega_j)}{(\alpha_j,\omega_j)}$.
In particular, for $v=\omega_i$, we obtain the needed factor
$\frac{(\omega_i,\omega_j)}{(\alpha_j,\omega_j)}$.
By symmetry, all facets $w( P_{W_{(j)}}(x)$ come with the same
factors.
\end{proof}

\begin{proof}[Proof of Theorem~\ref{th:Phi_mixed_eulerian_binary_trees}]
Fix a sequence $\mathbf{i}=(i_1,\dots,i_n)$ such that
$u_{i_1}\cdots u_{i_n} = u_1^{c_1}\cdots u_n^{c_n}$.
Then, by the definition,
$$
A_{c_1,\dots,c_n}^\Phi = 
\frac{\partial}{\partial\,u_{i_n}}
\cdots 
\frac{\partial}{\partial\,u_{i_1}}\cdot V_{\Phi}(u_1,\dots,u_n).
$$
Applying Proposition~\ref{prop:general_recur_weight_poly} 
repeatedly, we deduce
that $A_{c_1,\dots,c_n}^\Phi$ equals the weighted sum over 
$\mathbf{i}$-compatible
increasing $\Phi$-trees $(T,v)$,
where each tree comes with the weight
$$
\prod_{k=1}^n \left(\frac{|W_{I_k}|}{\prod_l |W_{I_{k,l}}| }\cdot
\frac {2} {(\alpha_{j_k},\alpha_{j_k})}\,
(\omega_{i_k}^{I_k},\omega_{j_k}^{I_k})\right),
$$
where $j_1,\dots,j_n$ is the inverse permutation to $v$,
$I_k = \desc(j_k,T)$, and $I_{k,l}$, $l=1,2,\dots$, are the vertex sets 
of the branches of the vertex $j_k$ in $T$.
Note that all terms in the first quotient,
except the term $|W|$, cancel each other.
Thus we obtain the expression in the right-hand side of 
Theorem~\ref{th:Phi_mixed_eulerian_binary_trees}.
\end{proof}

\begin{proof}[Proof of Theorem~\ref{th:Phi_volume_trees}]
The volume $V_\Phi(u_1,\dots,u_n)$ is obtained by 
multiplying the right-hand side of 
Theorem~\ref{th:Phi_mixed_eulerian_binary_trees}
by $\frac{1}{n!}\, u_{i_1}\cdots u_{i_n}$ and summing 
over all sequences $(i_1,\dots,i_n)\in[n]^n$.
Thus we obtain
$$
V_\Phi(u_1,\dots,u_n) = 
\frac{2^n\cdot |W|}
{n!\cdot \prod_{i=1}^n (\alpha_i,\alpha_i)}\,
\sum_{T} \mathrm{incr}(T)\,\prod_{j=1}^n (|\desc(j,T)|
\cdot f_{\desc(j,T),j}(u)),
$$
where the sum is over all $\Phi$-trees $T$ and
$\mathrm{incr}(T)$ is the number of increasing labeling
of $T$.  Using Lemma~\ref{lem:num_incr_lab}, which says
that $\mathrm{incr}(T)=n!/\prod|\desc(j,T)|$, 
we derive the needed statement.
\end{proof}

For the Lie type $A_n$, Proposition~\ref{prop:general_recur_weight_poly} 
specializes to the following claim.
Let us write $\Vol P_{n+1}$ as a polynomial $V_{n+1}(u_1,\dots,u_n)$ in 
$u_1,\dots,u_n$.

\begin{proposition}
\label{prop:V_uuu_recurrence}
For any $i=1,\dots,n$, we have
$\frac{\partial }{\partial u_i} V_{n+1}(u_1,\dots,u_n) =$
$$\sum_{j=1}^n
\binom{n+1}{j}\,\frac{j\,(n+1-j)}{n+1}\, 
wt_{i,j,n}\,V_{j}(u_1,\dots,u_{j-1})\, V_{n-j+1}(u_{j+1},\dots,u_n),
$$
where
$wt_{i,j,n} = 
\min(\frac{i}{j},\frac{n+1-i}{n+1-j})$.
\end{proposition}

\begin{proof}
In this case, we have $W = S_{n+1}$, $V_W = V_{n+1}(u_1,\dots,u_n)$,
$W_{(j)} = S_{j}\times S_{n+1-j}$, $P_{W_j}=P_{j}\times P_{n+1-j}$,
and $V_{W_{(j)}} = V_{j}(u_1,\dots,u_{j-1})\,V_{n-j+1}(u_{j+1},\dots,u_n)$.
Thus $\frac{|W|}{|W_{(j)}|}= \binom{n+1}{j}$.
The root system lives in the space
$\{(t_1,\dots,t_{n+1})\in\R^{n+1}\mid t_1+\cdots t_{n+1}=0\}$
with the inner product induced from $\R^{n+1}$.
In this space, the simple roots are $\alpha_i = e_i-e_{i+1}$
and the fundamental weights are $\omega_i =
e_1+\cdots + e_i - \frac{i}{n+1}(1,\dots,1)$,
for $i=1,\dots,n$.
We have $(\alpha_j,\alpha_j)=2$ and $(\alpha_j,\omega_j)=1$.
Thus 
$\frac{(\omega_i,\omega_j)}{(\alpha_j,\omega_j)} 
= (\omega_i,\omega_j) = \min(i,j) - \frac{i\cdot j}{n+1} = \frac{j\,(n+1-j)}{n+1}\,wt_{i,j,n}$.
\end{proof}

\begin{proof}[Proof of Theorems~\ref{th:vol_perm_binary_trees}
and~\ref{th:mixed_eulerian_binary_trees}]
By Theorem~\ref{th:Phi_mixed_eulerian_binary_trees}
and proof of Proposition~\ref{prop:V_uuu_recurrence},
the mixed Eulerian number $A_{c_1,\dots,c_n}$ equals the weighted sum over 
$\mathbf{i}$-compatible increasing binary trees, where each tree 
$(T,v)$ comes with the weight
$$
(n+1)!\cdot \prod_{j=1}^n \frac{(j-l_{j}+1)\,(h_j+1-j)}{h_j+1}
\cdot
\min\left(\frac{i_{v(j)}-l_{j}+1}{j-l_{j}+1},
\frac{r_{j}-i_{v(j)} + 1}{r_{j} - j +1}\right),
$$
where $l_j\leq  r_j$ are defined as in 
Section~\ref{sec:weighted_binary_trees}
and $h_j = |\desc(j,T)|=r_{j}-l_{j}+1$.
All terms in the first quotient, except the term $\frac{1}{n+1}$,
cancel each other.  Note that the product
$\prod_{j=1}^n 
\min(\frac{i_{v(j)}-l_j+1}{j-l_j+1},\frac{r_j-i_{v(j)} + 1}{r_j - j +1})$
is exactly $wt(\mathbf{i},T,v)$.
Thus the total weight of $(T,v)$ equals $(n+1)!\,\frac{1}{n+1}\,
wt(\mathbf{i},T,v)$, as needed.
\end{proof}

\section{Appendix: Lattice points and Euler-MacLaurin formula}
\label{sec:appendix-lattice-points}

In this section, we review some results of Brion~\cite{Bri},
Khovanskii-Pukhlikov~\cite{KP1, KP2}, Guillemin~\cite{Gui},
and Brion-Vergne~\cite{BV1, BV2}
related to counting lattice points and volumes of polytopes.
For the completeness sake, we included short proofs of these results.

Instead of calculating the volume or counting the number of lattice points 
in a polytope, let us sum monomials over the lattice points in the polytope.  
We can work with unbounded polyhedra, as well.

Recall that a polytope in $\R^n$ is a convex hull of a finite set
of vertices.  A {\it rational polyhedron} in $\R^n$ is an intersection 
of a finite set of half-spaces with rational (equivalently, integer) 
coordinates. 
In particular, rational polyhedra include polytopes with rational
vertices and {\it rational cones}, i.e., cones with a rational vertex
and integer generating vectors.

Let $\chi_P:\Z^n\to\Q$ be the {\it characteristic function} 
(restricted to the integer lattice) of a polyhedron $P$ given by
$\chi_P(x) = 1$, if $x\in P$, and $\chi_P(x)=0$, if $x\not\in P$.
The {\it algebra of rational polyhedra} $A$ is the linear space 
of functions $\Z^n\to\R$ spanned by the characteristic functions 
$\chi_P$ of rational polyhedra.  The space $A$ is closed
under multiplications of functions,
because $\chi_P\cdot \chi_Q = \chi_{P\cap Q}$.
The algebra $A$ is generated by the 
{\it Heaviside functions} $H_{h,c}=\chi_{\{x\mid h(x)\geq c\}}$, 
where $h$ is an integer linear form 
and $c\in\Z$.

The group algebra of the integer lattice $\Z^n$ is 
the algebra of Laurent polynomials $\Q[t_1^{\pm 1},\dots,t_n^{\pm1}]$. 
Let $\Q(t_1,\dots,t_n)$ be the {\it field of rational functions},
which is the field of fractions of the group algebra.
For a vector $a\in\Z^n$, let $t^a := t_1^{a_1}\cdots t_n^{a_n}$.

\begin{theorem}
\label{th:map-S}
{\rm Khovanskii-Pukhlikov~\cite{KP1}} \
There exists a unique linear map $S:A\to \Q(t_1,\dots,t_n)$ such that 
\begin{enumerate}
\item[(a)] $S(\delta) = 1$, where $\delta=\chi_{\{0\}}$ 
is the delta-function.
\item[(b)] For any $\nu \in A$ and $a\in\Z^n$, we have $S(\nu(x-a)) = 
t^a\,S(\nu)$.
\end{enumerate}
The map $S$ has the following properties:
\begin{enumerate}
\item For a function $\nu$ on $\Z^n$ with a finite support, we have 
$S(\nu)  = \sum_{a} \nu(a)\,t^a$.   In particular, for a polytope $P$, we have
$S(\chi_P) = \sum_{a\in P\cap\Z^n} t^a$.
\item If $\nu\in A$ is a $b$-periodic function for some
nonzero vector $b\in\Z^n$, i.e., $\nu(x)\equiv \nu(x-b)$, then $S(\nu)=0$.
Thus, for a rational polyhedron $P$ that contains a line,
we have $S(\chi_P)=0$.
\item For a simple rational cone 
$C=v+\R_{\geq 0} g_1 + \cdots + \R_{\geq 0} g_m$, where
$v\in\Q^n$ and $g_1,\dots,g_m\in\Z^n$ are linearly independent, we have
$$
S(\chi_C) = \left(\sum_{a\in\Pi\cap\Z^n} t^a \right) \prod_{i=1}^m (1-t^{g_i})^{-1},
$$
where $\Pi$ is the parallelepiped
$\{v + c_1 g_1 + \cdots + c_m g_m \mid 0\leq c_i < 1\}$.
\end{enumerate}
\end{theorem}

\begin{proof}
Let us first check that conditions (a) and (b) imply properties (1), (2), and (3).
We have $S(\nu) = S(\sum_a \nu(a) \delta(x-a))=\sum_{a} \nu(a)\,t^a$, for 
a function $\nu$ with a finite support.
For a $b$-periodic function $\nu\in A$, we have $S(\nu) = t^b S(\nu)$ by (b), 
and, thus, $S(\nu)=0$.
Let us write, using the inclusion-exclusion principle,
$\chi_\Pi = \chi_C - \sum_i \chi_{C + g_i} + \sum_{i<j} \chi_{C+g_i+g_j} -
\cdots$.  Thus by (b), 
we have $S(\chi_\Pi) = S(\chi_C) - (\sum_i t^{v_i}) S(\chi_C)
+(\sum_{i<j} t^{g_i+g_j})S(\chi_C) -\cdots = S(\chi_C) \prod_i (1-t^{g_i})$,
which is equivalent to (3).

Let us now prove the existence and uniqueness of the map $S$.  
We can subdivide any rational polyhedron $P$ into rational simplices and simple
rational cones.  Furthermore, we can present the characteristic function of a
simplex as an alternating sum of characteristic functions of simple rational
cones.  Thus we can write $\chi_P$ as a linear combination of characteristic
functions of simple rational cones.  Since conditions (a) and (b) imply
expression (3) for $S(\chi_C)$ for each simple rational cone, the
expression $S(\chi_P)$ is uniquely determined by linearity.

Let us verify that this construction for $S$ is consistent.
In other words, we need to check that, for any linear dependence
$b_1\chi_{C_1} + \cdots + b_N\chi_{C_N} = 0$
of characteristic functions of simple rational cones, we have
$b_1 S(\chi_{C_1}) + \cdots + b_N S(\chi_{C_N})=0$, where each term 
$S(\chi_{C_i})=f_i\cdot \prod_j (1-t^{v_{ij}})^{-1}$ is given 
by expression (3).  Here $f_i$ are certain Laurent polynomials.
Let us assume that $b_1\chi_{C_1} + \cdots + b_N\chi_{C_N} = 0$
and $b_1 S(\chi_{C_1}) + \cdots + b_N S(\chi_{C_N})=f/D$, 
where $f$ is a nonzero Laurent polynomial and $D=\prod_{ij}(1-t^{v_{ij}})$ 
is the common denominator of the terms $S(\chi_{C_i})$.  Let us select a norm 
on $\Z^n$, for example, 
$|a|:=\sqrt{a_1^2+\cdots+a_n^{2}}$.  Let $R$ be a sufficiently large number
such that $R> |a|$ for any monomial $t^a$ that occurs in $f$ or $D$ with
a nonzero coefficient.  We can write each term as
$S(\chi_{C_i}) = \sum_{|a|\leq 3R} \chi_{C_i}(a)\,t^a + 
\tilde f_i \cdot \prod_j (1-t^{v_{ij}})^{-1}$,
where, for any monomial $t^a$ that occurs in $\tilde f_i$, we have $|a|> 2R$.
Let us  sum the right-hand sides of these expressions with 
the coefficients $b_i$. 
Then the first terms cancel and we obtain 
$b_1 S(\chi_{C_1}) + \cdots + b_N S(\chi_{C_N}) = \sum_i 
\tilde f_i \prod_j (1-t^{v_{ij}})^{-1} = f/D$.  We deduce that $f$ is a 
linear combination of monomials $t^a$ with $|a|>R$, which contradicts to 
our choice of $R$.  This proves the existence and 
uniqueness of the map $S$.
\end{proof}

Let $A'$ be the subspace in the algebra of rational polyhedra $A$ 
spanned by characteristic functions $\chi_P$ of
rational polyhedra $P$ that contain lines.  
According to Theorem~\ref{th:map-S}, we have $S(f)=0$, for any $f\in A'$.
Thus we obtain a well-defined linear map $S:A/A'\to B$.

For a rational polyhedron $P$ and a point $u\in\Z^n$, 
let $C_{P,u}$ denote the rational cone with the vertex at $u$ 
such that $P\cap B = C_{P,u}\cap B$ for a sufficiently small
open neighborhood $B$ of $u$.  Notice that $\chi_{C_{P,u}} \not \in A'$ 
if and only if $u$ is a vertex of the polyhedron $P$.

For an analytic function $f(t)$ defined in a neighborhood of 0, let $[t^n]\,f(t)$ 
denote the coefficient of $t^n$ in its Taylor expansion.
Notice that $\frac{t}{1-e^{-t}} = 1 + \frac{t}{2} + 
\sum_{k=1}^\infty (-1)^{k-1}\,B_k\,\frac{t^{2k}}{(2k)!}$,
is an analytic function at $t=0$, where $B_k$ are the Bernoulli numbers.

\begin{theorem}
\label{th:S-any-polytope} 
{\rm Brion~\cite{Bri}, Khovanskii-Pukhlikov~\cite{KP1}} \

{\rm (1)}  For any rational polyhedron $P$, we have
$\chi_P \equiv \sum_{v\in V} \chi_{C_{P,v}}$ modulo the subspace $A'$,
where the sum is over the vertex set $V$ of $P$.

{\rm (2)}  We have $S(P) = \sum_{v\in V} S(C_{P,v})$.
In particular, for a simple rational polyhedron $P$, we
have
$$
S(P) = \sum_{v\in V} \frac {\sum_{a\in\Pi_v\cap\Z^n} z^a}
{\prod_{i=1}^n(1-z^{g_{i,v}})},
$$
where the sum is over vertices $v$ of $P$, $g_{1,v},\dots,g_{n,v}\in\Z^n$
are the integer generators of the cone $C_{P,v}$,
and $\Pi_v = \{v+c_1 g_{1,v} + \cdots + c_n g_{n,v} \mid 0\leq c_i < 1\}$.

{\rm (3)}  For a simple rational polytope $P$, 
the number of lattice points in $P$ equals
$$
\#\{P\cap \Z^n\} = 
[t^n] \left\{\sum_{v\in V} \left(\sum_{a\in\Pi_v\cap\Z^n} e^{t\cdot h(a)}\right)
\prod_{i=1}^n \frac {t} { 1 - e^{t\cdot h(g_{i,v})}}\right\}.
$$
where $h\in (\R^n)^*$ is any linear form such that $h(g_{i,v})\ne 0$,
for all vectors $g_{i,v}$.

{\rm (4)}  The volume of a simple rational polytope $P$ equals
$$
\Vol P = \frac{1} {n!} \sum_{v\in V} \frac{|\det(g_{1,v},\dots,g_{n,v})| \, h(v)^n}
{(-1)^n \prod_{i=1}^n h(g_{i,v})},
$$
where $\det(g_{1,v},\dots,g_{n,v})$ is the determinant of the $n\times n$-matrix 
with the row vectors $g_{i,v}$ and $h\in (\R^n)^*$ is any linear form such that 
$h(g_{i,v})\ne 0$, for all vectors $g_{i,v}$.
\end{theorem}

The formula for the sum of exponents $S(P)$ was first obtained by M.~Brion~\cite{Bri}. 
The formula for $\Vol P$ was given by Khovanskii-Pukhlikov~\cite{KP2}
(in case of Delzant polytopes) and by Brion-Vergne~\cite{BV1} in general.

\begin{proof}
(1) As we have mentioned in the proof of Theorem~\ref{th:map-S},
we can write the characteristic function of a rational polyhedron as
a finite linear combination of characteristic functions of rational cones:
$\chi_P  = \sum_i b_i\, \chi_{C_i}$.  Let $U\supseteq V$ be the set of vertices
of all cones $C_i$.  For $u\in U$, let $I_u$ be the collection of indices
$i$ such that the cone $C_i$ has the vertex $u$. 
Then  $\sum_{i\in I_u} b_i \,\chi_{C_i}\equiv \chi_{C_{P,u}}\pmod{A'}$.
Also $\chi_{C_{P,u}}\in A'$, for $u\in U\setminus V$.
This proves the claim.  

(2) This claim follows from (1) and Theorem~\ref{th:map-S}.

(3) Let us pick a linear form $h$ that does not annihilate any of the vectors
$g_{i,v}$.  Let $B$ be the subalgebra of $\Q(t_1,\dots,t_n)$
generated by the $z^a$ and $\frac{1}{1-z^b}$, for $a,b\in\Z^n$ such that $h(b)\ne 0$. 
Let $e_h:B\to \R((q))$ be the homomorphism from $B$ to the ring of formal Laurent
series in one variable $q$ given by $z^a\mapsto e^{q\cdot h(a)}$
and $\frac{1}{1-z^b}\mapsto \frac{1}{1-e^{q\cdot h(b)}}$. 
Let us apply the homomorphism $e_h$ to the expression for $S(P)$ given by (2).
Then the number of lattice points $\#\{P\cap \Z^n\}$ is the constant
coefficient of the resulting Laurent series.  This is exactly the need claim.

(4) The volume of a polytope $P$ can be calculated by counting the number
of lattice points in the inflated polytope $kP$ for large $k$.
Explicitly,  $\Vol P = \lim_{k\to\infty}\#\{kP\cap \Z^n\}/k^n$.
The vertices of the inflated polytope $kP$ are the vectors $k\,v$, for $v\in V$,
and the generators of the cone $C_{kP,kv}$ are exactly the same vectors
$g_{i,v}$ as for the original polytope $P$.  We may assume that the limit
is taken over $k$'s such that all vectors $k\,v$ are integer.
Each term in the expression for $\#\{kP\cap \Z^n\}$ given by (3)
has the form
$[t^n] \left\{ e^{t\cdot h(kv + a')}
\prod_{i=1}^n \frac {t} { 1 - e^{t\cdot h(g_{i,v})}}\right\}
= [t^n] \left\{ e^{t\cdot h(kv + a')}
\prod_{i=1}^n ( - \frac{1}{h(g_{i,v})} + O(t))\right\}$,
where $a'\in(\Pi_v- v)\cap\Z^n$.
Since $k$ appears only in the first exponent, this expression 
is a polynomial in $k$ of degree $n$ with the top term 
$k^n\, \left(\frac{1}{n!} h(v)^n (-1)^n \prod_{i=1}^n \frac{1}{h(g_{i,v})}\right)$.
There are 
$|\det(g_{1,v},\dots,g_{n,v})| = | \Pi_v|$ choices for $a'$.
Thus summing these expressions over all $v$ and $a'$ we obtain the needed 
expression for $\Vol P$.
\end{proof}

For a polytope $P$ with the vertices $v_1,\dots,v_M$, we say 
that a {\it deformation} of $P$ is a polytope of the form 
$P'=\mathrm{ConvexHull}(v_1',\dots,v_M')\in\R^n$ such that $v_i'-v_j' = k_{ij} (v_i-v_j)$,
for some nonnegative $k_{ij}\in\R_{\geq 0}$, whenever $[v_i,v_j]$ is a 1-dimensional
edge of $P$.  
A generic deformation of $P$ has the same combinatorial structure as $P$.
However in degenerate cases some of the vertices $v_i'$ may merge with each other.

Deformations of $P$ are obtained by parallel translations of its facets.
Suppose that the polytope $P$ has $N$ facets and is given by the linear inequalities
$P=\{x\in\R^n\mid h_i(x)\leq c_i, \ i=1,\dots,N\}$, for 
some $h_i\in(\R^n)^*$ and $c_i\in\R$.  Then any deformation $P'=
\mathrm{ConvexHull}(v_1',\dots,v_M')$ has the 
form
$$
P(z_1,\dots,z_n) := \{x\in\R^n\mid h_i(x)\leq z_i, \ i=1,\dots,N\},
\textrm{ for some } z_1,\dots,z_N\in \R,
$$
where $h_i(v_j')=z_i$ whenever $i$-th facet of $P$ contains
the $j$-th vertex $v_j$.
For this polytope we will write $v_i(z_1,\dots,z_N) = v'_i$.
Let $\D_P\subset\R^N$ be the set of $N$-tuples $(z_1,\dots,z_N)$
corresponding to deformations of $P$.  Then $\D_P$ is a certain
polyhedral cone in $\R^N$ that we call the {\it deformation cone}.
If $P$ is a simple polytope then $\D_P$ has dimension $N$, 
because any sufficiently small parallel translations of the facets of $P$ 
give a deformation of $P$.
Deformations $P(z_1,\dots,z_n)$ for interior points $(z_1,\dots,z_n)
\in \D_P\setminus \partial\D_P$ of the cone $\D_P$ are exactly 
the polytopes whose associated fan coincides with the fan of $P$.

A simple integer polytope $P$ is called a {\it Delzant polytope} if,
for each vertex $v$ of $P$, the cone $C_{P,v}$ is generated by 
an integer basis of the lattice $\Z^n$.
Such polytopes are associated with smooth toric varieties.
Formulas in Theorem~\ref{th:S-any-polytope} are especially simple 
for Delzant polytopes.  Indeed, in this case $\Pi_v\cap\Z^n$ consists
of a single element $v$ and 
$|\det(g_{1,v},\dots,g_{n,v})| = 1$. 
For Delzant polytopes, we assume that we pick the linear forms $h_i$
corresponding to the facets of $P$ so that $h_i$ are integer and 
are not divisible by a nontrivial integer factor.

Let $I_P(z_1,\dots,z_N)=\#\{P(z_1,\dots,z_N)\cap\Z^n\}$ 
be the number of lattice points
and $V_P(z_1,\dots,z_N)= \Vol P(z_1,\dots,z_N)$ be the volume 
of a deformation of $P$.

Let $\Todd(q) = \frac{q}{1-e^{-q}}$.  Since $\Todd(q)$ expands
as a Taylor series at $q=0$, we have the well-defined
operators $\Todd\left(\frac{\partial}{\partial\, z_i}\right)$
acting on polynomials in $z_1,\dots,z_N$.

\begin{theorem}
\label{th:Todd-Euler-Maclaurin}
{\rm (1)}
For an integer polytope $P$, 
and $(z_1,\dots,z_N)\in \D_P\cap\Z^N$,
the number of lattice points $I_P(z_1,\dots,z_N)$
and the volume $V_P(z_1,\dots,z_N)$ are given by polynomials in $z_1,\dots,z_N$
of degree $n$.
The polynomial $V_P(z_1,\dots,z_N)$ is the top homogeneous component
of the polynomial $I_P(z_1,\dots,z_N)$.

{\rm (2)}
If $P$ is a Delzant polytope then we have
$$
I_P(z_1,\dots,z_N) = \left(\prod_{i=1}^N 
\Todd\left(\frac{\partial}{\partial\, z_i}\right) \right)\, V_P(z_1,\dots,z_N).
$$
\end{theorem}

We will call the polynomial $I_P(z_1,\dots,z_N)$ the {\it generalized 
Ehrhart polynomial} of the polytope $P$.

\begin{proof}
(1) 
Assume $P$ is a simple polytope.  The vertices $v_i(z_1,\dots,z_N)$ 
of the deformation $P(z_1,\dots,z_N)$ linearly depend on $z_1,\dots,z_N$.
According to formulas (3) and (4) in Theorem~\ref{th:S-any-polytope},
$I_P(z_1,\dots,z_N)$ and $V_P(z_1,\dots,z_N)$ are polynomials in 
$z_1,\dots,z_N$, because each term in these formulas for $P(z_1,\dots,z_N)$
polynomially depend on $v$.   This remains true for 
degenerate deformations $P(z_1,\dots,z_N)$ when some of the vertices $v_i(z_1,\dots,z_N)$
merge.  Indeed, all claims of Theorem~\ref{th:S-any-polytope} remain valid
(and proofs are exactly the same) if, instead of summation over actual
vertices of $P(z_1,\dots,z_N)$, we sum over $v_i(z_1,\dots,z_N)$.
If $P$ is not simple then a generic small parallel translation of its facets 
results in a simple polytope.  Thus $P$ can be thought of as a degenerate
deformation of a simple polytope and the above argument works.

(2)  For a simple polytope $P$, we have
$$
\frac{\partial }{\partial z_i} v_j(z_1,\dots,z_N) =
\left\{
\begin{array}{cl}
-\alpha_{ij} \,g_{k,v_j} & \textrm{if $v_j$ belongs to the $i$-th facet},\\
0  & \textrm{otherwise},
\end{array}
\right.
$$
for some positive constants $\alpha_{ij}$,
where $g_{k,v_j}$ is the only  generator of the cone $C_{P,v_j}$ that is not 
contained in the $i$-th facet.  Indeed, a small parallel translation
of the $i$-th facet, moves each vertex $v_j$ in this facet in the direction opposite 
to the generator $g_{k,v_j}$ and does not change all other vertices.
If $P$ is a Delzant polytope then all constants $\alpha_{ij}$ are equal to $1$.
In this case,
by Theorem~\ref{th:S-any-polytope}(4), we have 
$$
V_P(z_1,\dots,z_N) = \frac{1}{n!}\sum_{j=1}^M 
\frac{h(v_j(z_1,\dots,z_N))^n}{(-1)^n\prod_{i=1}^n h(g_{i,v_j})}
= [t^n]\left\{\sum_{j=1}^M   \frac{e^{t\cdot h(v_j(z_1,\dots,z_N))}}
{ (-1)^n\prod_{i=1}^n h(g_{i,v_j})}
 \right\}
$$
The only term in this expression that involves $z_i$'s is the exponent
$e^{t\cdot h(v_j(z_1,\dots,z_N))}$.
For an analytic function $f(q)$, 
the operator $f\left(\frac{\partial }{\partial z_i}\right)$
maps this exponent to
$$
e^{t\cdot h(v_j(z_1,\dots,z_N))}
\mapsto
\left\{
\begin{array}{cl}
e^{t\cdot h(v_j(z_1,\dots,z_N))}\,f(-t\,h(g_{k,v_j}))
& \textrm{if $v_j$ lies in the $i$-th facet},\\
e^{t\cdot h(v_j(z_1,\dots,z_N))}\,f(0) & \textrm{otherwise},
\end{array}
\right.
$$
where $k$ is the same as above.
Using this for Todd operators, we obtain the expression for
$I_P(z_1,\dots,z_N)$ given by Theorem~\ref{th:S-any-polytope}(3).
\end{proof}

\end{document}